\input amstex
\documentstyle{amsppt}
\magnification=\magstep1
\hsize=5.2in
\vsize=6.7in

\catcode`\@=11
\loadmathfont{rsfs}
\def\mycal{\mathfont@\rsfs}
\csname rsfs \endcsname

\topmatter
\title  COCYCLE AND ORBIT EQUIVALENCE SUPERRIGIDITY \\
FOR MALLEABLE ACTIONS OF $w$-RIGID GROUPS
\endtitle
\author SORIN POPA \endauthor

\rightheadtext{Cocycle superrigidity}

\affil University of California, Los Angeles\endaffil

\address Math.Dept., UCLA, LA, CA 90095-155505\endaddress
\email popa@math.ucla.edu\endemail

\thanks Supported in part by NSF Grant 0601082.\endthanks

\abstract We prove that if a countable discrete group $\Gamma$ is
{\it w-rigid}, i.e. it contains an infinite normal subgroup $H$ with
the relative property (T) (e.g. $\Gamma= SL(2,\Bbb Z) \ltimes \Bbb
Z^2$, or $\Gamma = H \times H'$ with $H$ an infinite Kazhdan group
and $H'$ arbitrary), and $\Cal V$ is a closed subgroup of the group
of unitaries of a finite separable von Neumann algebra (e.g. $\Cal
V$ countable discrete, or separable compact), then any $\Cal
V$-valued measurable cocycle for a measure preserving action $\Gamma
\curvearrowright X$ of $\Gamma$ on a probability space $(X,\mu)$
which is weak mixing on $H$ and {\it s-malleable} (e.g. the
Bernoulli action $\Gamma \curvearrowright [0,1]^\Gamma$) is
cohomologous to a group morphism of $\Gamma$ into $\Cal V$. We use
the case $\Cal V$ discrete of this result to prove that if in
addition $\Gamma$ has no non-trivial finite normal subgroups then
any orbit equivalence between $\Gamma \curvearrowright X$ and a free
ergodic measure preserving action of a countable group $\Lambda$ is
implemented by a conjugacy of the actions, with respect to some
group isomorphism $\Gamma \simeq \Lambda$.
\endabstract

\endtopmatter

\document

\heading 0. Introduction
\endheading

There has recently been increasing interest in the study of measure
preserving actions of groups on probability measure spaces up to
{\it orbit equivalence} ({\it OE}), i.e. up to isomorphisms of
probability spaces taking the orbits of one action onto the orbits
of the other. While the early years of this subject concentrated on
the amenable case, culminating with the striking result that all
ergodic m.p. actions of all countable amenable groups are
undistinguishable under orbit equivalence ([Dy], [OW], [CFW]), the
focus is now on proving {\it OE rigidity} results showing that, for
special classes of (non-amenable) group actions, OE of $\Gamma
\curvearrowright X, \Lambda \curvearrowright Y$ is sufficient to
insure isomorphism of groups $\Gamma \simeq \Lambda$, or even
conjugacy of the actions. Ideally, one seeks to prove this under
certain conditions on the ``source'' group-action $\Gamma
\curvearrowright X$ but no condition at all (or very little) on the
``target'' action $\Lambda \curvearrowright Y$, a type of result
labeled {\it OE superrigidity}.

The first OE rigidity phenomena were discovered by Zimmer, who used
his celebrated cocycle superrigidity theorem to prove that any free
ergodic m.p. actions of the groups SL$(n,\Bbb Z)$, $n=2,3,...$, are
orbit inequivalent for different $n$'s ([Z1,2]). A parallel
discovery in von Neumann algebra theory, due to Connes ([C3]),
showed that II$_1$ factors from groups with property (T) of Kazhdan
have rigid symmetry structure. Ideas from ([C3], [Z2]) were used to
produce more OE rigidity results in ([Ge], [GeGo], [P8], etc), while
rigidity aspects of free and Gromov's word-hyperbolic groups began
to emerge in ([A1,2]). An important development came with the
remarkable work of Gaboriau on cost and $\ell^2$-cohomology
invariants for OE relations, showing in particular that free m.p.
actions of the free groups $\Bbb F_n$ are OE inequivalent, for
different $n\geq 1$ ([G1,2]). At the same time OE superrigidity
phenomena started to unveil in the work of Furman who proved the
striking result that, more than just being rigid, actions of higher
rank lattices such as SL$(n,\Bbb Z) \curvearrowright \Bbb T^n$, for
$n\geq 3$ odd, are in fact {\it OE superrigid}, i.e. any orbit
equivalence between such an action and an arbitrary free m.p. action
of a discrete group $\Lambda$ comes from a conjugacy ([Fu1,2]). A
startling new set of OE rigidity results was then established by
Monod and Shalom, for doubly ergodic actions of products of
word-hyperbolic groups ([MoS1,2]). The latest in this line is a
result in ([P3]), showing that if $\Gamma$ is an infinite conjugacy
class (ICC) Kazhdan group then any {\it s-malleable} mixing action
$\Gamma \curvearrowright (X,\mu)$ (such as the Bernoulli action
$\Gamma \curvearrowright [0,1]^\Gamma$) is OE superrigid. Unlike
previous OE rigidity results, which are all obtained in a measure
theoretic framework (albeit each using different techniques), in
([P3]) this is a consequence of a ``purely'' von Neumann algebra
rigidity result, in which the conjugacy class of $\Gamma
\curvearrowright X$ is recovered from the {\it group measure space}
von Neumann algebra $L^\infty X \rtimes \Gamma$, which apriori
contains less information than the OE class of the action (cf [CJ]).

We prove in this paper an OE superrigidity result covering a much
larger family of s-malleable actions $\Gamma \curvearrowright X$
than in ([P3]), with $\Gamma$ merely required to contain an infinite
normal subgroup with the relative property (T) (i.e. $\Gamma$ is
{\it w-rigid}) and the action assumed weak (rather than strong)
mixing. We in fact prove a stronger form of OE superrigidity for the
actions $\Gamma \curvearrowright X$, showing that if $\Lambda
\curvearrowright Y$ is an arbitrary free m.p. action and $\Delta:X
\simeq Y$ is an isomorphism of probability spaces taking each
$\Gamma$-orbit into a $\Lambda$-orbit (not necessarily onto), then
there exists a subgroup $\Lambda_0\subset \Lambda$ such that the
$\Gamma$ and $\Lambda_0$ actions are conjugate, i.e. a suitable
perturbation $\Delta_0$ of $\Delta$ by an automorphism in the full
group of $\Lambda$ satisfies $\Delta_0 \Gamma \Delta_0^{-1} =
\Lambda_0.$

The main result of this paper though is a {\it cocycle
superrigidity} result for the actions $\Gamma \curvearrowright X$,
from which the OE superrigidity is just a consequence. Thus, we show
that any measurable cocycle for $\Gamma \curvearrowright X$ with
values in an arbitrary discrete group $\Lambda$ is equivalent to a
group morphism of $\Gamma$ into $\Lambda$. This result provides the
first examples of {\it cocycle superrigid} group actions, i.e.
actions having ALL cocycles with values into ANY discrete group
cohomologous to group morphisms. The proof is very similar to (4.2
in [P1]), which served as a model for ([P2,3], [PSa]) as well. We
use a von Neumann algebra framework, as this is particularly
suitable for the ``deformation/rigidity'' arguments involved. In
fact, this setting allows us to prove that given any closed subgroup
$\Cal V$ of the unitary group of a finite von Neumann algebra, all
$\Cal V$-valued cocycles for $\Gamma \curvearrowright X$ can be
untwisted. We use the same framework to show that if a group action
is cocycle superrigid then it is OE superrigid. Partial cases of
this general ``principle'' have been known for some time (see 4.2.9,
4.2.11 in [Z2], 3.3 in [Fu2], 2.4 in [Fu3]).

By using cocycles to study orbit equivalence of actions, we adopt in
this paper an approach pioneered by Zimmer in the late 70's and
which has been used, in one form or another, in many OE rigidity
results obtained so far ([Z2], [Ge], [GeGo], [Fu1,2,3], [MoS1,2]).
The past effectiveness of this method strongly motivated us to seek
a suitable cocycle superrigidity result behind the OE superrigidity
phenomena for malleable actions found in [P3]. In this respect,
cocycle superrigidity with discrete targets is best suited for OE
rigidity applications. In fact, as explained above, it can be viewed
as a direct generalization (and strengthening) of OE superrigidity.
Zimmer's cocycle superrigidity ([Z2]), in turn, is for linear
algebraic groups as targets, thus extending Margulis' superrigidity
but making it non-trivial to apply towards OE superrigidity. It was
only after developing  a series of new techniques that Furman could
take full advantage of it to prove OE superrigidity of higher rank
lattices ([Fu1,2,3]).

To state in more details the results in this paper, we need to
recall some definitions and terminology (also used in [P2,3,4]).
Thus, an inclusion of groups $H\subset \Gamma$ has the {\it relative
property} (T) of Kazhdan-Margulis ([M], [K]) if any unitary
representation of $\Gamma$ that almost contains the trivial
representation of $\Gamma$ must contain the trivial representation
of $H$. We also call such $H$ a {\it rigid subgroup} of $\Gamma$.
We'll consider two ``weak normality'' conditions for infinite
subgroups $H \subset \Gamma$. Thus, $H$ is {\it w-normal} (resp.
{\it wq-normal}) in $\Gamma$ if there exists a sequence of
consecutive normal inclusions $H_i \subset H_{i+1}$ (resp. with
$H_{i+1}$ generated by elements $g\in \Gamma$ with $|gH_ig^{-1}\cap
H_i|=\infty$) reaching from $H$ to $\Gamma$ (see 5.1 for the precise
definitions). Infinite property (T) groups and the groups $\Gamma_0
\ltimes \Bbb Z^2$ with $\Gamma_0 \subset SL(2,\Bbb Z)$ non-amenable
have infinite normal rigid subgroups by ([K], [B]), and so are
groups of the form $\Gamma_0 \ltimes \Bbb Z^n$, with $\Gamma_0$
arithmetic lattice in a classical Lie group and suitable $n$ ([Va],
[Fe]). If $H \subset \Gamma$ is wq-normal rigid then $H \subset
(\Gamma * \Gamma_0) \times \Gamma_1$ is wq-normal rigid, $\forall
\Gamma_1$ infinite. Both properties are closed to normal extensions.

Roughly, an action $\Gamma \curvearrowright^\sigma X$ is {\it
s-malleable} if the flip automorphism on $X \times X$ is in the
connected component of the identity in the centralizer of the double
action $\sigma_g \times \sigma_g, g\in \Gamma$ (an additional
``symmetry'' requirement is in fact needed, see 4.3 for the exact
definition). A typical example of s-malleable action is the {\it
Bernoulli $\Gamma$-action} on $(X,\mu)=\Pi_{g\in \Gamma}
(X_0,\mu_0)_g$, with $(X_0,\mu_0)_g$ identical copies of a standard
probability space $(X_0,\mu_0)$, given by
$\sigma_g((t_h)_h)=(t_{g^{-1}h})_h$, $(t_h)_h \in X$. More
generally, if $\Gamma$ acts on a countable set $K$ and
$(X,\mu)=\Pi_{k\in K} (X_0,\mu_0)_k$ with $\Gamma$ acting on
$(t_k)_k \in X$ by $\sigma_g((t_k)_k)=(t_{g^{-1}k})_k$, then
$\sigma$ is called a {\it generalized Bernoulli $\Gamma$-action} and
it is still s-malleable. Note that Bernoulli actions are always
(strongly) mixing, while a generalized Bernoulli action is weak
mixing iff $|\Gamma k|=\infty, \forall k$.

If $\sigma$ is a $\Gamma$-action on $(X,\mu)$ and $\Cal V$ a Polish
group, then a (right) $\Cal V$-valued measurable {\it cocycle} for
$\sigma$ is a measurable map $w: X \times \Gamma \rightarrow \Cal V$
satisfying for each $g_1,g_2\in \Gamma$ the identity
$w(t,g_1)w(g_1^{-1}t,g_2)=w(t,g_1g_2)$, $\forall t\in X$ (a.e.).
(N.B. All results below hold true, of course, for left cocycles as
well. We choose to work with right measurable cocycles in this paper
because in the von Neumann algebra framework, which is used for the
proofs, they become left cocycles.) Two cocycles $w,w'$ for $\sigma$
are {\it cohomologous} (or {\it equivalent}) if there exists a
measurable map $u: X \rightarrow \Cal V$ such that for each $g\in
\Gamma$ one has $w'(t,g)=u(t)^{-1} w(t,g)u(g^{-1}t)$, $\forall t\in
X$ (a.e.). Note that a cocycle $w$ is independent of the $X$
variable iff it is a group morphism of $\Gamma$ into $\Cal V$.

A Polish group $\Cal V$ is {\it of finite type} if it is isomorphic
to a closed subgroup of the group of unitary elements $\Cal U(N)$ of
a countably generated finite von Neumann algebra $N$, equivalently a
von Neumann algebra $N$ having a faithful normal trace state $\tau$
such that $N$ is separable in the Hilbert norm
$\|x\|_2=\tau(x^*x)^{1/2}, x\in N$. Countable discrete groups and
separable compact groups are of finite type, as they can be embedded
as closed subgroups of their group von Neumann algebra ([MvN1,2]).
However, the only connected locally compact groups of finite type
are the groups $\Cal V=K \times V$ with $K$ compact and $V\simeq
\Bbb R^n$ a vector group (cf. [KaSi], [vNS]).

\proclaim{0.1. Theorem (Cocycle superrigidity)} Let $\Gamma
\curvearrowright^\sigma X$ be a s-malleable action (e.g. a
generalized Bernoulli $\Gamma$-action) and assume $\Gamma$ has an
infinite rigid subgroup $H$ such that either $H$ is wq-normal with
$\sigma$ mixing, or that $H$ is w-normal with $\sigma_{|H}$ weak
mixing. Let $\Cal V$ be a Polish group of finite type. Then any
$\Cal V$-valued cocycle for $\sigma$ is cohomologous to a group
morphism of $\Gamma$ into $\Cal V$. More generally, if $\sigma'$ is
an action of the form $\sigma'_g=\sigma_g \times \rho_g \in
\text{\rm Aut}(X\times Y, \mu\times \nu), g\in\Gamma$, where $\rho$
is an arbitrary $\Gamma$-action on a standard probability space
$(Y,\nu)$, then any $\Cal V$-valued cocycle $w$ for $\sigma'$ is
cohomologous to a $\Cal V$-valued cocycle $w'$ which is independent
on the $X$-variable (i.e. $w'$ comes from a cocycle of $\rho$).
\endproclaim

\proclaim{0.2. Corollary} Let $\Gamma$ be a discrete group having
infinite wq-normal rigid subgroups and $\Cal V$ a Polish group of
finite type. Then any $\Cal V$-valued cocycle for a Bernoulli
$\Gamma$-action is cohomologous to a group morphism of $\Gamma$ into
$\Cal V$.
\endproclaim

For the next result we denote by  $\Cal R_\theta$ the equivalence
relation given by the orbits of a m.p. action $\Lambda
\curvearrowright^\theta Y$ of a countable group $\Lambda$ on a
probability space $(Y,\nu)$. More generally, if $Y_0\subset Y$ is a
subset of positive measure then $\Cal R^{Y_0}_\theta$ denotes the
equivalence relation on $Y_0$ given by the intersection of the
orbits of $\theta$ and the set $Y_0$. If $\theta$ is free ergodic
then this is easily seen to only depend on $\mu(Y_0)$ ([Dy]), up to
{\it isomorphism} of equivalence relations, i.e. up to isomorphism
of probability spaces taking the orbits of one relation onto the
orbits of the other (a.e.). As a consequence, it follows that if
$t>0$ and we take $m \geq t$, $\tilde{\theta}$ the product action of
$\tilde{\Lambda}=\Lambda \times \Bbb Z/m\Bbb Z$ on the product
probability space $\tilde{Y}=Y \times \Bbb Z/m\Bbb Z$ and
$Y_0\subset \tilde{Y}$ a subset of (product) measure $t/m$, then the
isomorphism class of $\Cal R^{Y_0}_{\tilde{\theta}}$ only depends on
$t$, not on the choice of $m$ and $Y_0\subset \tilde{Y}$. We call it
the {\it amplification of $\Cal R_\theta$ by $t$} and denote it
$\Cal R_\theta^t$.

Unlike Theorem 0.1, where the $\Gamma$-action $\sigma$ doesn't need
to be free, in the OE rigidity results below we have to assume
freeness. If $\sigma$ is a generalized Bernoulli action coming from
an action of $\Gamma$ on a set $K$ as before, then the condition
$|\{k\in K \mid gk \neq k\}|=\infty$, $\forall g\in \Gamma\setminus
\{e\}$, insures that $\sigma$ is free.

\proclaim{0.3. Theorem (OE superrigidity)} Let $\Gamma
\curvearrowright^\sigma X$ be as in $0.1.$ and assume in addition
that $\Gamma$ has no nontrivial finite normal subgroups and $\sigma$
is free. Let $\theta$ be an arbitrary free ergodic measure
preserving action of a countable discrete group $\Lambda$ on a
standard probability space $(Y,\nu)$. If $\Delta$ is an isomorphism
of probability spaces which takes $\Cal R_\sigma$ onto $\Cal
R_\theta^t$, for some $t > 0$, then $n=t^{-1}$ is an integer and
there exist a subgroup $\Lambda_0 \subset \Lambda$ of index
$[\Lambda: \Lambda_0]=n$, a subset $Y_0 \subset Y$ of measure
$\nu(Y_0)=1/n$ fixed by $\theta_{|\Lambda_0}$, an inner automorphism
$\alpha \in \text{\rm Inn}(\Cal R_{\theta})$ and a group isomorphism
$\delta: \Gamma \simeq \Lambda_0$ such that $\alpha \circ \Delta$
takes $X$ onto $Y_0$ and conjugates the actions $\sigma, \theta_0
\circ \delta$, where $\theta_0$ denotes the action of $\Lambda_0$ on
$Y_0$ implemented by $\theta$.

If moreover $\Gamma$ is infinite conjugacy class then any quotient
of $\Gamma \curvearrowright^\sigma X$ is OE superrigid, i.e. any OE
of $\sigma$ and an arbitrary free action $\Lambda \curvearrowright
Y$ comes from a conjugacy.
\endproclaim

\proclaim{0.4. Theorem (Superrigidity of embeddings)} Let $\Gamma
\curvearrowright^\sigma X$ be as in $0.1$, $\Lambda
\curvearrowright^\theta Y$ an arbitrary free ergodic action and
$t>0$.

If $\Delta: (X,\mu) \simeq (Y,\nu)^t$ is an identification of $\Cal
R_\sigma$ with a subequivalence relation of $\Cal R^t_\theta$ such
that any $\Gamma$-invariant finite subequivalence relation of $\Cal
R_\theta^t$ must be contained in $\Cal R_\sigma$, then $t\leq 1$ and
there exists an isomorphism $\delta:\Gamma\simeq \Lambda_0 \subset
\Lambda$ and $\alpha \in \text{\rm Inn}(\Cal R_{\theta})$ such that
$\alpha \circ \Delta$ takes $X$ onto a $\Lambda_0$-invariant subset
$Y_0\subset Y$, with $\nu(Y_0)=t$, and conjugates the actions
$\sigma, \theta_{|\Lambda_0}$ with respect to the identification
$\delta:\Gamma \simeq \Lambda_0$.
\endproclaim

An interesting application of 0.4 is as follows: Let $\mu_0$ be a
probability measure on $\{0,1\}$ with unequal weights, i.e.
$s=\mu_0(\{0\})/\mu_0(\{1\})\neq 1$, let
$(X,\mu)=(\{0,1\},\mu_0)^\Gamma$ and denote by $\Cal R$ the
countable equivalence relation on $X$ generated by the Bernoulli
$\Gamma$-action and the relation $\Cal R_0$ given by $(t_g)_g \sim
(t'_g)_g$ if $\exists F \subset \Gamma$ finite such that $t_g=t'_g,
\forall g\in \Gamma \setminus F$, $\Pi_{g\in F} \mu_0(t_g) =
\Pi_{g\in F} \mu_0(t'_g)$. If $\Gamma$ is w-rigid and has no finite
normal subgroups then $\mycal F(\Cal R) \supset s^{\Bbb Z}$ (see
[P2]; equality is in fact shown for certain $\Gamma$, such as
$\Gamma=SL(2,\Bbb Z)\ltimes \Bbb Z^2$) and it is easy to see that
any finite subequivalence relation of $\Cal R$ invariant to the
action of $\Gamma$ on $\Cal R$ given by $g(t,t')=(gt,gt')$ must be contained in
$\Cal R_{\Gamma}$. Thus, by 0.4 it follows that $\Cal R^t$ cannot be
implemented by a free action of a group, $\forall t>0$. The first
example of an equivalence relation with the property that all its
amplifications are not implementable by free group actions was
obtained in ([Fu2]).

The ideas behind the proofs of 0.1-0.4 (notably the
deformation/rigidity arguments used in 0.1) are quintessentially
``von Neumann algebra'' in spirit. This made us favor a von Neumann
algebra framework for the presentation, rather than a measure
theoretical one. The two points of view are in fact equivalent, due
to a well known observation showing that if $\Gamma
\curvearrowright^\sigma X$, $\Lambda \curvearrowright^\theta Y$ are
free m.p. actions then an isomorphism of probability spaces
$\Delta:X \simeq Y$ gives an OE of $\sigma, \theta$ if and only if,
when regarded as an algebra isomorphism $\Delta : L^\infty X \simeq
L^\infty Y$, $\Delta$ extends to an isomorphism of the von Neumann
algebras $L^\infty X \rtimes \Gamma\simeq L^\infty Y \rtimes
\Lambda$ (cf. [Si], [Dy1,2], [FM]).

We summarize in Section 1 the tools from the theory of von Neumann
algebras that we need in this paper, for convenience. In Section 2
we introduce the class of Polish groups of finite type and explain
how measurable cocycles with values in such groups can be viewed as
cocycles for actions on von Neumann algebras. Also, we discuss a
{\it relative weak mixing} condition for actions on von Neumann
algebras considered in ([P2]), which generalizes a concept
introduced in the measure theoretic context by Furstenberg ([F]) and
Zimmer ([Z3]) in the 1970's. It plays an important role in this
paper. We show for instance that if $\Gamma \curvearrowright^\sigma
X$ is a quotient of a cocycle superrigid action $\Gamma
\curvearrowright^{\sigma'} X'$ and the latter is weak mixing
relative to the former, then $\Gamma \curvearrowright X$ is cocycle
superrigid as well. In Section 3 we prove a key criterion for
``untwisting'' cocycles with values in the unitary group of a finite
von Neumann algebra $N$, extracted from proofs in ([P1,2]). It shows
that a $\Cal U(N)$-valued cocycle $w$ for a weak mixing action
$\Gamma \curvearrowright X$ is equivalent to a group morphism of
$\Gamma$ into $\Cal U(N)$ iff the cocycles $w^l(t,s,g)=w(t,g)$ and
$w^r(t,s,g)=w(s,g)$ for the double action $\sigma_g \times \sigma_g$
are equivalent. A similar statement holds true for arbitrary
(non-commutative) finite von Neumann algebras. We also prove a
hereditary result for weak mixing actions showing that in order to
untwist a $\Cal V$-valued cocycle, for some Polish group $\Cal V$ of
finite type, it is sufficient to untwist it as a cocycle with values
in a larger Polish group of finite type.

In Section 4 we recall from ([P1-4]) the notion of {\it s-malleable}
action and the proof that (generalized) Bernoulli actions have this
property, then use a deformation/rigidity argument from ([P1,2]) to
show that if a group $\Gamma$ contains a large {\it rigid} part then
any s-malleable $\Gamma$-action satisfies the criterion for
untwisting cocycles from Section 3. In Section 5 we derive the
cocycle superrigidity of s-malleable (e.g. Bernoulli) actions with
values in Polish groups of finite type. We then prove a general
principle showing that if $\Gamma \curvearrowright X, \Lambda
\curvearrowright X$ have the same orbits, $w$ denotes the associated
cocycle and $v:X \rightarrow \Lambda$ implements an equivalence of
$w$ with a group morphism $\delta:\Gamma \rightarrow \Lambda$, then
$v$ gives rise to an automorphism $\alpha$ with graph in $\Cal
R_\sigma=\Cal R_\theta$ which conjugates $\sigma$ and $\theta\circ
\delta$. In particular, this shows that cocycle superrigidity
implies OE superrigidity, thus yielding the OE applications. Another
general principle we prove is that any quotient of a cocycle
superrigid weak mixing action of an ICC group is OE superrigid.

A first account with proofs of the results in this paper has been
presented in the Seminaire d'Alg\`ebres d'Op\'erateurs in Paris,
Sept. 2005. It is a pleasure for me to thank Stefaan Vaes and the
members of the seminar for their stimulating comments. I am also
very grateful to Alex Furman and Stefaan Vaes for many useful
discussions and remarks on the preliminary version of the
paper (math.GR/0512646). Last but not least I would like to thank  
the referee for a very careful  
reading and useful suggestions.

\heading 1. Preliminaries
\endheading

Although we are interested in actions of groups on the probability
space, our approach will be functional analytical, using von Neumann
algebra framework. This section is intended for readers who are less
familiar with this field. Thus, we'll summarize here some basic
tools in von Neumann algebras such as: the standard representation
of a finite von Neumann algebra with a trace; the crossed product
(resp. the group measure space) construction of a von Neumann
algebra starting from an action of a discrete group on a finite von
Neumann algebra (resp. on a probability space); orbit equivalence of
actions as isomorphisms of group measure space von Neumann algebras;
von Neumann subalgebras and the basic construction. Some knowledge
in functional analysis and the spectral theorem should be sufficient
to recover the omitted proofs.

\vskip .05in

\noindent {\bf 1.1. Probability spaces as abelian von Neumann
algebras}. The ``classical'' measure theoretical approach is
equivalent to a ``non-classical'' operator algebra approach due to a
well known observation of von Neumann, showing that measure
preserving isomorphisms between standard probability spaces
$(X,\mu)$, $(Y,\nu)$ are in natural correspondence with $^*$-algebra
isomorphisms between their function algebras $L^\infty X$, $L^\infty
Y$ preserving the functional given by the integral, $\tau_\mu=\int
\cdot \text{\rm d}\mu$, $\tau_\nu=\int \cdot \text{\rm d}\nu$. More
generally:

\vskip .05in \noindent $(1.1.1)$. Let $\Delta: (X,\mu) \rightarrow
(Y,\nu)$ be a measurable map with $\nu \circ \Delta=\mu$. Then
$\Delta^*:L^\infty Y \rightarrow L^\infty X$ defined by
$\Delta^*(x)(t)=x(\Delta t), t \in X$, is an injective $^*$-algebra
morphism satisfying $\tau_{\mu} \circ \Delta^*=\tau_{\nu}$.
Conversely, if $(X,\mu), (Y,\nu)$ are probability spaces and
$\rho:L^\infty Y \rightarrow L^\infty X$ is an injective $*$-algebra
morphism such that $\tau_\mu \circ \rho = \tau_\nu$, then there exists
a measurable map $\Delta:X \rightarrow Y$ such that $\rho=\Delta^*$.
Moreover, $\Delta$ is unique and onto, modulo a set of measure $0$,
and the correspondence $\Delta \mapsto \Delta^*$ is
``contravariant'' functorial, i.e. $(\Delta \circ
\Delta')^*={\Delta'}^* \circ \Delta^*$. Also, $\Delta$ is a.e. $1$
to $1$  if and only if $\Delta^*$ is onto and if this is the case
then $\Delta^{-1}$ is also measurable and measure preserving.

\vskip .05in

There are two norms on $L^\infty X$ that are relevant for us in this
paper, namely the ess-sup norm $\|\cdot \|= \|\cdot \|_\infty$ and
the norm $\|\cdot \|_2$. Note that the unit ball $(L^\infty X)_1$ of
$L^\infty X$ (in the norm $\|\cdot \|$) is complete in the norm
$\|\cdot \|_2$. We will often identify $L^\infty X$ with the von
Neumann algebra of (left) multiplication operators $L_x, x\in
L^\infty X$, where $L_x(\xi)=x\xi$, $\xi \in L^2 X$. The
identification $x \mapsto L_x$ is a $^*$-algebra morphism, it is
isometric (from $L^\infty X$ with the ess-sup norm into $\Cal B(L^2
X)$ with the operatorial norm) and takes the $\|\cdot \|_2$-toplogy
of $(L^\infty X)_1$ onto the strong operator topology on the image.
Also, the integral $\tau_{\mu}(x)$ becomes the vector state $\langle
L_x (1),1 \rangle$, $x\in L^\infty X$. Moreover, if $\Delta:(X,\mu)
\simeq (Y,\nu)$ for some other probability space $(Y,\nu)$, then
$\rho={\Delta^{-1}}^*$ extends to an (isometric) isomorphism of
Hilbert spaces $L^2 X\simeq L^2 Y$ which conjugates the von Neumann
algebras $L^\infty X\subset \Cal B(L^2 X)$, $L^\infty Y\subset \Cal
B(L^2 Y)$ onto each other.

With this in mind, let us denote by Aut$(X,\mu)$ the group of
(classes modulo null sets of) {\it measure preserving automorphisms}
$T:(X,\mu)\simeq (X,\mu)$ of the standard probability space
$(X,\mu)$. Denote Aut$(L^\infty X,\tau_\mu)$ the group of
$^*$-automorphisms of the von Neumann algebra $L^\infty X$ that
preserve the functional $\tau_\mu$, and
identify Aut$(X,\mu)$ and Aut$(L^\infty X,\tau_\mu)$ via the map $T
\mapsto (T^{-1})^*$. One immediate benefit of the functional
analysis framework and of this identification is that it gives a
natural Polish group topology on Aut$(X,\mu)$, given by pointwise
$\|\cdot \|_2$-convergence in Aut$(L^\infty X, \tau_\mu)$, i.e.
$\vartheta_n \rightarrow \vartheta$ in Aut$(L^\infty X,\tau_\mu)$ if
$\lim_n \|\vartheta_n(x)-\vartheta(x)\|_2 =0$, $\forall x\in
L^\infty X$.

An {\it action} of a discrete group $\Gamma$ on the standard
probability space $(X,\mu)$ is a group morphism $\sigma: \Gamma
\rightarrow \text{\rm Aut}(X,\mu)$. Using the identification
between Aut$(X,\mu)$ and Aut$(L^\infty X, \tau_\mu)$, we
alternatively view $\sigma$ as an action  of $\Gamma$ on $(L^\infty
X,\tau_\mu)$, i.e as a group morphism  $\sigma: \Gamma \rightarrow
\text{\rm Aut} (L^\infty X,\tau_\mu)$. Although we use the same
notation for both actions, the difference will be  clear from the
context. Furthermore, when viewing $\sigma$ as an action on the
probability space $(X,\mu)$, we'll use the simplified notation
$\sigma_g(t)=gt$, for $g\in \Gamma, t\in X$. The relation between
$\sigma$ as an action on $(X,\mu)$ and respectively on $(L^\infty X,
\tau_\mu)$ is then given by the equations
$\sigma_g(x)(t)=x(g^{-1}t)$, $\forall t\in X$ (a.e.), which hold
true for each $g\in \Gamma$, $x\in L^\infty X$.

The action $\sigma$ is {\it free} if for any $g\in \Gamma, g\neq e$,
the set $\{t\in X \mid gt=t\}$ has $\mu$-measure 0. The action is
{\it ergodic} if $X_0 \subset X$ measurable with $gX_0=X_0$ (a.e.)
for all $g\in \Gamma$, implies $X_0=X$ or $X_0=\emptyset$ (a.e.), in
other words $\mu(X_0)=0,1$.

\vskip .05in

\noindent {\bf 1.2. Finite von Neumann algebras}. When  viewed as an
abelian von Neumann algebra with its integral functional, the
natural generalization of a probability space is a von Neumann
algebra $N$ with a linear
functional $\tau:N \rightarrow \Bbb C$
satisfying the conditions: $\tau(x^*x) \geq 0, \forall x\in N$ and
$\tau(1)=1$ ($\tau$ is a {\it state}); $\tau(x^*x)=0$ iff $x=0$
($\tau$ is {\it faithful}); the unit ball $(N)_1=\{ x\in B \mid
\|x\|\leq 1\}$ of $N$ is complete in the Hilbert norm
$\|x\|_2=\tau(x^*x)^{1/2}$ ($\tau$ is {\it normal});
$\tau(xy)=\tau(yx), \forall x,y \in N$
($\tau$ is a {\it trace}).
Such $(N,\tau)$ is called a {\it finite} von Neumann algebra (with
its {\it trace} $\tau$). We say that $(N,\tau)$ is {\it separable}
if it is separable in the Hilbert norm $\|x\|_2=\tau(x^*x)^{1/2},
x\in N$. The trace $\tau$ on a
finite von Neumann algebra $N$ is in general not unique, but
there does exist a unique expectation of $N$ onto its center $\Cal Z(N)$,
denoted $Ctr_N$ and called the {\it central trace} on $N$,
with the property that $Ctr_N(xy)=Ctr_N(yx), \forall x,y\in N$ ([D2]).
Any (faithful normal) trace $\tau$ on $N$ is the composition of $Ctr_N$ with
a (faithful normal) state on $\Cal Z(N)$. In particular, if $N$ is a
{\it factor}, i.e. $\Cal Z(N)=\Bbb C$, then $N$ has a unique
trace state $\tau=Ctr_N$, which is automatically normal and faithful.

If $N$ is a finite von Neumann algebra and $p,q$ are
(selfadjoint) projections
in $N$ then $p \prec q$
(resp. $p \sim q$) in $N$, i.e. there exists a partial isometry $v\in N$
with $vv^*=p, v^*v \leq q$ (resp. $vv^*=p, v^*v=q$), if and only
if $Ctr_N(p)\leq Ctr_N(q)$ (resp $Ctr_N(p)=Ctr_N(q)$).
In particular, if $1 \sim q$ for some projection $q$ in $N$
then $q=1$. By a celebrated theorem of Murray and von Neumann ([MvN1]),
this purely algebraic condition is sufficient to ensure that $N$ has
a central trace, and thus a faithful normal trace state as well.

The representation of $L^\infty X$ as an algebra of left
multiplication operators on $L^2 X$ generalizes to the
non-commutative setting of finite von Neumann algebras $(N,\tau)$ as
follows: Denote by $L^2 N = L^2(N,\tau)$ the completion of $N$ in
the norm $\|\cdot \|_2$. Then each element $x\in N$ defines an
operator of left multiplication $L_x$ on $L^2 N$, by $L_x
(\xi)=x\xi$, $\xi\in L^2 N$. The map $N \ni x \mapsto L_x \in \Cal
B(L^2 X)$ is clearly a $^*$-algebra morphism. Due to the
faithfulness of $\tau$, it preserves the operatorial norm on $N$,
with the normality of $\tau$ insuring that the image $L_N=\{L_x \mid
x\in N\}$ is weakly closed in $\Cal B(L^2 X)$, i.e. $L_N$ is a von
Neumann algebra. It can be easily shown that the commutant $L_N'$ of
$L_N$ in $\Cal B(L^2 N)$ is equal to the algebra $R_N$ of operators
of right multiplication by elements in $N$, and conversely
$R_N'=L_N$. Also, by the definition, we have $\langle x 1,1
\rangle=\tau(x)$. We will always identify $N$ with its image
$L_N\subset \Cal B(L^2 N)$ and call this the {\it standard
representation} of $(N,\tau)$.

Note that if $(N,\tau)=(L^\infty X,\tau_\mu)$ then $L^2(N,\tau)$
coincides with the Hilbert space $L^2 X$ of square integrable
functions on $(X,\mu)$. In this case $\xi \in L^2 X$ can be viewed
as the closed (densely defined) operator of (left) multiplication by
$\xi$, whose spectral resolution lies in $L^\infty X$.

When $(N,\tau)$  is an arbitrary finite von Neumann algebra one can
still interpret the elements in $L^2 N$ as closed linear operators
on $L^2 N$, as follows: For each $\xi \in L^2 N$ let $L^0_\xi$ be
the linear operator with domain $N$ (regarded as a vector subspace
of $L^2N$) defined by $L_\xi(x)=\xi x$, $x\in N$. This operator
extends to a unique closed operator $L_\xi$ which commutes with
$R_N$, equivalently its polar decomposition $L_\xi=u |L_\xi|$ has
both the partial isometry $u$ and the spectral resolution
$\{e_s\}_{s > 0}$ of $|L_\xi|=(L_\xi^* L_\xi)^{1/2}$ lying in $N$.
Closed operators satisfying this property are said to be {\it
affiliated with} $N$. In addition, $L_\xi$ is {\it square
integrable}, i.e. $\tau(L_\xi^* L_\xi) \overset \text{\rm def} \to =
\int s^2 \text{\rm d} \tau(e_s) = \| L_\xi(1)\|^2_2 =\|\xi\|_2^2 <
\infty$. Noticing that $L_\xi(1)=\xi$, it follows that $\xi \mapsto
L_\xi$ gives a 1 to 1 correspondence between $L^2 N$ and the space
of square integrable operators affiliated with $N$. We will always
view elements of $L^2 N$ in this manner.

If $\xi, \eta \in L^2 N$ then their product (composition) as closed
operators $\xi \cdot \eta$ is also a densely defined operator
affiliated with $N$, i.e. it is of the form $u b$ with $u$ partial
isometry in $N$ and $b = |\xi \eta|$ a positive operator with
spectral resolution $\{e_s\}_{s > 0}$ in $N$ and $\tau(b) \overset
\text{\rm def} \to = \int s \text{\rm d} \tau(e_s) < \infty $. We
denote by $L^1 N$ the set of all such operators, which is easily
seen to be a Banach space when endowed with the norm $\|ub\|_1 =
\tau(b)$. Also, it has  $L^2 N \supset N$ as dense subspaces. Any
$\zeta =  \xi\cdot\eta^* \in L^1(M,\tau)$ defines a functional on
$M$ by $\tau( x\zeta)\overset \text{\rm def} \to =\langle x\xi, \eta
\rangle$, $x\in M$, which is positive iff $\zeta \geq 0$  as an
operator. The norm of $\zeta$ as a functional on $M$ coincides with
$\|\zeta\|_1$ and in fact, as a space of functionals on $M$, $L^1 M$
is the predual of $M$, i.e.  $(L^1 M)^*=M.$

\vskip .05in \noindent {\bf 1.3. Actions of groups and crossed
product algebras}. Like in the commutative case, we denote by
Aut$(N,\tau)$ the group of automorphisms of the finite von Neumann
algebra $(N,\tau)$ (i.e. the $\tau$-preserving $^*$-algebra
isomorphisms of $N$ onto $N$), and endow it with the Polish group
topology given by point-wise $\|\cdot \|_2$-convergence. Note that
any $\vartheta \in {\text{\rm Aut}}(N,\tau)$ preserves the Hilbert
norm $\|\cdot \|_2$ and thus extends to a unitary operator on $L^2
N$ which, if no confusion is possible, will still be denoted
$\vartheta$.

Given a discrete group $\Gamma$, an action $\sigma$ of $\Gamma$ on
$(N,\tau)$ is a group morphism $\sigma: \Gamma \rightarrow \text{\rm
Aut}(N,\tau)$. Since any $\sigma_g$ extends to a unitary operator on
$L^2N$, $\sigma$ extends to a unitary representation of $\Gamma$ on
the Hilbert space $L^2N$, still denoted $\sigma$.

We use the (rather standard) notation $N^\sigma$ to designate the
{\it fixed point algebra} of the action, $N^\sigma=\{x\in N \mid
\sigma_g(x)=x, \forall g\in \Gamma\}$.

A key tool in the study of actions is the {\it crossed product
construction}, which associates to $(\sigma,\Gamma)$ the von Neumann
algebra $N \rtimes_\sigma \Gamma$ generated on the Hilbert space
$\Cal H = L^2N \overline{\otimes} \ell^2(\Gamma)$ by a copy of the
algebra $N$, acting on $\Cal H$ by left multiplication on $L^2N$,
and a copy of the group $\Gamma$, acting on $\Cal H$ as the
operators $u_g = \sigma_g \otimes \lambda_g$, where $\sigma_g, g\in
\Gamma,$ is viewed here as unitary representation. In fact,
$\{u_g\}_g$ is easily seen to be a multiple of left regular
representation. In case $(N,\tau)=(L^\infty X, \mu)$ with $\sigma$
coming from an action of $\Gamma$ on $(X,\mu)$, this amounts to
Murray-von Neumann's {\it group measure space construction}
([MvN1]).

The following more concrete description of $M=N \rtimes_\sigma
\Gamma$ and its standard representation is quite useful: Identify
$\Cal H=\ell^2(\Gamma, L^2N)$ with the Hilbert space of
$\ell^2$-summable formal sums $\Sigma_g \xi_g u_g$, with
``coefficients'' $\xi_g$ in $L^2N$ and ``indeterminates''
$\{u_g\}_g$ labeled by the elements of the group $\Gamma$. Define a
$^*$-operation on $\Cal H$ by $(\Sigma_g \xi_g u_g)^*=\Sigma_g
\sigma_g(\xi_{g^{-1}}^*)u_g$ and let both $N$ and the $u_g$'s act on
$\Cal H$ by left multiplication, subject to the product rules $y(\xi
u_g)=(y\xi)u_g$, $u_g (\xi u_h)=\sigma_{g}(\xi)u_{gh}, \forall g,h
\in G$, $y\in N$, $\xi \in L^2N$. In fact, given any $\xi=\Sigma_g
\xi_g u_g, \zeta = \Sigma_h \zeta_h u_h \in \Cal H$ one can define
the product $\xi \cdot \zeta$ as the formal sum $\Sigma_k \eta_k
u_k$ with coefficients $\eta_k = \Sigma_g \xi_g \zeta_{g^{-1}k}$,
the sum being absolutely convergent in the norm $\|\cdot \|_1$, with
estimates $\|\eta_k\|_1 \leq \|\xi\|_2 \|\zeta\|_2$, $\forall k\in
\Gamma$, by the Cauchy-Schwartz inequality.

$\xi \in \Cal H$ is a {\it convolver} if $\xi \zeta \in \Cal H$
(i.e. with the above notations $\eta_k \in L^2 N$ and $\Sigma_k
\|\eta_k\|_2^2 < \infty$) for all $\zeta \in \Cal H$. $M$ is then
the algebra of all left multiplication operators $\zeta \mapsto \xi
\cdot \zeta$ by convolvers $\xi$. Its commutant in $\Cal B(\Cal H)$
is the algebra of all right multiplication operators $\zeta \mapsto
\zeta\xi$ by convolvers $\xi$. If $T \in M$ then $\xi=T(1)\in \Cal
H$ is a convolver and $T$ is the operator of left multiplication by
$\xi$, with $T^*$ corresponding to the left multiplication by
$\xi^*$. The left multiplication by convolvers supported on $N=Nu_e$
gives rise to a multiple of the standard representation of $N$,
while the left multiplication by the convolvers $\{u_g\}_g$ gives
rise to a multiple of the left regular representation of $\Gamma$.
The trace $\tau$ on $N$ extends to all $N \rtimes_\sigma G$ by
$\tau(\Sigma_g y_g u_g)=\tau(y_e)=\langle \xi, 1\rangle = \langle
\xi \cdot 1, 1\rangle$, where $\xi=\Sigma_g y_g u_g$. The Hilbert
space $\Cal H$ naturally identifies with $L^2 M$, while its subspace
$M\subset L^2 M$ identifies with the set of convolvers and the
standard representation of $M$ with the algebra of left
multiplication by convolvers.

In case $N=L^\infty X$ and $\sigma$ comes from an m.p. action
$\Gamma \curvearrowright (X,\mu)$, the condition $\sigma$ free is
equivalent to $L^\infty X$ being maximal abelian in $M=L^\infty X
\rtimes \Gamma$, i.e. if $x\in M$, $[x,L^\infty X]=0$ then $x\in
L^\infty X$. If $\sigma$ is free, then $M$ is a factor if and only
if $\sigma$ is ergodic. If $\sigma$ is free and ergodic
(equivalently $L^\infty X$ maximal abelian in $M$ and $M$ is a
factor) then there are two possibilities: $\Gamma$ infinite, in
which case $M$ is a II$_1$ factor; $|\Gamma|=n < \infty$, in which
case $M=M_{n \times n}(\Bbb C)$ with the subalgebra $L^\infty X
\subset M_{n\times n}(\Bbb C)$ corresponding to a diagonal
subalgebra of $M_{n\times n}(\Bbb C)$.

\vskip .05in

\vskip .05in \noindent {\bf 1.4. Isomorphism of algebras from
equivalence of actions}. Let $\Gamma \curvearrowright^\sigma
(X,\mu)$ and $\Lambda \curvearrowright^\theta (Y,\nu)$ be free m.p.
actions of discrete groups on probability spaces. It is trivial to
see that if $\Delta:(X,\mu) \simeq (Y,\nu)$ gives a {\it conjugacy}
of $\sigma,\theta$, i.e. $\Delta\circ \sigma_g = \theta_{\delta(g)}
\circ \Delta$, $\forall g\in \Gamma$, for some group isomorphism
$\delta: \Gamma \simeq \Lambda$, then $(\Delta^{-1})^*: L^\infty X
\simeq L^\infty Y$ extends to an isomorphism of the group measure
space algebras $L^\infty X \rtimes \Gamma\simeq L^\infty Y \rtimes
\Lambda$, which takes a formal sum $\Sigma_g a_g u_g$ onto $\Sigma_g
(\Delta^{-1})^*(a_g) v_{\delta(g)}$, where $u_g\in L^\infty X
\rtimes \Gamma$, $v_h\in L^\infty Y \rtimes \Lambda$ are the
canonical unitaries.

It has been shown by Singer ([Si]), Dye ([Dy1,2]) and Feldman-Moore
([FM]) that in fact much less than conjugacy is sufficient: If
$\Delta:(X,\mu) \simeq (Y,\nu)$ is an isomorphism of probability
spaces then $(\Delta^{-1})^*: L^\infty X \simeq L^\infty Y$ extends
to an isomorphism of the corresponding group measure space von
Neumann algebras if and only if $\Delta$ is an {\it orbit
equivalence} (OE) of $\sigma,\theta$, i.e. if for almost all $t\in
X$ we have $\Delta(\Gamma t)=\Lambda \Delta(t)$.

This observation leads to the consideration of the {\it orbit
equivalence relation} $\Cal R_\sigma$ implemented on the probability
space $(X,\mu)$ by the orbits of a free m.p. action $\Gamma
\curvearrowright^\sigma (X,\mu)$, i.e. $(t, t')\in \Cal R_\sigma$ if
$\Gamma t = \Gamma t'$ ([Si], [FM]). More generally, if $X_0 \subset
X$ is a measurable subset then one denotes by $\Cal R^{X_0}_\sigma$
the equivalence relation on $X_0$ given by the intersection of the
orbits of $\sigma$ and the set $X_0$: $t,t'\in X_0$ are equivalent
if $\Gamma t \cap X_0 = \Gamma t' \cap X_0$. If $\Lambda
\curvearrowright^\theta (Y,\nu)$ is another free m.p. action and
$Y_0\subset Y$, $p_0=\chi_{X_0}, q_0=\chi_{Y_0}$, then an
isomorphism of probability spaces $\Delta:(X_0, \mu_{X_0}) \simeq
(Y_0, \nu_{Y_0})$ extends to a von Neumann algebra isomorphism
$p_0(L^\infty X \rtimes \Gamma)p_0\simeq q_0(L^\infty Y \rtimes
\Lambda)q_0$ if and only if for almost all $t\in X$ we have $\Delta
(\Gamma t \cap X_0)=\Lambda (\Delta(t)) \cap Y_0$, i.e. iff $\Delta$
takes the equivalence relation $\Cal R_\sigma \cap (X_0\times X_0)$
onto the equivalence relation $\Cal R_\theta \cap (Y_0 \times Y_0)$.

To explain this result, it is convenient to consider a more general
notion of equivalence relation (cf [FM]) and construct its
associated von Neumann algebra. Thus, an equivalence relation $\Cal
R$ on $X$ is called a {\it countable m.p. equivalence relation} if
there exists a countable subgroup $\Gamma \subset \text{\rm
Aut}(X,\mu)$ and a subset $N_0\subset X$ of measure 0 such that for
all $t\in X\setminus N_0$ the orbit of $t$ under $\Cal R$ coincides
with $\Gamma t$. The {\it full group} $[\Cal R]$ of the equivalence
relation $\Cal R$ is the set of all $\phi \in \text{\rm Aut}(X,\mu)$
with the property that there exists a null set $N_0\subset X$
such that the graph $G_\phi=\{(t,\phi(t))\mid t\in X\setminus N_0\}$ is
contained in $\Cal R$. In this same spirit, if $\Gamma\subset
\text{\rm Aut}(X,\mu)$ then $[\Gamma]$ denotes the group of
automorphisms $\phi$ of $(X,\mu)$ which are locally implemented by
elements on $\Gamma$, i.e. for which there exist a partition of $X$ with
measurable subsets $\{X_n\}_n$ and automorphisms $\phi_n\in \Gamma$
such that $\phi_{|X_n}={\phi_n}_{|X_n}, \forall n$. Note that if
$\Gamma\subset \text{\rm Aut}(X,\mu)$ is a countable group, then
$\Gamma$ implements $\Cal R$ iff $[\Gamma]=[\Cal R]$.

Similarly, the {\it full pseudogroup} of a countable m.p.
equivalence relation $\Cal R$ (resp. of a subgroup $\Gamma\subset
\text{\rm Aut}(X,\mu)$) is the set $_{_{p}}[\Cal R]$ (resp.
$_{_{p}}[\Gamma]$) of all measurable $\mu$-m.p isomorphisms $\phi:
r(\phi) \simeq l(\phi)$, with $r(\phi), l(\phi) \subset X$
measurable and graph contained in $\Cal R$ (resp. locally
implemented by elements in $\Gamma$). This is easily seen to
coincide with the set of all ``local isomorphisms'' of the form
$\phi_{|Y_0}$ with $\phi \in {_{_{p}}[\Cal R]}$ (resp $\phi\in
{_{_{p}}[\Gamma]}$) and $Y_0\subset X$. This set is endowed with a
product given by $\phi \psi(t)\overset \text{\rm def} \to
=\phi(\psi(t))$, $t \in r(\phi \psi)\overset \text{\rm def} \to
=\{t\in r(\psi) \mid \psi(t) \in r(\phi) \}$ and inverse given by
inverse of functions.

If $\Cal G$ is a full pseudogroup on $(X,\mu)$ (coming from either a
countable equivalence relation or a countable subgroup of
Aut$(X,\mu)$) then its associated von Neumann algebra $L(\Cal G)$ is
defined as follows: For each $\phi \in \Cal G$, $a\in L^\infty X$
let $\phi(a) \in L^\infty X$ be defined by $\phi(a)(t)=a
(\phi^{-1}(t))$, if $t\in l(\phi)$, $\phi(a)(t)=0$ if $t\not\in
l(\phi)$. Denote $L_0(\Cal G)$ the algebra of formal finite sums
$\Sigma_\phi a_\phi u_\phi$, with $a_\phi \in (L^\infty X) r(\phi)$
and ``indeterminates'' $u_\phi$, product rule given by $(a_\phi
u_\phi)(a_\psi u_\psi)=a_\phi \phi(a_\psi)u_{\phi\psi}$ and
$*$-operation given by $(a_\phi u_\phi)^*= \phi^{-1}(a_\phi)
u_{\phi^{-1}}$. Define $\tau(a_\phi u_\phi)=\int a_\phi i(\phi)
\text{\rm d} \mu$, where $i(\phi)$ is the characteristic function of
the largest set on which $\phi$ acts as the identity, then extend
$\tau$ to all $L_0(\Cal G)$ by linearity. Denote by $L^2(\Cal G)$
the Hilbert space obtained by completing $L_0(\Cal G)/I_\tau$ in the
norm $\|x\|_2=\tau(x^*x)^{1/2}$, where $I_\tau=\{x \mid \langle x,x
\rangle =0\}$. The $*$-algebra $L_0(\Cal G)$ acts on $L^2(\Cal G)$
by left multiplication and $L(\Cal G)$ is defined to be its weak
closure. In case $\Cal G=  {_{_{p}}[\Cal R]}$ for a countable m.p.
equivalence relation, we denote $L( _{_{p}}[\Cal R])$ by $L(\Cal R)$
(in the sense of 1.5 below).

Note that the algebra of coefficients $L^\infty X$ is maximal
abelian in $L(\Cal R)$ and that any unitary element $u\in L(\Cal R)$
which normalizes $L^\infty X$, i.e. $uL^\infty X u^*=L^\infty X$, is
of the form $u=au_\phi$, where $a$ is a unitary element in $L^\infty
X$ and $\phi \in [\Cal R]$. It is useful to notice that, by
maximality, there exist $\phi_n \in {_{_{p}}[\Cal R]}$ such that
$i(\phi_n^{-1}\phi_m) =0$, $\forall n\neq m$, and $\forall \phi \in
{_{_{p}}[\Cal R]}$, $\exists X^\phi_n$ partition of $r(\phi)$ such
that $\phi(t)=\phi_n(t)$, $\forall t\in X_n^\phi$ (a.e.). Note that
these conditions amount to saying that $\{u_{\phi_n}\}_n$ is an
orthonormal basis of $L([\Cal R])$ over $L^\infty X$ (in the sense
of 1.5 below).

\vskip .05in \noindent {\it 1.4.1. Definition}. If $\sigma$ is free
and ergodic then $\Cal R_\sigma^{X_0}$ only depends on $t=\mu(X_0)$,
up to isomorphism of equivalence relations (see e.g. [Dy1]). More
generally, let $t > 0$ then choose an integer $n \geq t$ and denote
by $\tilde{\sigma}$ the action of $\Gamma \times (\Bbb Z/n\Bbb Z)$
on $\tilde{X}=X \times \Bbb Z/n\Bbb Z$ given by the product of
$\sigma$ and the left translations by elements in $\Bbb Z/n\Bbb Z$.
Let $X_0\subset \tilde{X}$ be a subset of measure $t/n$. It is then
trivial to see that, up to isomorphism of equivalence relations,
$\Cal R_{\tilde{\sigma}}^{X_0}$ only depends on $t$ (not on the
choice of $n$ and $X_0$). We denote the isomorphism class of $\Cal
R^{X_0}_{\tilde{\sigma}}$ by $\Cal R^t_\sigma$ and call it the {\it
amplification of $\Cal R_\sigma$ by $t$}.

\vskip .05in

Let now $\Cal S$ be another countable m.p. equivalence relation on
the probability space $(Y, \nu)$ and let $\Delta:(X,\mu) \simeq
(Y,\nu)$ be an isomorphism of probability spaces. It is now trivial
from the definitions that $\Delta$ takes almost every orbit of $\Cal
R$ onto an orbit of $\Cal S$ iff $\Delta$ conjugates the full groups
$[\Cal R], [\Cal S]$, and also iff it conjugates the corresponding
full pseudogroups. From the definition of $L(\Cal R)$, this later
condition is clearly equivalent to the fact that $(\Delta^{-1})^*:
L^\infty X \simeq L^\infty Y$ extends to a von Neumann algebra
isomorphism of $L(\Cal R)$ onto $L(\Cal S)$. Such $\Delta$ is called
an {\it orbit equivalence} ({\it OE}) of $\Cal R, \Cal S$. More
generally:

\vskip .05in \noindent {\it 1.4.2. Definition}. Let $\Delta: (X,\mu)
\rightarrow (Y,\nu)$ be a measurable measure preserving map.
$\Delta$ is a {\it m.p. morphism} of $\Cal R$ into $\Cal S$ if there
exists $N_0 \subset X_0$ with $\mu(N_0)=0$ such that for all $t\in X
\setminus N_0$, $\Delta$ takes the $\Cal R$-orbits of $t$ into the
$\Cal S$-orbits of $\Delta(t)$. $\Delta$ is called an {\it
embedding} of $\Cal R$ into $\Cal S$ if it is an isomorphism of $X$
onto $Y$ and takes almost every orbit orbit of $\Cal R$ into an
orbit of $\Cal S$. $\Delta$ is a {\it local OE} (or {\it local
isomorphism}) of $\Cal R$, $\Cal S$ if there exists $N_0 \subset X$,
$\mu(N_0)=0$, such that $\forall t\in X \setminus N_0$, $\Delta$ is
a bijection between the $\Cal R$-orbit of $t$ and the $\Cal S$-orbit
of $\Delta(t)$.

\vskip .05in

By $(1.1.1)$, since any morphism $\Delta$ is surjective (a.e.), i.e.
$\nu(Y \setminus\Delta(X))=0$, $\Delta^*:L^\infty Y \rightarrow
L^\infty X$ is a faithful integral preserving von Neumann algebra
embedding. It is trivial from the definitions that if $\Delta$ is an
isomorphism of probability spaces then $\Delta$ is an embedding on
$\Cal R$ into $\Cal S$ iff $\Delta [\Cal R] \Delta^{-1}\subset [\Cal
S]$ and iff $(\Delta^{-1})^*:L^\infty X \simeq L^\infty Y$ extends
to a von Neumann algebra isomorphism of $L(\Cal R)$ into $L(\Cal
S)$.

If $\Delta:(X,\mu) \rightarrow (Y,\nu)$ is a local OE of $\Cal
R,\Cal S$, then for each $\psi \in  {_{_{p}}[\Cal S]}$ let
$\Delta^*(\psi)$ be the {\it pull back} of $\psi$, i.e. the local
isomorphism from $\Delta^{-1}(r(\psi))$ onto $\Delta^{-1}(l(\psi))$
which takes $t$ onto the unique element $t' \in \Cal R t$ with
$\Delta(t')=\psi(\Delta(t))$. Thus, $\Delta^*(\psi)\in {_{_{p}}[\Cal
R]}$ satisfies $\Delta \circ \Delta^*(\psi)=\psi \circ \Delta$.
Another way of describing the pull back map is as follows: Let
$\phi_n \in [\Cal R]$ be an orthonormal basis of $\Cal R$ and denote
$X^\psi_n=\{t \in \Delta^{-1}(l(\psi)) \mid
\psi^{-1}(\Delta(t))=\Delta (\phi_n^{-1}(t))\}$; $\{X^\psi_n\}_n$
are then measurable, give a partition of $\Delta^{-1}(l(\psi))$ and
we have $\psi^{-1} (\Delta(t)) = \phi_n^{-1}(t), \forall t\in
X_n^\psi.$ The pull back is clearly multiplicative (by the
definitions), thus giving an embedding of full pseudogroups
$\Delta^*: {_{_{p}}[\Cal S]} \rightarrow {_{_{p}}[\Cal R]}$. By the
above definition of the von Neumann algebra associated to a full
pseudogroup, this implies that $\Delta^*$ induces an isomorphism
from $L(\Cal S)$ into $L(\Cal R)$, still denoted $\Delta^*$. From
the above remarks, if we denote $v_\psi$, $\psi\in {_{_{p}}[\Cal
S]}$, the canonical partial isometry implementing $\psi$ on
$L^\infty Y$, i.e. $v_\psi a {v_\psi}^* = \psi(a)$, $a\in L^\infty
Y$, and we let $p^\psi_n = \chi_{X^\psi_n} \in L^\infty X$, then
with the above notations we trivially have:

\proclaim {1.4.3. Proposition} $\Delta^*(v_\psi)=\Sigma_n p^\psi_n
u_{\phi_n}$ defines a von Neumann algebra embedding of $L(\Cal S)$
into $L(\Cal R)$. Moreover, if we identify $L(\Cal S)$ with its
image under $\Delta^*$, then we have $(L^\infty Y)'\cap L(\Cal
R)=L^\infty X$, any element in $L(\Cal S)$ normalizing $L^\infty Y$
normalizes $L^\infty X$ and any orthonormal basis $v_n$ of $L(\Cal
S)$ over $L^\infty Y$ is an orthonormal basis of $L(\Cal R)$ over
$L^\infty X$.
\endproclaim

\vskip .05in

\noindent {\bf 1.5. von Neumann subalgebras and basic construction}.
If $(Q,\tau)$ is a finite von Neumann algebra then a $^*$-subalgebra
$N\subset Q$ closed in the weak operator topology (equivalently,
$(N)_1$ closed in $\|\cdot \|_2$) and with the same unit as $Q$ is
called a {\it von Neumann subalgebra} of $Q$. For instance, if $Q =
N \rtimes_\sigma \Gamma$ is the crossed product algebra
corresponding to some actions $\sigma$ of a discrete group $\Gamma$
on the finite von Neumann algebra $(N,\tau)$ as in 1.3, then $N$
identifies naturally with a von Neumann subalgebra of $Q$ by viewing
$a\in N$ as the ``polynomial'' $au_e$. Another important subalgebra
of $N \rtimes \Gamma$ is the von Neumann subalgebra $L\Gamma$
generated by the canonical unitaries $\{u_g\}_g \subset L^\infty X
\rtimes \Gamma$, i.e. the algebra of all convolvers $\Sigma_g c_g
u_g$ with scalar coefficients $c_g\in \Bbb C$.

The restriction of functionals from $Q$ to $N$ implements a positive
$N$-bimodular projection of $L^1(Q,\tau)$ onto $L^1(N,\tau_{|N})$
(Radon-Nykodim type theorem), whose restriction to $Q$ gives the
(unique) $\tau$-preserving conditional expectation of $Q$ onto $N$,
denoted $E_N$. Restricted to $L^2 Q$ it implements the orthogonal
projection of $L^2 Q$ onto $L^2 N$, denoted $e_N$. Identifying
$Q=L_Q\subset \Cal B(L^2 Q)$ gives $e_N x e_N = E_N(x)e_N, x\in Q$.

We denote by $\langle Q, e_N \rangle$ the von Neumann algebra
generated in $\Cal B (L^2 Q)$ by $Q=L_Q$ and $e_N$. Since $e_Nxe_N =
E_N(x)e_N, \forall x\in Q$, and $\vee \{x(e_N(L^2Q)) \mid x\in
Q\}=L^2 Q$, it follows that span$Qe_NQ$ is a *-algebra with support
equal to $1$ in $\Cal B(L^2 Q)$ (i.e. if $p\in \Cal B(L^2 Q)$ is a
projection with $pT=T$, $\forall T \in Qe_NQ$ then $p=1$). Thus,
$\langle Q, e_N \rangle= {\overline{\text{\rm sp}}}^{\text{\rm w}}
\{xe_Ny \mid x, y \in Q\}$. Also $e_N \langle Q, e_N \rangle e_N =
Ne_N$ implying that $\langle Q, e_N \rangle$ is a semifinite von
Neumann algebra. This is called the (Jones) {\it basic construction}
for the inclusion $N\subset Q$, with $e_N$ its {\it Jones
projection} ([J]).

We endow $\langle Q, e_N \rangle$ with a densely defined trace $Tr$
by $Tr(\Sigma_i x_ie_Ny_i)=\Sigma_i \tau(x_iy_i),$ for $x_i, y_i$
finite sets of elements in $Q$. We denote by $L^2(\langle Q, e_N
\rangle, Tr)$ the completion of sp$Qe_NQ$ in the norm $\|x\|_{2,Tr}=
Tr(x^*x)^{1/2}, x\in \text{\rm sp}Qe_NQ$. Exactly as in the case
of finite von Neumann algebras with a trace, $\langle Q,e_N \rangle$
acts (as a von Neumann algebra)
on $L^2(\langle Q,e_N \rangle, Tr)$ by left multiplication, and we call
this the {\it standard representation} of $(\langle Q,e_N \rangle, Tr)$.
Also, like in the finite case, the elements
in $L^2(\langle Q, e_N \rangle, Tr)$ can be viewed as square summable
(with respect to the semifinite trace $Tr$)
operators affiliated with $\langle Q, e_N \rangle$,  i.e. as
closed operators $T\in \Cal B(L^2(\langle Q, e_N \rangle, Tr))$
whose polar decomposition $T=u|T|$ has the partial isometry $u$
and the spectral resolution $e_s, s > 0,$ of $|T|$ lying
in $\langle Q, e_N \rangle$ and satisfying $Tr(T^*T)
\overset \text{\rm def} \to =
\int s^2 \text{\rm d} Tr(e_s) < \infty$.  The space
$L^1(\langle Q, e_N \rangle, Tr)$ is defined similarly
and for an operator $T$ affiliated with $\langle Q, e_N \rangle$
we have $T\in L^2(\langle Q, e_N \rangle, Tr)$ iff
$T^*T \in L^1(\langle Q, e_N \rangle, Tr)$, like in the finite case.

Any Hilbert subspace of $L^2 Q$ which is invariant under
multiplication to the right by elements in $N$ is a right Hilbert
$N$-module. If $\Cal H \subset L^2Q$ is a Hilbert subspace and $f$
is the orthogonal projection onto $\Cal H$ then $\Cal HN=\Cal H$
(i.e. $\Cal H$ is a right $N$-module) iff $f$ lies in $\langle Q,
e_N \rangle$. An {\it orthonormal basis} over $N$ for $\Cal H$ is a
subset $\{\eta_i\}_i \subset \Cal H$ such that $\Cal H =
\overline{\Sigma_k \eta_k N}$ and $E_N(\eta_i^*\eta_{i'}) =
\delta_{ii'}p_i \in \Cal P(N), \forall i, i'$. Note that in this
case we have $\xi = \Sigma_i \eta_i E_N(\eta_i^* \xi), \forall \xi
\in \Cal H$. A set $\{\eta_j\}_j \subset L^2 Q$ is an orthonormal
basis of $\Cal H$ over $N$ iff  the orthogonal projection $f$ of
$L^2 Q$ on $\Cal H$ satisfies $f = \Sigma_j \eta_j e_N \eta_j^*$,
with $\eta_je_N\eta_j^*$ projection $\forall j$. A simple maximality
argument shows that any left Hilbert $N$-module $\Cal H \subset L^2
Q$ has an orthonormal basis (see [P6] for all this).

If $\xi \in L^2Q$ satisfies $E_N(\xi^*\xi)\in N$ (i.e. it
is a bounded operator) then
the closed operator $\xi e_N \xi^*$ lies in $\langle Q, e_N \rangle$
(i.e. it is bounded).  If we denote by
$\Phi$ the $Q$-bimodule map from sp$Q
e_N Q \subset \langle Q, N \rangle $ into $Q$ defined by
$\Phi(xe_Ny) = xy, \forall x, y \in Q$, then $\tau\circ \Phi= Tr$
and $\Phi$ extends to a
linear map from $L^1(\langle Q, e_N \rangle, Tr)$
onto $L^1 Q$ satisfying
$\|\Phi (T)\|_1 \leq
\|T\|_{Tr,1}$.
If $\xi \in L^2 Q$ then $\Phi(\xi e_N \xi^*)=\xi \xi^*$.

An action $\sigma$ of a group $\Gamma$ on $(Q,\tau)$ which leaves
$N$ invariant extends to a $Tr$-preserving action $\sigma^N$ on
$\langle Q, e_N \rangle$ by $\sigma_g^N(xe_N y)=\sigma_g(x) e_N
\sigma_g(y),$ $\forall x,y\in Q$, $g\in \Gamma$. Moreover, since it
preserves $Tr$, $\sigma^N$ extends to a representation of $\Gamma$
on the Hilbert space $L^2(\langle Q, e_N \rangle, Tr)$.

\heading 2. Some generalities on cocycles
\endheading

\vskip .05in \noindent {\bf 2.1. Definition}.  Let $\sigma$ be an
action of a countable discrete group $\Gamma$ on a standard
probability space $(X,\mu)$ and  $\Cal V$ a Polish group. A (right)
measurable $1$-{\it cocycle} for $\sigma$ with values in $\Cal V$ is
a measurable map $w: X \times \Gamma \rightarrow \Cal V$ with the
property that for all $g_1, g_2 \in \Gamma$ the equation
$$
w(t, g_1) w({g_1}^{-1}t, g_2) = w(t, g_1g_2) \tag 2.1.1
$$
holds true $\mu$-almost everywhere in $t \in X$. Two $\Cal V$-valued
cocycles $w, w'$ are (measurably) {\it cohomologous} (or {\it
equivalent}) if there exists $u: X \rightarrow \Cal V$ measurable
such that for all $g\in \Gamma$ we have
$$
w'(t, g)=u(t)^{-1} w(t, g) u(g^{-1}t), \forall t\in X a.e.
\tag 2.1.2
$$
We then write $w'\sim w$. We denote by $\text{\rm Z}^1(\sigma;\Cal
V)$ the space of $\Cal V$-valued cocycles for $\sigma$ endowed with
the topology of convergence in measure. Note that if $\Cal V$ is
abelian then plain multiplication of cocycles (as $\Cal V$ valued
functions) gives a group structure on Z$^1(\sigma;\Cal V)$, which is
clearly Polish with respect to the above topology.

\vskip .05in

{\it Untwisting} a cocycle $w$ means showing it is equivalent to a
cocycle which is independent on $t \in X$ (in the a.e. sense). It is
immediate to see that such cocycles correspond precisely to group
morphisms of $\Gamma$ into $\Cal V$: If $\delta : \Gamma \rightarrow
\Cal V$ is a group morphism and we define $w^\delta : \Gamma \times
X \rightarrow \Cal V$ by $w^\delta(t,g)=\delta(g), \forall t$, then
$w^\delta$ is a cocycle; and conversely, if a cocycle $w: X \times
\Gamma \rightarrow \Cal V$ is so that for each $g\in \Gamma$ the map
$t \mapsto w(t,g)$ is constant in $t\in X$ (a.e.), then there exists
a unique group morphism $\delta : \Gamma \rightarrow \Cal V$ such
that $w=w^\delta$. We denote by Z$^1_0(\Gamma; \Cal V)$ the subset
of cocycles in Z$^1(\Gamma;\Cal V)$ which are equivalent to group
morphisms.

We only study in this paper cocycles with values in some closed
subgroup of the group $\Cal U(\Cal H)$ of unitaries acting on a
separable Hilbert space $\Cal H$, endowed with the usual Polish
group structure given by the $s^*$-topology. Such cocycles are of
particular interest because they can be interpreted in operator
algebra framework, as elements of the von Neumann algebra of bounded
measurable functions on $X$ with values in $\Cal U(\Cal H)$. Even
more so, such cocycles will become (left) cocycles for a certain
action of $\Gamma$ on a von Neumann algebra. To explain this in
details, we need some notations and general considerations.

Thus, if $(\Cal Y, d_{\Cal Y})$ is a separable complete metric space
with finite diameter and $(X,\mu)$ a standard probability space,
then we denote by $\Cal Y^X$ the set of classes (modulo null sets)
of measurable functions on $X$ with values in $\Cal Y$. We endow
this set with the metric given by the $L^2$-norm of the distance
function, i.e. $d(f_1,f_2)=(\int d_{\Cal Y}(f_1(t),f_2(t))^2
\text{\rm d} \mu (t))^{1/2}$, $f_1,f_2 \in \Cal Y^X)$. (N.B. Any
other $L^p$-norm, $p \geq 1$, gives an equivalent metric.) $(\Cal
Y^X, d)$ is clearly a complete metric space, which is separable
whenever $\Cal Y$ is separable. This follows immediately by
approximating elements in $\Cal Y^X$ by step functions. Since such
approximation will be used repeatedly in this paper, we mention it
as a lemma. Its proof is standard and is thus omitted:

\proclaim{2.2. Lemma} Let $(\Cal Y,d_{\Cal Y})$, $\Cal Y^X$ be as
above and $f_1, ..., f_m \in \Cal Y^X$.

$1^\circ$. For any $\varepsilon
> 0$ there exist a finite partition $X_0,X_1,...,X_n$ of
$X$ and elements $v_i^j \in \Cal Y$, $1\leq i\leq n, 1\leq j \leq
m$,  such that $\mu(X_0) \leq \varepsilon$ and $d_{\Cal Y}(f_j(t),
v^j_i) \leq \varepsilon$, $\forall 1\leq i \leq n, 1\leq j \leq m,$
$\forall t \in X_i$ a.e.

$2^\circ$. If $Y \subset X$ is a set of positive measure and
$(v_1,...,v_n)\in \Cal Y^n$ is an essential value of $(f_1,...,f_n)$
on $Y$ then there exists a decreasing sequence of subsets of
positive measure $Y \supset Y_1\supset...$ such that $\text{\rm
d}_{\Cal Y}(f_j(t),v_j)\leq 2^{-n}$, $\forall t \in Y_n$, $n \geq
1$, $1 \leq j \leq m$.
\endproclaim

Let now $\Cal B$ be a von Neumann algebra acting on the separable
Hilbert space $\Cal H$. Note that if $(\Cal B)_1$ denotes the unit
ball of $\Cal B$ (with respect to the operatorial norm $\| \cdot \|$
on $\Cal B$) then the Borel structures on $(\Cal B)_1$ corresponding
to the $w$, $s$ and $s^*$ topologies coincide (see e.g. [D1]). Thus,
if we equip  $(\Cal B)_1$ with either of these topologies and choose
a Borel structure on $X$, then the Borel functions on $X$ with
values in $(\Cal B)_1$ will be the same. We denote by
$L^\infty(X;\Cal B)$ the corresponding set of classes (modulo null
sets) of (essentially) bounded $\mu$-measurable functions on $X$ with
values in $\Cal B$. It has a natural $^*$-algebra structure given by
point addition, multiplication and $^*$-operation, and a C$^*$-norm
given by the ess-sup norm. This algebra acts in an obvious way on
the Hilbert space $L^2(X;\Cal H)$ of square integrable functions on
$X$ with values in $\Cal H$, as a von Neumann algebra.

There is a natural spatial isomorphism between the von Neumann
algebra $L^\infty(X;\Cal B)$ acting this way on $L^2(X;\Cal H)$ and
the tensor product von Neumann algebra $L^\infty(X,\mu)
\overline{\otimes} \Cal B$ acting on $L^2X \overline{\otimes} \Cal
H$, which sends a step function $f: X \rightarrow \Cal B$ taking
constant value $y_i \in \Cal B$ on $Y_i\subset X$, for $\{Y_i\}_i$ a
finite partition with measurable subsets of $X$, into the element
$\Sigma_i \chi_{Y_i} \otimes y_i$.

Notice now that any embedding of a separable Polish group $\Cal V$
as a closed subgroup of $\Cal U(\Cal B)$ implements an embedding of
$\Cal V^X$ as a closed subgroup of the unitary group of the von
Neumann algebra $L^\infty(X;\Cal B)$, and thus of the unitary group
of $L^\infty X \overline{\otimes} \Cal B$. By 2.2, when regarded as
a subgroup of $\Cal U(L^\infty X\overline{\otimes} \Cal B)$, $\Cal
V^X$ is the closure in the $s^*$-topology of the group of unitaries
of the form $\Sigma_i \chi_{Y_i} \otimes v_i$, with $v_i \in \Cal
V\subset \Cal U(\Cal B)$ and $\{Y_i\}_i$ finite partitions of $X$.

Notice also that if $\Cal V$ is countable discrete, then any
measurable map $f: X \rightarrow \Cal V$ is given by a partition
$\{Y_g\}_g$ of $X$ into a countable family of measurable sets such
that $f(t)=g, t\in Y_g, g\in \Cal V$. Equivalently, when regarded as
an element in $L^\infty X \overline{\otimes} \Cal B$, $f$ is of the form
$f=\Sigma_g \chi_{Y_g} \otimes v^0_g$, where $\{v^0_g\}_g = \Cal V
\subset \Cal U(\Cal B)$.

With $(X,\mu)$, $\Cal V\subset \Cal U(\Cal B)$ as above, let now
$\sigma: \Gamma \rightarrow \text{\rm Aut}(X,\mu)$ be an action of a
countable discrete group $\Gamma$ on $(X,\mu)$. We still denote by
$\sigma: \Gamma \rightarrow \text{\rm Aut}(L^\infty X, \tau_\mu)$
the action it implements on $L^\infty X$, as well as the $\Cal
B$-{\it amplification} of $\sigma$, i.e. the action of $\Gamma$ on
$L^\infty X \overline{\otimes} \Cal B$ given by $\sigma_g \otimes
id_{\Cal B}$, $g\in \Gamma$. If $w: X \times \Gamma \rightarrow \Cal
V$ is a measurable map and we denote $w_g=w(\cdot, g)\in \Cal V^X$,
$g\in \Gamma$, with $\Cal V^X$ viewed as a closed subgroup of the
unitary group of $L^\infty(X;\Cal B)=L^\infty X \overline{\otimes}
\Cal B$ as explained above, then conditions $(2.1.1), (2.1.2)$
translate into properties of $w_g$ as follows:

\proclaim{2.3. Lemma} A measurable map $w: X \times \Gamma
\rightarrow \Cal V$ is a cocycle for the action $\sigma$ of $\Gamma$
on $(X,\mu)$ $($i.e. it satisfies $(2.1.1))$ if and only if
$w_g=w(\cdot, g) \in \Cal V^X, g\in \Gamma,$ satisfies with respect
to the action $\sigma=\sigma \otimes id_{\Cal B}$ of $\Gamma$ on
$L^\infty X \overline{\otimes} \Cal B$ the relations

$$
w_g \sigma_g(w_h)=w_{gh}, \forall g,h\in \Gamma \tag 2.3.1
$$
Moreover, if $w,w'\in \text{\rm Z}^1(\sigma; \Cal V)$ and $u\in \Cal
V^X$ then $u$ satisfies $(2.1.2)$ if and only if, when viewed as an
element in $\Cal U(L^\infty X \overline{\otimes} \Cal B)$, it
satisfies the relations
$$
u^* w_g \sigma_g(u)=w'_g, \forall g\in \Gamma \tag 2.3.2
$$
\endproclaim
\noindent {\it Proof}. Trivial by the definitions. \hfill
$\square$

\vskip .05in

\noindent {\bf 2.4. Definition}. Let $\Gamma$ be a discrete group,
$\Cal N$ a von Neumann algebra and $\sigma : \Gamma \rightarrow
\text{\rm Aut}(\Cal N)$ an action of $\Gamma$ on $\Cal N$ (i.e. a
group morphism of $\Gamma$ into the group of automorphisms
$\text{\rm Aut}(\Cal N)$ of the von Neumann algebra $\Cal N$). A
{\it cocycle} for $\sigma$ is a map $w: \Gamma \rightarrow \Cal
U(\Cal N)$ satisfying equations $(2.3.1)$ above. Two cocycles $w,
w'$ for $\sigma$ are {\it cohomologous} (or {\it equivalent}) if
there exists $u\in \Cal U(\Cal N)$ such that $(2.3.2)$ above holds
true. More generally, a {\it local cocycle} for $\sigma$ is a map
$w$ on $\Gamma$ with values in the set of partial isometries of
$\Cal N$ satisfying $w_g \sigma_g(w_h)=w_{gh}, \forall g,h\in
\Gamma$. It is easy to see that if $w$ is such a local cocycle, then
$p=w_e$ is a projection and for each $g\in \Gamma$ the element $w_g$
belongs to the set $\Cal U(p\Cal N\sigma_g(p))$ of partial
isometries with left support $p$ and right support $\sigma_g(p)$.

\vskip .05in With this terminology, Lemma 2.3 shows that a
measurable map $w: X \times \Gamma \rightarrow \Cal U(\Cal B)$ is a
cocycle for an action $\sigma$ of $\Gamma$ on the probability space
$(X,\mu)$ if and only if $w_g=w(\cdot, g), g\in \Gamma,$ is a
cocycle for the amplified action $\sigma_g \otimes id_{\Cal B}$ of
$\Gamma$ on the algebra $L^\infty(X;\Cal B)=L^\infty
X\overline{\otimes} \Cal B$. Also, equivalence of measurable
cocycles $w,w'$ for $\sigma$ corresponds to their equivalence as
algebra cocycles for $\sigma \otimes id_{\Cal B}$.

In the rest of the paper, we will in fact concentrate on measurable
cocycles with values in the following class of Polish groups:

\vskip .05in
\noindent {\bf 2.5. Definition}. A Polish group $\Cal
V$ is of {\it finite type} if it is isomorphic (as a Polish group)
to a closed subgroup of the group of unitary elements of a separable
(equivalently countably generated) finite von Neumann algebra. We
denote by $\mycal U_{fin}$ the class of all such groups.

\proclaim{2.6. Lemma} Let $\Cal V$ be a Polish group. The following
conditions are equivalent:

\vskip .05in $(i)$. $\Cal V\in \mycal U_{fin}.$

\vskip .05in $(ii)$. $\Cal V$ is isomorphic to a closed subgroup
of the group of unitary elements $\Cal U(Q)$ of a separable type
$\text{\rm II}_1$ factor $Q$.

\vskip .05in $(iii)$. $\Cal V$ is isomorphic to a closed subgroup
$\Cal V\subset \Cal U(\Cal H)$ of the group of unitary operators
on a separable Hilbert space $\Cal H$ such that there exists $\xi
\in \Cal H$, $\xi \neq 0$, with the properties:
$$
\langle u_1u_2 \xi, \xi \rangle=\langle u_2u_1 \xi, \xi \rangle,
\forall u_1, u_2 \in \Cal V, \tag a
$$
$$
\overline{\text{\rm sp}} \Cal V' \xi = \Cal H, \tag b
$$
where $\Cal V'$ denotes the commutant of $\Cal V$ in $\Cal B(\Cal
H)$.

\endproclaim
\noindent {\it Proof}. $(iii) \Rightarrow (i)$. If $\Cal V \subset
\Cal U(\Cal H)$ is a closed subgroup and there exists a vector $\xi
\in \Cal H$ such that conditions $(a)$, $(b)$ are satisfied, then
let $B=\Cal V'' \subset \Cal B(\Cal H)$ be the von Neumann algebra
generated by $\Cal V$ in $\Cal B(\Cal H)$. By $(a)$, $\tau (x)=
\langle x \xi,\xi \rangle$, $x\in B$, is a normal trace on $B$ which
by $(b)$ is faithful. Since the Hilbert space $\Cal H$ on which it
acts is separable, $B$ follows finite and separable.

$(i)\Rightarrow (ii)$. If $\Cal V$ is a closed subgroup of $\Cal
U(B)$ for a separable finite von Neumann algebra $(B,\tau)$, then
let $Q$ be the free product $Q=B
* R$, with $R$ the hyperfinite II$_1$ factor. By ([P7]) $Q$ is a
(separable) II$_1$ factor and since $\Cal V$ is closed in $\Cal
U(B)$, it is closed in $\Cal U(Q)\supset \Cal U(B)$ as well.

$(ii)\Rightarrow (iii)$. If $\Cal V \subset \Cal U(Q)$ for a
separable II$_1$ factor $(Q,\tau)$ and $\Cal H = L^2(Q,\tau)$ is
the standard representation of $Q$, then $\xi= 1$ is a trace
vector for $\Cal U(Q)$, thus for $\Cal V\subset \Cal U(Q)$ as
well, and one has $\overline{\text{\rm sp}} \Cal V' \xi \supset
\overline{\text{\rm sp}}Q' \xi =\Cal H$, showing that both
conditions $(a), (b)$ of $(i)$ are satisfied. \hfill $\square$

\vskip .05in The next result provides some easy examples:

\proclaim{2.7. Lemma} $1^\circ$. If a group $\Cal V$ is either
countable discrete or separable compact, then $\Cal V\in \mycal
U_{fin}$.

$2^\circ$. If $\Cal V_n \in \mycal U_{fin}, n \geq 1$, then $\Pi_{n
\geq 1} \Cal V_n\in \mycal U_{fin}$.

$3^\circ$. If $\Cal V \in \mycal U_{fin}$ and $(X,\mu)$ is a
standard probability space then $\Cal V^X\in \mycal U_{fin}.$
\endproclaim
\noindent {\it Proof}. $1^\circ$. In both cases ($\Cal V$ discrete,
or $\Cal V$ separable compact), the group von Neumann algebra $L\Cal
V$ of $\Cal V$ is finite and has a normal faithful trace state (see
e.g. [D]). Since $\Cal V$  embeds into the unitary group of $L\Cal
V$ as the (closed) group of canonical unitaries $\{u_g\}_{g\in \Cal
V}$, it is of finite type by 2.6.

$2^\circ$. If $\Cal V_n$ is a closed subgroup in $\Cal U(B_n)$, for
some finite von Neumann algebra $B_n$ with normal faithful trace
state $\tau_n$, $n \geq 1$, and we let $(B,\tau)$ be given by $B =
\oplus_n B_n$, $\tau (\oplus_n x_n)=\Sigma_n \tau(x_n)/2^n$, then
$(B,\tau)$ is finite with normal faithful trace and $\Pi_n \Cal V_n$
embeds into $\Cal U(B)$ as the closed subgroup $\Pi_n \Cal V_n$.

Part $3^\circ$ is trivial, since $\Cal V$ closed in $\Cal U(B)$
implies $\Cal V^X$ closed in the unitary group of the (separable)
finite von Neumann algebra $L^\infty(X,\mu) \overline{\otimes} B.$
 \hfill $\square$

\vskip .05in \noindent {\bf 2.8. Notation}. If $\sigma$ is an action
of a discrete group $\Gamma$ on a finite von Neumann algebra
$(Q,\tau)$ then we denote by Z$^1(\sigma)$ the set of cocycles for
$\sigma$. Note that if $N\subset Q$ is a von Neumann algebra such
that $\sigma_g(N)=N, \forall g\in \Gamma$, and we denote
$\rho_g=\sigma_{g|N}$ the restriction of $\sigma$ to $N$, then we
have a natural embedding Z$^1(\rho) \subset \text{\rm Z}^1(\sigma)$.
This embedding is in general not compatible with the equivalence of
cocycles, i.e. two cocycles in Z$^1(\rho)$ may be equivalent as
cocycles in Z$^1(\sigma)$ without being equivalent in Z$^1(\rho)$.
However, mixing properties of $\sigma$ can entail the compatibility
of the two equivalence relations. The following property from ([P2])
is quite relevant in this respect:

\vskip .05in \noindent {\bf 2.9. Definition}. Let $(Q,\tau)$ be a
finite von Neumann algebra and $N\subset Q$ a von Neumann
subalgebra. Let $\sigma: \Gamma \rightarrow \text{\rm Aut}(Q,\tau)$
be an action of a discrete group $\Gamma$ on $(Q,\tau)$ that leaves
$N$ globally invariant. The action $\sigma$ is {\it weak mixing}
(resp. {\it mixing}) {\it relative to} $N$ if for any finite set $F
\subset Q \ominus N$ and any $\varepsilon > 0$ there exists $g \in
\Gamma$ (resp. there exists $K_0 \subset \Gamma$ finite) such that
$\|E_N(\eta^* \sigma_g(\eta'))\|_2 \leq \varepsilon, \forall \eta,
\eta' \in F$ (resp. $\forall g \in \Gamma \setminus K_0$).

\vskip .05in

It is easy to check that if $N \subset Q$ are abelian von Neumann
algebras then 2.9 becomes the weak mixing property introduced in the
1970's by Furstenberg [F] and Zimmer [Z3]. (N.B. I am grateful to
Alex Furman for drawing my attention on this work.)

To give alternative characterizations of property 2.9, recall some
notations from 1.4. Thus, let $N \subset Q \subset \langle Q, e
\rangle$ be the basic construction for $N\subset Q$, with $e=e_N$
denoting the Jones projection. Let $Tr$ be the canonical trace on
$\langle Q, e\rangle$ and $\sigma^N$ the action of $\Gamma$ on
$(\langle Q, e \rangle, Tr)$ given by
$\sigma_g^N(xey)=\sigma_g(x)e\sigma_g(y)$, $x,y\in Q$.

The next two results and their proofs are in the spirit of (Sec. 3
and 5.2 in [P2]).

\proclaim{2.10. Lemma} The following conditions are equivalent:

\vskip .03in \noindent $(i)$. $\sigma$ is weak mixing relative to
$N$.

\vskip .03in \noindent $(i')$. For any finite 
subset $F\subset L^2 Q \ominus L^2 N$, with
$E_N(\eta^*\eta)$ bounded $\forall \eta \in F$, and any $\varepsilon
> 0$ there exists $g \in \Gamma$
such that $\|E_N(\eta^* \sigma_g(\eta'))\|_2 \leq \varepsilon,
\forall \eta, \eta' \in F$.

\vskip .03in \noindent $(ii)$. There exist $g_n \rightarrow \infty$
in $\Gamma$ and an orthonormal basis $\{1\} \cup \{\xi_j\}_{j\geq
1}$ of $L^2 Q$ over $N$ such that $\underset n \rightarrow \infty
\to \lim \| E^{Q}_{N}(\xi^*_i \sigma_{g_n}(\xi_j))\|_2  = 0$,
$\forall i,j.$

\vskip .03in \noindent $(iii)$. Any $\xi\in L^2(\langle Q, e
\rangle, Tr)$ fixed by $\sigma^N$ lies in the subspace $L^2(e\langle
Q, e \rangle e, Tr)=L^2(Ne)$.

\vskip .05in

Moreover, if $(Q,\tau)=(P\overline{\otimes} N, \tau_P\otimes
\tau_N)$ and $\sigma$ leaves both $N=1\otimes N$ and $P=P\otimes 1$
globally invariant, then $\sigma$ is weak mixing (resp. mixing)
relative to $N$ if and only if $\sigma_{|P}$ is weak mixing (resp.
mixing). Also, if $\sigma_{|P}$ is weak mixing then any fixed point
of $\sigma$ lies in $N$.
\endproclaim
\noindent {\it Proof}. $(ii) \Rightarrow (iii)$. To prove that all
fixed points of $\sigma^N$ as an action on $L^2(\langle Q, e\rangle,
Tr)$ are under $e$ it is sufficient to prove that any projection $f
\in \langle Q, e \rangle$ with $Tr(f)<\infty$ which is fixed by
$\sigma^N$ must satisfy $f \leq e$. Assume there exists a
$\sigma^N$-invariant $f$ with $fe \neq f$, equivalently
$(1-e)f(1-e)\neq 0$. Since $(1-e)$ is $\sigma^N$-invariant, by
replacing $f$ by an appropriate spectral projection of
$(1-e)f(1-e)\neq 0$ we may thus assume $0\neq f \leq 1-e$.

Denote $f_n = \Sigma_{i=1}^n \xi_i e \xi_i^*$. Then $Tr (f_n) <
\infty$, $\forall n,$ and $e+f_n \nearrow 1$ in $\langle Q,
e\rangle$. Thus $\|f_n f - f\|_{2,Tr} =\|(e+f_n)f-f\|_{2,Tr}
\rightarrow 0$.

Let $\varepsilon > 0$ and $n$ be large enough such that $\|f_n f -
f\|_{2,Tr} < \varepsilon/3$. Since $\sigma$ satisfies $(ii)$, there
exists $g\in \Gamma$ such that $\|E_N(\xi_j^*\sigma_g(\xi_i)) \|_2 <
\varepsilon/(3n), 1\leq i,j \leq n$. We have
$\sigma^N(f_n)=\Sigma_{i=1}^n \sigma_g(\xi_i) e \sigma_g(\xi_i^*)$
and thus
$$
\|\sigma^N_g(f_n)f_n\|^2_{2,Tr}=Tr(f_n \sigma^N_g(f_n))=
\Sigma_{i,j=1}^n Tr(\xi_j e \xi^*_j \sigma_g(\xi_i) e
\sigma_g(\xi_i^*))
$$
$$
=\Sigma_{i,j=1}^n \|E_N(\xi_j^*\sigma_g(\xi_i))\|_2^2 \leq
\varepsilon^2/9.
$$
This further implies
$$
Tr(f)=Tr(\sigma^N_g(f)f) = Tr(\sigma^N_g(f)f_nf) +
Tr(\sigma^N_g(f)(f-f_nf))
$$
$$
\leq |Tr(\sigma^N_g(f-ff_n)f_nf)|+ |Tr(\sigma^N_g(ff_n)f_nf)| +
|Tr(\sigma^N_g(f)(f-f_nf)|
$$
$$
\leq \|f-ff_n\|_{2,Tr} \|f_nf\|_{2,Tr} +
\|\sigma^N_g(f_n)f_n\|_{2,Tr} \|f\sigma^N_g(f)\|_{2,Tr} +
\|f-f_nf\|_{2,Tr} \|\sigma^N_g(f)\|_{2,Tr}
$$
$$
\leq (\varepsilon/3+ \varepsilon/3 + \varepsilon/3)
\|f\|_{2,Tr}=\varepsilon \|f\|_{2,Tr}.
$$
Taking $\varepsilon < 1$, this shows that $f=0$, a contradiction.
Thus, $f \leq e$, finishing the proof of $(ii) \Rightarrow (iii)$.

$(i') \Rightarrow (i)$ and $(i') \Rightarrow (ii)$ are trivial. To
prove $(i) \Rightarrow (i')$, note first that for $\eta \in L^2 Q$
the condition $E_N(\eta^*\eta)\in N$ (i.e. bounded) is equivalent to
$\eta e\eta^*$ bounded. Let $F\subset L^2 Q \setminus L^2 N$ be a
finite set such that $E_N(\eta^*\eta)\in N$, $\forall \eta \in F$.
By $(i)$ there exists $g_n \in \Gamma$ such that $\lim_n \|E_N(x^*
\sigma_{g_n}(y))\|_2=0$, for all $x, y$ in the (separable) von
Neumann algebra $Q_0$ generated by $\eta |\eta|^{-1}$ and by the
spectral projections of $\eta\eta^*$ for $\eta\in F$. Showing that
$\lim_n \|E_N(\xi^* \sigma_{g_n}(\eta))\|_2=0$, for some $\xi,
\eta\in F$, is equivalent to showing that $\lim_n \|\xi e \xi^*
\sigma_{g_n}^N(\eta e \eta^*)\|_{2,Tr}=0$. But since for spectral
projections $p$ of $\xi \xi^*$ we have $\|p \xi e\xi^* p\|\leq \|\xi
e\xi^*\|$, the latter is trivially implied by the fact that $\lim_n
\|p\xi e \xi^*p \sigma_{g_n}^N(q\eta e \eta^*q)\|_{2,Tr}=0$ for all
$\xi, \eta\in F$ and $p, q$ spectral projections of $\xi\xi^*,
\eta\eta^*$ corresponding to finite intervals.

$(iii) \Rightarrow (i')$. Assume by contradiction that there exists
a finite subset $F \subset L^2 Q \ominus L^2 N$, with
$E_N(\eta^*\eta)\in N, \forall \eta\in F,$ and $\varepsilon_0 > 0$
such that $\Sigma_{\eta, \eta'} \|E_N(\eta^* \sigma_g(\eta'))\|^2_2
\geq \varepsilon_0$, $\forall g\in \Gamma.$ This is equivalent to
the fact that $b = \Sigma_{\eta \in F} \eta e \eta^*$ satisfies
$$
Tr(b\sigma^N_g(b)) \geq \varepsilon_0, \forall g\in \Gamma.
$$

Note that $b \in L^2(\langle Q, e\rangle, Tr)$. Denote by
$K=\overline{\text{\rm co}}\{ \sigma^N_g(b) \mid g\in \Gamma
\}\subset L^2(\langle Q, e \rangle, Tr)$, the closure being in the
weak topology on the Hilbert space.  Let $b_0\in K$ be the unique
element of minimal norm $\|\cdot \|_{2,Tr}$. Noticing that
$\sigma^N_g(K)=K$ and $\|\sigma^N_g(b_0)\|_{2,Tr}=\|b_0\|_{2,Tr}$,
$\forall g$, by the uniqueness of $b_0$ it follows that
$\sigma^N_g(b_0)=b_0$, $\forall g\in \Gamma$. Thus $b_0 = eb_0e$.
But by the condition $F\perp L^2N$ it follows that $be=0$, thus
$0=Tr(beb_0)=Tr(bb_0) \geq \varepsilon_0$, a contradiction.

To prove that last part of the statement, notice first that if
$\sigma_{|P}$ is not weak mixing then it has an invariant finite
dimensional subspace $\Cal H_0 \subset L^2P\ominus \Bbb C$. If
$\{\eta_i\}_i$ is an orthonormal basis of $\Cal H_0$, then $\Sigma_i
\eta_i \otimes N$ is $\sigma$-invariant. Equivalently, $f = \Sigma_i
\eta_i e_N \eta_i^*\in \langle Q,e \rangle$ is fixed by $\sigma^N$.
It is also finite and satisfies $fe=0$, showing that $\sigma$ is not
weak mixing relative to $N$ (by the equivalence of $(i)$ and
$(iii)$).

Conversely, if $\sigma_{|P}$ is weak mixing then any orthonormal
basis $\{1\} \cup \{\xi_i\}$ of $L^2 P$ is an orthonormal basis of
$Q$ over $N$ and condition $(ii)$ above is clearly satisfied. By the
equivalence of $(ii)$ and $(i)$, it follows that $\sigma$ is weak
mixing relative to $N$.

Now if $\sigma_{|P}$ is weak mixing and $x\in P\overline{\otimes} N$
is fixed by $\sigma$ then $xex^*$ is fixed by $\sigma^N$ and, since
$\sigma$ is weak mixing relative to $N=1\otimes N$, by $(iii)$ it
follows that $xex^* \in Ne$, i.e. $x\in N$.
\hfill $\square$.

\vskip .05in

Note that if $N=\Bbb C$ then Definition 2.8 amounts to $\sigma$
being weak mixing and $(i) \Leftrightarrow (iii)$ amounts to the
characterization of this property stating that the only invariant
finite dimensional vector subspace of $L^2 N$ is $\Bbb C$.

\proclaim{2.11. Lemma} Let $\sigma, \rho$ be actions of the discrete
group $\Gamma$ on finite von Neumann algebras $(P,\tau), (N,\tau)$
and $w$ a cocycle for the diagonal product action $\sigma_g\otimes
\rho_g$, $g \in \Gamma$, on $(P\overline{\otimes} N,\tau)$. Let
$(P',\tau')$ be a finite von Neumann algebra containing $(P,\tau)$
and $\sigma'$ an extension of $\sigma$ to an action of $\Gamma$ on
$(P',\tau')$  such that $\sigma'$ is weak mixing relative to $P$.
If $w_g (\sigma'_g\otimes \rho_g)(v)=vw'_g, \forall g\in
\Gamma$, for some cocycle $w'$ with values in $\Cal U(N)\subset \Cal
U(P'\overline{\otimes} N)$ and some $v \in \Cal
U(P'\overline{\otimes} N)$, then $v\in P\overline{\otimes} N$, i.e.
$w\sim w'$ as cocycles for $\sigma_g\otimes \rho_g, g\in \Gamma$. In
particular, it follows that if $\Gamma \curvearrowright^\sigma
(X,\mu)$ is the quotient of an action $\Gamma
\curvearrowright^{\sigma'} (X',\mu')$ such that $\sigma'$ is weak
mixing relative to $L^\infty X$ and $\Cal V \subset \Cal U(N)$ is a
closed subgroup, then any $\Cal V$-valued cocycle $w$ for $\sigma$ which
can be untwisted to a group morphism $w':\Gamma \rightarrow \Cal V$
as a cocycle for $\sigma'$, with (un)twister $v: X' \rightarrow \Cal
V$, then $v$ comes from a m.p. map $v: X \rightarrow \Cal V$, and
thus $w$ can be untwisted as a cocycle for $\sigma$.
\endproclaim
\noindent {\it Proof}. If we consider the basic construction
$P\overline{\otimes} N \subset P'\overline{\otimes} N \subset
\langle P'\overline{\otimes} N, e \rangle$ and denote by $\lambda_g
= \sigma'_g \otimes (\text{\rm Ad}w'_g \circ \rho_g)$ then
$\lambda_g(v P\overline{\otimes} N)=v P\overline{\otimes} N$,
$\forall g$. But by 2.10, $\sigma'$ weak mixing relative to $P$
implies $\lambda$ weak mixing relative to $P\overline{\otimes} N$,
so by 2.10 again we have $v \in P\overline{\otimes} N$.  \hfill
$\square$

\vskip .05in We end this section with a result showing that the
equivalence class of a cocycle for an action of a group on a finite
von Neumann algebra is closed in the topology of uniform $\|\cdot
\|_2$-convergence. We also show that if two cocycles are
cohomologous then any partial isometry giving a ``partial
equivalence'' can be extended to a unitary element that implements
an actual equivalence.

\proclaim{2.12. Lemma} Let $w,w'$ be cocycles for the action
$\sigma$ of a group $\Gamma$ on a finite von Neumann algebra
$(Q,\tau)$.

$1^\circ$. If $\|w_g-w'_g\|_2 \leq \delta, \forall g\in \Gamma$,
then there exists a partial isometry $v\in Q$ such that $\|v-1\|_2
\leq 4\delta^{1/2}$ and $w_g \sigma_g(v) = v w'_g$, $\forall g\in
\Gamma$.

$2^\circ$. If for any $\varepsilon >0$ there exists $u\in \Cal U(Q)$
such that $\|w_g \sigma_g(u) - u w'_g\|_2\leq \varepsilon$, $\forall
g\in \Gamma$, then $w,w'$ are cohomologous.

$3^\circ$. If $w,w'$ are cohomologous and $v\in Q$ is a partial
isometry satisfying $w_g \sigma_g(v) = v w'_g$, $\forall g\in
\Gamma$, then there exists  $u\in \Cal U(Q)$ such that $uv^*v=v$ and
$w_g \sigma_g(u) = u w'_g$, $\forall g\in \Gamma$.
\endproclaim
\noindent {\it Proof}. 1$^\circ$. Note first that if we let
$\pi_g(\xi)=w'_g \sigma_g(\xi)w_g^*$, for each $g\in \Gamma$, $\xi
\in L^2 Q$, then $\pi_g$ are unitary elements and $g \mapsto \pi_g$
is a unitary representation of $\Gamma$ on $L^2 Q$. Indeed, this is
an immediate consequence of the cocycle relation $(2.4.1)$.

Let $K=\overline{\text{\rm co}}^w \{w'_gw_g^* \mid g\in \Gamma\}$,
the closure being in the weak topology on the Hilbert space $L^2 Q$.
Note that $K$ is actually contained in $Q$, more precisely $\|y\|
\leq 1$, $\forall y\in K$. In particular $K$ is bounded in the norm
$\|\cdot \|_2$, so it is compact in the weak topology. Also, since
$\|w'_gw_g^*-1\|_2=\|w'_g-w_g\|_2$, $\forall g\in \Gamma$, by taking
convex combinations and weak closure we get $\|y-1\|_2 \leq \delta$,
$\forall y\in K$. Let $y_0 \in K$ be the unique element of minimal
norm-$\|\cdot \|_2$. Noticing that $\pi_g(K)=K$ and
$\|\pi_g(\xi)\|_2=\|\xi\|_2$, $\forall g\in \Gamma, \xi \in L^2 Q$,
by the uniqueness of $y_0$ it follows that
$w'_g\sigma_g(y_0)w_g^*=\pi_g(y_0)=y_0, \forall g\in \Gamma$. Thus
$w'_g \sigma_g(y_0)=y_0w_g, \forall g\in \Gamma$. In addition,
$y_0\in K$ implies $\|y_0-1\|_2 \leq \delta$. But then the partial
isometry $v$ in the polar decomposition of $y_0$ also satisfies
$w'_g \sigma_g(v)=vw_g, \forall g\in \Gamma$, while by ([C2]) we
have $\|v-1\|_2 \leq 4\delta^{1/2}$.

$2^\circ$. For the  proof of this part and part $3^\circ$ below, we
use Connes' ``2 by 2 matrix trick'' ([C1]). Thus, let
$\tilde{\sigma}$ be the action of $\Gamma$ on $\tilde{Q}=M_{2 \times
2}(Q)=Q \otimes M_{2 \times 2}(\Bbb C)$ given by
$\tilde{\sigma}_g=\sigma_g \otimes id$. If $\{e_{ij} \mid 1\leq i,j
\leq 2\}$ is a matrix unit for $M_{2 \times 2} (\Bbb C)\subset
\tilde{Q}$, then $\tilde{w}_g=w_g e_{11}+ w'_g e_{22}$ is a cocycle
for $ \tilde{\sigma}$. If $B\subset \tilde{Q}$ denotes the fixed
point algebra of the action Ad$\tilde{w}_g \circ \tilde{\sigma}$,
then $e_{11}, e_{22}\in B$ and the existence of a unitary element
$u\in Q$ intertwining $w,w'$ is equivalent to the fact that $e_{11},
e_{22}$ are equivalent projections in $B$. Moreover, any partial
isometry $v\in Q$ satisfying $w_g \sigma_g(v) = v w'_g$, $\forall
g\in \Gamma$, gives a partial isometry $v \otimes e_{12}$ in $B$
with left, right supports given by $vv^* \otimes e_{11}$, $v^*v
\otimes e_{22}$.

Thus, by part $1^\circ$ and the hypothesis, there exist partial
isometries $v_n \in B$ such that $v_nv_n^* \leq e_{11}, v_n^*v_n
\leq e_{22}$ and $\|v_nv_n^* - e_{11}\|_2=\|v_n^*v_n - e_{22}\|_2
\rightarrow 0$. But this implies $e_{11}, e_{22}$ are equivalent in
$B$ (for instance, because they have the same central trace in $B$).

$3^\circ$. Since $w,w'$ are equivalent, the projections $e_{11},
e_{22}$ are equivalent in $B$. On the other hand, if $v$ satisfies
the given condition, then $\tilde{v}=v \otimes e_{12}\in B$ and
$\tilde{v}\tilde{v}^*\leq e_{11}, \tilde{v}^*\tilde{v} \leq e_{22}$.
This implies the projections $e_{11}-\tilde{v}\tilde{v}^*$ and
$e_{22}-\tilde{v}^*\tilde{v}$ are equivalent in $B$, say by a
partial isometry $\tilde{v}'=v'\otimes e_{12}$. Then $u=v'+v$
clearly satisfies the condition. \hfill $\square$

\heading 3. Techniques for untwisting cocycles
\endheading

Let $\sigma$ be an action of a countable discrete group $\Gamma$ on
the standard probability space $(X,\mu)$. Denote by $\tilde{\sigma}$
the associated {\it double action} of $\Gamma$ on $(X \times X, \mu
\times \mu)$, given by the diagonal product
$\tilde{\sigma}_g(t_1,t_2)=(gt_1, gt_2)$, $t_1, t_2\in X$.

The main result in this section is a criterion for untwisting
cocycles, extracted from proofs in ([P1]) and ([P2]).

\proclaim{3.1. Theorem} Assume $\sigma$ is weakly mixing and let
$\rho$ be another action of $\Gamma$ on a standard probability space
$(Y,\nu)$. Let $\Cal V \in \mycal U_{fin}$ and $w\in \text{\rm
Z}^1(\sigma \times \rho;\Cal V)$, where $\sigma\times \rho$ is the
diagonal product action on $(X \times Y, \mu \times \nu)$ given by
$\sigma_g \times \rho_g, g\in \Gamma$. Denote by $w^l$, $w^r$ the
cocycles for the diagonal product action $\tilde{\sigma}\times \rho$
defined by $w^l (t_1, t_2, s, g)=w(t_1, s, g)$, $w^r (t_1, t_2, s,
g)=w(t_2, s, g)$, $g\in \Gamma, t_1, t_2 \in X, s\in Y$. Assume
there exists some separable finite von Neumann algebra $(N,\tau)$
such that $\Cal V$ can be realized as a closed subgroup of $\Cal
U(N)$ and such that $w^l\sim w^r$ as $\Cal U(N)$-valued cocycles for
$\tilde{\sigma}\times \rho$. Then $w$ is equivalent to a $\Cal
V$-valued cocycle which is independent on the $X$-variable.

Conversely, if $w\in \text{\rm Z}^1(\sigma \times \rho;\Cal V)$ is
equivalent to a cocycle in $\text{\rm Z}^1(\rho;\Cal V)$ then $w^l
\sim w^r$ in $\text{\rm Z}^1(\tilde{\sigma}\times \rho;\Cal V)$.

In particular, if $w\in \text{\rm Z}^1(\sigma;\Cal V)$ then $w\in
\text{\rm Z}_0^1(\sigma;\Cal V)$ if and only if $w^l \sim w^r$ in
$\text{\rm Z}^1(\tilde{\sigma};\Cal V)$, where
$w^l(t_1,t_2,g)=w(t_1,g)$, $w^r(t_1,t_2,g)=w(t_2,g)$, $t_1,t_2\in
X$, $g\in \Gamma$.
\endproclaim

The above theorem holds in fact true under more general assumptions, in a
non-commutative setting. To state the result, we need some
notations. Thus, let $\sigma: \Gamma \rightarrow \text{\rm
Aut}(P,\tau_P)$ be an action of a discrete group $\Gamma$ on a
finite von Neumann algebra $(P,\tau_P)$. We assume there exists an
extension of $\sigma$ to an action $\tilde{\sigma}$ of $\Gamma$ to a
larger finite von Neumann algebra $(\tilde{P},\tilde{\tau})$, with
the property that there exists $\alpha_1\in \text{\rm
Aut}(\tilde{P}, \tilde{\tau})$ commuting with the action
$\tilde{\sigma}$ and satisfying the properties:
$$
\overline{\text{\rm sp}}^w P \alpha_1(P)=\tilde{P}, \tilde{\tau}(x
\alpha_1(y))=\tau(x)\tau(y), \forall x,y\in P. \tag 3.2.0.a
$$

Let $(N,\tau_N)$ be another finite von Neumann algebra and $\rho$ an
action of $\Gamma$ on $(N,\tau)$. Denote by $\tilde{\sigma}'$ the
``diagonal product'' action of $\Gamma$ on $(\tilde{P}
\overline{\otimes} N,\tau)$, given by
$\tilde{\sigma}'_g=\tilde{\sigma}_g \otimes \rho_g$, $g\in \Gamma$.
Notice that the restriction to $P\overline{\otimes} N$ of
$\tilde{\sigma}'$ gives the action $\sigma'_g=\sigma_g \otimes
\rho_g$, $g\in \Gamma$. Also, denote $M=P\overline{\otimes} N
\rtimes_{{\sigma}'} \Gamma$, $\tilde{M} =
\tilde{P}\overline{\otimes} N \rtimes_{\tilde{\sigma}'} \Gamma$,
with the corresponding canonical unitaries denoted $\{u_g\}_g\subset
M$ respectively $\{\tilde{u}_g\}_g\subset \tilde{M}$. Whenever
identifying $M$ with a subalgebra of $\tilde{M}$ in the natural way,
we identify $u_g=\tilde{u}_g$.

Note that $[\alpha_1,\tilde{\sigma}]=0$ implies $[\alpha_1,
\tilde{\sigma}']=0$ (as usual, we still denote by $\alpha_1$ the
amplifications $\alpha_1\otimes id_N$). Thus, $\alpha_1$ extends to
an automorphism $\alpha_1$ of $\tilde{M}$, by $\alpha_1((x \otimes
y) \tilde{u}_g) =(\alpha_1(x)\otimes y) \tilde{u}_g$, $\forall x\in
\tilde{P}, y\in N$, $g\in \Gamma$.

The inclusion $M \subset \tilde{M}$ and the automorphism $\alpha_1$
of $\tilde{M}$ implement a $M-M$ Hilbert bimodule structure on
$L^2(\tilde{M})$, by
$$
z_1 \cdot \xi \cdot z_2 = z_1\xi \alpha_1(z_2), \forall z_1, z_2 \in
M, \xi \in L^2(\tilde{M}),  \tag 3.2.0.b
$$
which restricted to $P\overline{\otimes} N$ implements a
$P\overline{\otimes} N$ bimodule structure on
$L^2(\tilde{P}\overline{\otimes} N)$. With these notations at hand,
we have:

\proclaim{3.2. Proposition} Let $w:\Gamma \rightarrow \Cal
U(P\overline{\otimes} N)$ be a $1$-cocycle for $\sigma'$, which we
also view as a cocycle for $\tilde{\sigma}'$. Assume the action
$\sigma$ of $\Gamma$ on $P$ is weakly mixing. If there exists $b\neq
0$ in $\tilde{P} \overline{\otimes} N$ satisfying
$bw_g=\alpha_1(w_g)\tilde{\sigma}'_g(b)$, $\forall g\in \Gamma$,
with the left support of $b$ under some $p_0\in \Cal
P(P\overline{\otimes} N)$, then there exist a non-zero partial
isometry $v_0 \in P\overline{\otimes} N$, with $v_0v_0^* \leq p_0$,
$p=v_0^*v_0\in N$, and a local cocycle $w'_g\in \Cal
U(pN\rho_g(p))$, $g\in \Gamma$, for $\rho$ such that $w_g
\sigma'_g(v_0)=v_0 w'_g, \forall g\in \Gamma$.

Moreover, $w \sim \alpha_1(w)$ as cocycles for $\tilde{\sigma}'$ if
and only if the cocycle $w$ for $\sigma'$ is equivalent to a cocycle
$w': \Gamma \rightarrow \Cal U(N)$ for $\rho$.
\endproclaim

To prove this result we need some further considerations. With the
notations introduced above, note that if $g \mapsto v_g \in \Cal
U(M)$ is a unitary representation of a group $\Gamma$ then $\xi
\mapsto v_g \cdot \xi \cdot v_g^*$, $\xi \in L^2(\tilde{M})$, $g\in
\Gamma$, gives a representation of $\Gamma$ on the Hilbert space
$L^2(\tilde{M})$. Hence, if $\{w_g\}_g \subset \Cal U(P
\overline{\otimes} N)$ is a 1-cocycle for $\sigma'$ and we apply
this observation to $v_g = w_g u_g=w_g\tilde{u}_g$ then $\xi \mapsto
w_g\text{\rm Ad}(\tilde{u}_g) (\xi) \alpha_1(w_g^*)$, $\xi \in
L^2(\tilde{M})$, gives a representation $\tilde{\sigma}'_w$ of
$\Gamma$ on $L^2(\tilde{M})$. This representation clearly leaves
$L^2(\tilde{P} \overline{\otimes} N)$ invariant, acting on it by
$$
\tilde{\sigma}'_w(g)(\eta)=w_g\tilde{\sigma}_g'(\eta)\alpha_1(w_g^*),
\eta \in L^2(\tilde{P} \overline{\otimes} N). \tag 3.2.1
$$

We will now identify the representation $\tilde{\sigma}'_w$ in a
different way. Namely, we let $N=1\otimes N \subset
P\overline{\otimes} N \overset e_0 \to \subset \langle
P\overline{\otimes} N, e_0 \rangle$ be the Jones basic construction
for the inclusion $N \subset P\overline{\otimes} N$, where $e_0$
denotes the Jones projection implementing the trace preserving
conditional expectation $E_{N}$ of $P\overline{\otimes} N$ onto $N$.
We then denote by $Tr$ the canonical trace on $\langle
P\overline{\otimes} N, e_0 \rangle$ and let ${\sigma'}^N$ be the
action of $\Gamma$ on $(\langle P\overline{\otimes} N, e_0 \rangle,
Tr)$ given by ${\sigma'}^N(g)(z_1 e_0 z_2) = \sigma'_g(z_1)e_0
\sigma'_g (z_2)$ (see 1.4). We also denote

$$
{\sigma'}^N_w(g)(z_1 e_0 z_2)=w_g \sigma'_g(z_1)e_0
\sigma'_g(z_2)w_g^*, \forall z_1, z_2 \in P\overline{\otimes} N, g
\in \Gamma \tag 3.2.2
$$

\proclaim{3.3. Lemma} Each ${\sigma'}^N_w(g), g\in \Gamma$ defines a
$Tr$-preserving automorphism of the semifinite von Neumann algebra
$\langle P\overline{\otimes} N, e_0 \rangle$. The map $g \mapsto
{\sigma'}^N_w(g)$ gives an action of $\Gamma$ on $(\langle
P\overline{\otimes} N, e_0 \rangle, Tr)$, thus also a representation
of $\Gamma$ on the Hilbert space $L^2(\langle P\overline{\otimes} N,
e_0 \rangle, Tr)$. Moreover, the map $x_1 \alpha_1(x_2)\otimes y
\mapsto (x_1 \otimes 1) (ye_0) (x_2\otimes 1)$ extends to an
isomorphism $\vartheta_0$ from the Hilbert space $L^2(\tilde{P}
\overline{\otimes} N )$ onto the Hilbert space $L^2(\langle
P\overline{\otimes} N, e_0 \rangle, Tr)$ which intertwines the
representations $\tilde{\sigma}'_w$, ${\sigma'}_w^N$. Also,
$\vartheta_0$ intertwines the left $P\overline{\otimes} N$-module
structures on these Hilbert spaces.
\endproclaim

While this lemma can be easily given a proof by direct computation,
we prefer to derive it as a consequence of a more general result,
which has the advantage of giving a clear conceptual explanation of
the equivalence of representations in 3.3. More precisely, we will
identify $\vartheta_0$ as the restriction (to a subspace) of an
isomorphism of Hilbert $M$-bimodules.

To do this, note that the von Neumann algebra generated in $M$ by
$N=1\otimes N$ and $\{u_g\}_g$ is equal to $N\rtimes \Gamma=N
\rtimes_\rho \Gamma$. Let $N \rtimes \Gamma \subset M \overset e \to
\subset \langle M, e\rangle$ be the basic construction corresponding
to the inclusion $N \rtimes \Gamma \subset M$, where $e$ denotes the
Jones projection implementing the trace preserving conditional
expectation of $M$ onto $N \rtimes \Gamma$. Endow $L^2(\langle M,
e\rangle, Tr)$ with the Hilbert $M$-bimodule structure implemented
by the inclusion $M \subset \langle M, e\rangle$. Then we have:

\proclaim{3.4. Lemma} The map $\vartheta$ which takes $(x_1
\alpha_1(x_2)\otimes y) \tilde{u}_g$ onto $x_1 (y u_g) e
\sigma_g^{-1}(x_2)$, $x_1,x_2 \in P, y\in N, g\in \Gamma$, extends
to a well defined isomorphism between the Hilbert $M$-bimodules
$_ML^2(\tilde{M})_M$ and $_ML^2(\langle M, e\rangle, Tr)_M$.
\endproclaim
\noindent {\it Proof}. Note first that $\vartheta$ is consistent
with the $M$-bimodule structures on the two Hilbert spaces. Thus, in
order to prove the statement it is sufficient to show that
$\vartheta$ preserves the scalar product. Let $x_1, x_2, x_1', x_2'
\in P$, $y,y'\in N$, $g,g'\in \Gamma$. By the definition of the
scalar product in $\tilde{M}$ and $(3.2.0.a)$, we have:
$$
\langle (x'_1 \alpha_1(x_2') \otimes y')\tilde{u}_{g'},
(x_1\alpha_1(x_2) \otimes y)\tilde{u}_{g} \rangle
=\tau(x_1^*x_1')\tau(x_2^*x_2')\tau(y^*y') \delta_{g',g}
$$
On the other hand, by the definition of the scalar product in
$L^2(\langle M, e \rangle, Tr)$ we have
$$
\langle \vartheta((x'_1\alpha_1(x_2') \otimes y')\tilde{u}_{g'}),
\vartheta((x_1 \alpha_1(x_2) \otimes y)\tilde{u}_{g}) \rangle
$$
$$
=\langle x_1' (y'  u_{g'}) e \sigma_{g'}^{-1}(x_2'), x_1 (y u_{g})
e \sigma_{g}^{-1}(x_2) \rangle
$$
$$
=\langle x_1' y' e (x_2'u_{g'}), x_1 y e (x_2u_{g}) \rangle =
Tr(E_{N\rtimes \Gamma}(y^*x_1^*x_1'y') e E_{N\rtimes \Gamma}
(x_2'u_{g'}u_g^* x_2^*))
$$
$$
=Tr (y^* E_{N\rtimes \Gamma}(x_1^*x_1')y'e E_{N\rtimes \Gamma}(x_2'
\sigma_{g'g^{-1}}(x_2^*)) u_{g'g^{-1}})
$$
$$
=\tau(x_1^*x_1') \tau (x_2'\sigma_{g'g^{-1}}(x_2^*)) Tr(y^*y'e
u_{g'g^{-1}})
$$
$$
=\tau(x_1^*x_1') \tau(y^*y') \tau
(x_2'\sigma_{g'g^{-1}}(x_2^*)) \delta_{g',g}.
$$
This shows that $\vartheta$ preserves the scalar product, thus
finishing the proof. \hfill $\square$

\vskip .05in \noindent {\it Proof of} 3.3. We use the notations of
3.4. Since the Jones projection $e$ implements the trace preserving
expectation of $M$ onto $N \rtimes_\rho \Gamma$, it also implements
the trace preserving expectation of $P\overline{\otimes} N$ onto
$N=1\otimes N$. Moreover, since sp$P (N \rtimes\Gamma)$ is dense in
$L^2(M)$, it follows that the weakly closed $^*$ algebra generated
by $P\overline{\otimes} N$ and $e$ in $\langle M, e\rangle$ has the
same unit as $\langle M, e\rangle$ and is isomorphic to the von
Neumann algebra $\langle P\overline{\otimes} N, e_0 \rangle$ of the
basic construction for $N \subset P\overline{\otimes} N$, with the
canonical trace $Tr$ on $\langle P\overline{\otimes} N, e_0 \rangle$
corresponding to the restriction to $\langle P\overline{\otimes} N,
e \rangle \overset \text{\rm def} \to = \overline{\text{\rm sp}}^w
(P\overline{\otimes} N) e (P\overline{\otimes} N)$ of the canonical
trace of $\langle M, e\rangle$. In other words we have a
non-degenerate commuting square of inclusions ([P6]):

$$
\CD N \rtimes \Gamma \subset M \overset{e}\to\subset \langle M, e
\rangle\\
\noalign{\vskip-1pt}\cup \quad\quad \cup\quad\quad\quad \cup\\
\noalign{\vskip-6pt} N \subset P \overline{\otimes} N
\overset{e}\to\subset
\langle P \overline{\otimes} N, e \rangle\\
\endCD
$$

Thus, we can view $L^2(\langle P\overline{\otimes} N, e\rangle, Tr)$
as a Hilbert subspace of $L^2(\langle M, e\rangle, Tr)$. By the
definitions, $\vartheta_0$ clearly coincides with the restriction of
$\vartheta$ to $\tilde{P}\overline{\otimes} N$. The statement
follows then by noticing that, by Lemma 3.4, for $g \in \Gamma$ and
$\xi\in \tilde{P}\overline{\otimes} N$ we have
$$
\vartheta_0(\tilde{\sigma}'_w(g)(\xi))=\vartheta_0(w_g\tilde{u}_g
\xi \tilde{u}_g^* \alpha_1(w_g)^*)
$$
$$
=\vartheta(w_g\tilde{u}_g \xi \alpha_1(\tilde{u}_g^*w_g^*))=w_gu_g
\vartheta(\xi) u_g^* w_g^*
$$
$$
= w_g u_g \vartheta_0(\xi) u_g^*w_g^* =
{\sigma'}^N_w(g)(\vartheta_0(\xi)).
$$
\hfill $\square$

\vskip .05in \noindent {\it Proof of} 3.2. By Lemma 3.3,
$\xi=\vartheta_0(b) \in L^2(\langle P\overline{\otimes }N, e_0\rangle,
Tr)$ is a nonzero element which is fixed by ${\sigma'}^N_w(g),
\forall g\in \Gamma$. Moreover, since it preserves the trace $Tr$,
the action ${\sigma'}^N_w$ of $\Gamma$ on the semifinite von Neumann
algebra $(\langle P\overline{\otimes} N, e_0 \rangle, Tr)$ induces
actions of $\Gamma$ on all $L^p$-spaces associated with $(\langle
P\overline{\otimes }N, e_0\rangle, Tr)$, being multiplicative whenever
a product of elements in such spaces is well defined. This implies
that if $\eta = \eta^*$ is say in $L^1(\langle P\overline{\otimes
}N, e_0\rangle, Tr)$ and is fixed by ${\sigma'}^N_w$ then all its
spectral projections are fixed by ${\sigma'}^N_w$.

From these remarks it follows that $0\neq \xi\xi^* \in L^1(\langle
P\overline{\otimes}N, e_0\rangle, Tr)$ is fixed by ${\sigma'}^N_w$,
thus any spectral projection $f$ of $\xi\xi^*$ corresponding to an
interval of the form $[c,\infty)$ is fixed by ${\sigma'}^N_w$.
Moreover, $Tr(f) < \infty$ whenever $c>0$, with $f \neq 0$ if $c>0$
is sufficiently small. Also, since $\vartheta_0$ is a left
$P\overline{\otimes} N$-module isomorphism, if $p_0\in P
\overline{\otimes} N\subset \tilde{P}\overline{\otimes} N$ is a
projection satisfying $p_0b=p_0$, then $p_0\xi=\xi$ and thus
$p_0f=f$.

Let $f=\Sigma_i m_ie_0m_i^*,$  for some $\{m_i\}_i \subset
L^2(P\overline{\otimes} N, \tau)$ satisfying
$E_N(m_i^*m_j)=\delta_{ij}p_j \in \Cal P(N), \forall i,j$ (see e.g.
[P6]). In other words, $\{m_i\}_i$ is an orthonormal basis of $\Cal
H = f(L^2(P\overline{\otimes} N, \tau))$. We next show that by
``cutting'' $f$ with an appropriate projection in $\Cal Z(N)^\rho$,
we may assume $\{m_i\}_i$ is a finite set.

Since $\tau(\Sigma_i m_im_i^*)=Tr f <\infty$, the operator valued
weight $\phi(xe_0y)=xy$ satisfies $\phi(f)=\Sigma_i m_im_i^* \in
L^1(P\overline{\otimes} N, \tau)$. Since
$$
\Sigma_i m_ie_0m_i^*=f = {\sigma'_w}^N(f)=\Sigma_i
w_g\sigma'_g(m_i)e_0\sigma'_g(m_i)^*w_g^*,
$$
by applying $\phi$ to the first and last terms we get $\text{\rm Ad}
w_g \circ \sigma'_g(\Sigma_i m_im_i^*)=\Sigma_i m_im_i^*$, $\forall
g$. Thus $\Sigma_i m_im_i^*$ is fixed by the action
$\lambda_g=\text{\rm Ad} w_g \circ \sigma'_g$ of $\Gamma$ on
$(P\overline{\otimes} N, \tau)$. Also, $\lambda$ leaves invariant
the centers of $P\overline{\otimes} N$ and $N$ (the latter because
$\sigma'$ leaves $N$ invariant and Ad$w_g$ commutes with $\Cal
Z(N)\subset \Cal Z(P\overline{\otimes} N) = Z(N)\overline{\otimes}
\Cal Z(N)$). Thus, $a=E_{\Cal Z(N)}(\Sigma_i m_im_i^*)$ is fixed by
$\lambda$ and so are all its spectral projections. Let $z$ be a
spectral projection of $a$ such that $0\neq \|za \|< \infty$.

Since $\Cal Z(N)\subset \Cal Z(\langle P\overline{\otimes} N, e_0
\rangle)$, $z$ commutes with $\{m_i\}_i$, $e$ and $f$. Thus, by
replacing $f$ with $zf$ we may assume the orthonormal basis
$\{m_i\}_i$ is so that $E_{\Cal Z(N)}(\Sigma_i m_im_i^*)$ is a
bounded operator. But $E_{\Cal Z(N)}=E_{\Cal Z(N)}\circ E_{\Cal
Z(P\overline{\otimes} N)}$, the latter being in fact the central
trace $Ctr$ on $P\overline{\otimes} N$. Since $Ctr(xy)=Ctr(yx)$, we
get $Ctr(\Sigma_i m_im_i^*)=Ctr(\Sigma_i m_i^*m_i)$ and thus
$E_{\Cal Z(N)}(\Sigma_i m^*_im_i) = E_{\Cal Z(N)}(\Sigma_i
m_im_i^*)$ is bounded. This implies $p_i=E_N(m_i^*m_i)\in \Cal P(N)$
satisfy $\Sigma_i Ctr_N(p_i)$ bounded, where $Ctr_N$ denotes the
central trace on $N$. Thus, by using the ``cutting and gluing''
procedures in ($1.1.4$ of [P6]) we can choose the orhonormal basis
$m_1, m_2, ...$ such that the central supports $z(p_i)$ of $p_i$ in
$N$ satisfy $p_i\geq z(p_{i+1}), \forall i$. But this implies that
the cardinality of the set $\{m_2, m_3, ...\}$ is majorized by
$\|\Sigma_i Ctr(p_i)\|< \infty$, thus showing that it is a finite
set.

Let $t$ be the cardinality of the finite orthonormal basis
$\{m_i\}_i$. Denote $N^t = M_{t \times t} (N)=N \otimes M_{t \times
t} (\Bbb C)$ and notice that by tensoring all inclusions $N \subset
P\overline{\otimes} N \subset \langle P\overline{\otimes} N, e_0
\rangle$ by $M_{t \times t} (\Bbb C)$ one gets the non-degenerate
commuting squares of inclusions

$$
\CD N^t\overset{E_N^t}\to\subset P \overline{\otimes} N^t
\overset{e_0}\to\subset \langle P \overline{\otimes} N^t, e_0
\rangle\\
\noalign{\vskip-1pt}\cup \quad\quad \cup\quad\quad\quad \cup\\
\noalign{\vskip-6pt} N \overset{E_N}\to\subset P \overline{\otimes}
N   \overset{e_0}\to\subset
\langle P \overline{\otimes} N, e_0 \rangle\\
\endCD
$$
with ${\sigma'}^N$  extending to an action ${\sigma'}^{N^t}$ of
$\Gamma$ on $\langle P \overline{\otimes} N^t, e_0 \rangle$, which
acts as $\rho\otimes id$ on $N\otimes M_{t\times t}(\Bbb C)=N^t$.
Similarly, we denote ${\sigma'_w}^{N^t}$ the action Ad$w_g \circ
{\sigma'}^{N^t}(g), g\in \Gamma$.

Let $f_0=fe_{11}$, where $\{e_{ij}\}_{i,j}$ are matrix units for
$M_{t \times t}(\Bbb C) \subset 1\otimes N^t\subset P
\overline{\otimes} N^t.$ Note that $f_0 \prec  e_0$ in $\langle P
\overline{\otimes} N^t, e_0 \rangle$, so there exists an element $m\in
L^2(P\overline{\otimes}N^t)$ such that $f_0=me_0m^*$. In particular,
$q=E_{N^t}(m^*m)$ is a projection in $N^t=1_P\otimes N^t$

Since $me_0m^*={\sigma'}^{N^t}_w(g)(me_0m^*)=w_g{\sigma'}_g(m)e_0
\sigma'_g(m^*)w_g^*$, it follows that $me_0N^t=w_g\sigma'_g(m)e_0N^t$,
$\forall g\in \Gamma$. Taking the operator valued weight of $\langle P\overline{\otimes} N^t, e_0 \rangle$ onto $P\overline{\otimes} N^t$,
it follows that $mN^t=w_g\sigma'_g(m)N^t$, $\forall g\in \Gamma$, as well.
Thus, for each $g\in \Gamma$ 
there exists $w'_g\in N^t$ such that $mw'_g=w_g \sigma'_g(m)$ 
($w'_g=E_{N^t}(m^*w_g\sigma'_g(m)$ will do), and we have
$$
\sigma'_g(m^*m) =\sigma'_g(m^*)w_g^*w_g\sigma'_g(m)=
{w'_g}^*m^*mw'_g, \forall g.
$$
Since $\sigma'_g$ commutes with $E_{N^t}$, by applying $E_{N^t}$ to
these equalities we get $\sigma'_g(q)= {w'_g}^*qw'_g$, $\forall g$,
showing that $qw'_g\rho_g(q)$ are in the set $\Cal U(qN^t\rho_g(q))$
of partial isometries with left support equal to $q$ and right
support equal to $\rho_g(q)$. Thus, by replacing $w'_g$ by $qw'_g
\rho(q)$ we may assume $w'_g\in \Cal U(qN^t\rho_g(q))$, $\forall g$.

Since $w_h\sigma'_h(m)=mw'_h$, $w_{gh}\sigma'_{gh}(m)=m w'_{gh}$ and
$w_{gh}=w_g \sigma'_g(w_h)$, we then get for all $g,h\in \Gamma$:
$$
mw'_{gh}=w_{gh}\sigma'_{gh}(m)=w_{gh}\sigma'_g(\sigma'_h(m))
$$
$$
=w_g \sigma'_g(w_h)\sigma'_g(\sigma'_h(m))=w_g \sigma'_g(w_h
\sigma'_h(m))= w_g\sigma'_g(mw'_h)
$$
$$
=w_g \sigma'_g(m)\rho_g(w'_h)=m w'_g\rho_g(w'_h)
$$
showing that $w'_{gh}=w'_g\rho_g(w'_h)$, i.e. $w'$ is a local
cocycle for $\rho$.

Note that if one denotes by $\pi_1$ (resp. $\pi_2$) the
representations of $\Gamma$ into the unitary group of $M =
P\overline{\otimes} N^t \rtimes \Gamma$ (resp. into the unitary
group of $qMq$) given by $\pi_1(h)= w_hu_h$, $h\in \Gamma$ (resp.
$\pi_2(g)= w'_gu_g$, $g\in \Gamma$), then the relation
$w_g\sigma_g(m)=mw'_g, \forall g,$ simply states that $\pi_1, \pi_2$
are intertwined by $m$. This implies that the partial isometry $v\in
P\overline{\otimes} N^t$ in the polar decomposition of $m$
intertwines these representations as well and that $v^*v$ commutes
with $\pi_2(\Gamma)$. Thus, $v^*v$ belongs to the subalgebra of
$P\overline{\otimes} qN^tq$ fixed by the action $\sigma_g\otimes
({\text{\rm Ad}}(w'_g)\circ \rho_g), g\in \Gamma,$ which since
$\sigma_g$ is weakly mixing follows equal to the fixed point algebra
$B_0$ of the action Ad$w'_g \circ \rho_g$ of $\Gamma$ on $qN^tq$.

On the other hand, since $e_{11}m=m$ it follows that $v=e_{11}v$,
implying that the central trace of $v^*v$ in $N^t$
is $\leq 1/t$. Thus, $v^*v$ is equivalent in $N^t$
to a projection of the form $q_0e_{11}\in N^t$ with $q_0 \in N=N
\otimes 1\subset N^t$. Let $u_0 \in \Cal U(N^t)$ be so that
$u_0q_0e_{11}u_0^*=v^*v$.

Denote  $v_1=vu_0$, $w'_1(g)=u_0^* w'_g\rho_g(u_0)$. We then have
for all $g\in \Gamma$:
$$
v_1 w'_1(g)=vu_0u_0^* w'_g\rho_g(u_0)=vw'_g\rho_g(u_0)
$$
$$
=w_g \sigma'_g(v)\rho_g(u_0) = w_g \sigma'_g(vu_0)=w_g
\sigma'_g(v_1).
$$
Also,  $v_1 = e_{11}v_1e_{11}$, $v_1^*v_1\in e_{11}N^te_{11}=Ne_{11}$
and $w'_1(g)\in
\Cal U(v_1^*v_1N^t\rho_g(v_1^*v_1))$. This implies there exists a
partial isometry $v_0 \in N$ and a cocycle $w_0': \Gamma \rightarrow
\Cal U(v_0^*v_0N\rho_g(v_0^*v_0))$ for $\rho$ such that
$v_1=v_0e_{11}$ and $w'_1(g)=w'_0(g)e_{11}, \forall g$. But then, by
the isomorphism $N \simeq Ne_{11}$, we clearly have $v_0 w'_0(g)=w_g
\sigma'_g(v_0), \forall g$. Finally, notice that by the definition
of $m$ and the fact that $p_0f=p_0$, we have $p_0m=p_0$, implying
that $v_0v_0^* \leq p_0$.

This proves the first part of 3.2. To prove the last part, we use a
maximality argument. Thus, denote by $\mycal W$ the set of pairs
$(v,w')$ with $v\in P \overline{\otimes} N$ partial isometry
satisfying $v^*v\in N=1\otimes N$ and $w': \Gamma \rightarrow \Cal
U(v^*vN\rho_g(v^*v))$ a local cocycle for $\rho$ such that
$vw'(g)=w_g \sigma'_g(v)$, $\forall g\in \Gamma$. We endow $\mycal
W$ with the order: $(v_0,w_0') \leq (v_1,w'_1)$ iff
$v_0=v_1v_0^*v_0$, $v^*vw_1'(g)=w_0'(g), \forall g\in \Gamma$.

$(\mycal W, \leq)$ is clearly inductively ordered, so let
$(v_0,w'_0)\in \mycal W$ be a maximal element. We want to prove that
$v_0$ is a unitary element (automatically implying that $w'_0$ is a
a cocycle for $\rho$). Assume $v_0$ is not a unitary element and let
$v=v_0\alpha_1(v_0)^*\in \tilde{P} \overline{\otimes} N$. Since
$\alpha_1$ acts as the identity on $N$, it follows that
$vv^*=v_0v_0^*$. Also, for $g\in \Gamma$ we have
$$
w_g \tilde{\sigma}'_g(v) = w_g \sigma'_g(v_0)
\alpha_1(\sigma'_g(v_0))^* = v_0w'_0(g) \alpha_1(\sigma'_g(v_0^*))
$$
$$
=v_0 \alpha_1(\sigma'_g((v_0w'_0(g^{-1}))^*))=v_0
\alpha_1(\sigma'_g(w_{g^{-1}} \sigma'_{g^{-1}}(v_0))^*))
$$
$$
=v_0\alpha_1(\sigma'_g(\sigma'_{g^{-1}}(v_0^*)\sigma'_g(w_{g^{-1}})^*))
=v_0\alpha_1(v_0^*) \alpha_1(w_g),
$$
where for the last equality we have used the fact that the cocycle
relation $w_g \sigma'_g(w_{g^{-1}})=w_e=1$ implies
$\sigma'_g(w_{g^{-1}})^*=w_g$. Thus, since $w$ and $\alpha_1(w)$ are
equivalent cocycles for the action $\tilde{\sigma}'$, by 2.9 it
follows that there exists a partial isometry $v'\in \tilde{P}
\overline{\otimes} N$ such that $w_g \tilde{\sigma}'_g(v')=v'
\alpha_1(w_g)$, $\forall g\in \Gamma$, and $v'{v'}^*=1-vv^*$,
${v'}^*v'=1-v^*v$. The first part then implies there exists a
non-zero partial isometry $v_1 \in P\overline{\otimes} N$, with left
support majorized by $1-v_0v_0^*$ and right support in $N$, such
that $w_g \sigma'_g(v_1)=v_1 w'_1(g)$, $\forall g\in \Gamma$, for
some local cocycle $w'_1: \Gamma \rightarrow \Cal
U(v_1^*v_1N\rho_g(v_1^*v_1))$ for $\rho$.

Thus, in the finite von Neumann algebra $\tilde{P}\overline{\otimes}
N$ we have the equivalence of projections $v_1^*v_1\sim v_1v_1^*
\prec 1-v_0v_0^*\sim 1-v_0^*v_0$. Since the first and last
projections lye in $N$ and satisfy $v_1^*v_1 \prec 1-v_0^*v_0$ in
$\tilde{P}\overline{\otimes} N$, they satisfy the same relation in
$N$ (for instance, because the central trace of
$\tilde{P}\overline{\otimes} N$ is the tensor product of the central
traces on $\tilde{P}$ and $N$). Thus, by multiplying $v_1$ to the
right with a partial isometry in $N$ that makes $v_1^*v_1$
equivalent to a projection in $(1-v_0^*v_0)N(1-v_0^*v_0)$ and
conjugating $w'_1$ appropriately, we may assume $v_1^*v_1 \leq
1-v_0^*v_0$. But then $(v_0+v_1, w'_0 \oplus w'_1)\in \mycal W$
strictly majorizes $(v_0, w'_0)$, contradicting the maximality of
$(v_0, w'_0)$.

\hfill $\square$

The next result shows that untwisting cocycles with values in Polish
groups of finite type is a hereditary property. The rather
elementary argument is reminiscent of the proof of (5.2 in [P3]).

\proclaim{3.5. Proposition} Let $\sigma$ be a weakly mixing action
of the discrete group $\Gamma$ on the standard probability space
$(X,\mu)$. Let $\Cal W\in \mycal U_{fin}$ and $\Cal V \subset \Cal
W$ a closed subgroup. If a $\Cal V$-valued cocycle for $\sigma$ can
be untwisted as a $\Cal W$-valued cocycle then it can be untwisted
as a $\Cal V$-valued cocycle. More generally, let $\rho$ be another
action of $\Gamma$ on a standard probability space $(Y,\nu)$. If a
$\Cal V$-valued cocycle $w$ for the diagonal product action $\sigma
\times \rho$ is equivalent to a $\Cal W$-valued cocycle for $\rho$
then it is equivalent to $\Cal V$-valued cocycle for $\rho$.
\endproclaim
\noindent {\it Proof}. Embed $\Cal W$ as a closed subgroup of the
unitary group of a separable finite von Neumann algebra $(N,\tau)$.
Let $w: \Gamma \rightarrow \Cal V^{X\times Y}$ be a $\Cal V$-valued
1-cocycle for $\sigma$ and assume there exist a cocycle $w': \Gamma
\rightarrow \Cal W^Y$ for $\rho$ and a measurable map $v: X\times Y
\rightarrow \Cal W$ such that for each $g\in \Gamma$ we have
$$
w(t,s, g)v(g^{-1}t, g^{-1}s)=v(t,s)w'(s,g), \forall (t,s)\in X
\times Y (a.e.) \tag 3.5.1
$$

By the approximation Lemma 2.2 applied to the measurable map $X \ni
t \mapsto v(t, \cdot)\in \Cal W^Y$, there exists a decreasing
sequence of subsets of positive measure $X \supset X_1 \supset ...$
and an element $u_0\in \Cal W^Y$ such that $\|v(t)-u_0\|_2 \leq
2^{-n}$, $\forall t\in X_n$, $\forall n\geq 1$, where $\|\cdot \|_2$
corresponds here to the Hilbert norm on $L^\infty Y
\overline{\otimes} N$ given by $\tau_\nu \otimes \tau$.

This shows that by replacing $v$, regarded as a function in $t\in
X$, by $v(t)u_0^*, t\in X$, and $w'(\cdot, g)$ by $u_0w'(\cdot,
g)\rho(u_0^*)$, $g\in \Gamma$, we may assume $v\in \Cal W^{X\times
Y}$ and $w'\in \text{\rm Z}^1(\rho; \Cal W)$ satisfy both $(3.5.1)$
and
$$
\|v(t)-1\|_2 \leq 2^{-n}, \forall t \in X_n, \forall n \geq 1, \tag
3.5.2
$$
for some decreasing sequence of subsets of positive measure $X_n$.
We'll show that, together with the weak-mixing assumption, this
entails $v(t)\in \Cal V^Y$, $\forall t \in X$ (a.e.) and $w'(\Gamma)
\subset \Cal V^Y$, where $v(t)=v(t, \cdot)$ and $w'(g)=w'(\cdot,
g)$.

Thus, let $h_n \rightarrow \infty$ in $\Gamma$ be so that $\lim_n
\mu(h_nZ\cap Z')=\mu(Z)\mu(Z')$, $\forall Z, Z'\subset X$
measurable. Note that this implies that for any $g,g'\in \Gamma$ the
sequences $\{gh_ng'\}_n, \{gh_n^{-1}g'\}_n$ satisfy the same
condition.

Fix $h\in \Gamma$. For each $m \geq 1$ let $n_m$ be so that
$X'_m=X_m \cap h_{n_m}X_m$ and $X''_m=X_m \cap (h_{n_m}h)
X_m$ have positive measure. Then for $t'$ in $X'_m$ we have
$\|v(t')-1\|_2\leq 2^{-m}$, $\|v(h_{n_m}^{-1} t')-1\|_2 \leq 2^{-m}$
while for $t''\in X''_m$ we have $\|v(t'')-1\|_2 \leq 2^{-m}$ and
$\|v((h_{n_m}h)^{-1} t'')-1\|_2 \leq 2^{-m}$. By applying $(3.5.1)$
first for $g=h_{n_m}$ and then for $g=h_{n_m}h$, if we denote
$w(t,g)=w(t, \cdot, g)\in \Cal V^Y$, this implies
$$
\|w(t', h_{n_m})- w'(h_{n_m})\|_2 \leq 2^{-m+1}, \forall t'\in X'_m
$$
and respectively
$$
\|w(t'', h_{n_m}h)- w'(h_{n_m}h) \|_2 \leq 2^{-m+1}, \forall t''\in
X''_m.
$$
But this entails
$$
\|\rho_{h_{n_m}}(w'(h))-w(t', h_{n_m})^*w(t'', h_{n_m}h)\|_2
$$
$$
=\|w(t', h_{n_m})\rho_{h_{n_m}}(w'(h))-w(t'', h_{n_m}h)\|_2
$$
$$
\leq \|w(t', h_{n_m})\rho_{h_{n_m}}(w'(h)) -
w'(h_{n_m})\rho_{h_{n_m}}(w'(h))\|_2
$$
$$
+ \|w'(h_{n_m}h)-w(t'',
h_{n_m}h)\|_2 \leq 2^{-m+2}
$$
showing that $w'(h)$ can be approximated arbitrarily well with
elements in $\Cal V$. Thus $w'(h)\in \Cal V$.

We still need to prove that $v(t)\in \Cal V^Y$, $\forall t\in X
(a.e.)$. By using the approximation Lemma 2.2, it is clearly
sufficient to prove that for any set $Z' \subset X$ of positive
measure and any $\varepsilon > 0$ there exists a subset of positive
measure $Z_0 \subset Z'$ such that d$(v(t), \Cal V^Y)< \varepsilon$,
$\forall t\in Z_0$. Let $h_n \in \Gamma$ be so that  $\lim_n
\mu(h_nZ\cap Z')=\mu(Z)\mu(Z')$, $\forall Z, Z'\subset X$
measurable, as before. By $(3.5.2)$ we can take $m$ sufficiently
large  such that $\|v(t)-1\|_2 \leq \varepsilon$, $\forall t \in
X_m$. Let $n=n(m)$ be so that $\mu(h_nX_m \cap Z') \neq 0$ and
denote $Z_0 = h_nX_m \cap Z'$. If $t_0$ belongs to $Z_0$ then
$t=h_n^{-1}t_0 \in X_n$ so if we apply $(3.5.1)$ to $g=h_n^{-1}$ we
get the following identity in $\Cal W^Y$:
$$
v(t_0, h_n \cdot) = w(t, \cdot, h_n^{-1})^{-1}v(t, \cdot)w'(\cdot,
h_n^{-1})
$$
Since both $w(t, \cdot, h_n^{-1})$ and $w'(\cdot, h_n^{-1})$ lie in
$\Cal V^Y$ and $v(t, \cdot)$ is $\varepsilon$-close to an element in
$\Cal V^Y$, it follows that $v(t_0, h_n\cdot)$ is
$\varepsilon$-close to an element in $\Cal V^Y$, with $t_0\in Z_0$
arbitrary. Thus, $v(t_0, \cdot) \in \Cal V^Y$, $\forall t_0 \in
Z_0$,  as well \hfill $\square$

\vskip .05in \noindent {\it Proof of Theorem 3.1}. Since $\Cal V \in
\mycal U_{fin}$, the group $\Cal V$ is isomorphic to a closed
subgroup of $\Cal U(N)$, for some separable II$_1$ factor $N$. By
Proposition 3.2 it follows that, as cocycles for $\sigma \times
\rho$, $w$ is equivalent to a $\Cal U(N)$-valued cocycle $w'$ for
$\rho$. But by Propositon 3.5, $w$ is then equivalent in Z$^1(\sigma
\times \rho; \Cal V)$ to a  cocycle lying in Z$^1(\rho;\Cal V)$.
\hfill $\square$

\proclaim{3.6. Proposition} Let $\sigma, \rho$ be  actions of the
discrete group $\Gamma$ on finite von Neumann algebras $(P,\tau),
(N,\tau)$ and $w$ a cocycle for $\sigma\otimes \rho$. Let $H\subset
\Gamma$ be an infinite subgroup and assume $w_h\in \Cal U(N), h\in
H$. Let $H'$ denote set set of all $g\in \Gamma$ for which $w_g\in
N$. Then we have:

$1^\circ$. $H'$ is a subgroup of $\Gamma$.

$2^\circ$. If $g\in \Gamma$ is so that $gHg^{-1} \cap H$ is infinite
and $\sigma$ is weak mixing on $gHg^{-1} \cap H$ then $g\in H'$.

$3^\circ$. If $H$ is normal in $\Gamma$ and $\sigma$ is weak mixing
on $H$ then $H'=\Gamma$.

$4^\circ$. If $\sigma$ is weak mixing on $\Gamma$ then two cocycles
$w_1, w_2: \Gamma \rightarrow \Cal U(N)$ for $\rho$ are equivalent
as cocycles for the diagonal product action $\sigma \otimes \rho$ if
and only if they are equivalent as cocycles for $\rho$.
\endproclaim
\noindent {\it Proof}. Part $1^\circ$ is trivial and $3^\circ$
follows from $2^\circ$.

To prove $2^\circ$ we apply Lemma 2.10 to the action $\lambda$ of
$H_0=H \cap gHg^{-1}$ on $P\overline{\otimes} N$ given by
$\lambda_{h_0} =\text{\rm Ad}(w_{h_0}) \circ (\sigma_{h_0} \otimes
\rho_{h_0}), h_0\in H_0$. To this end, let us note that if $h_0$ is
an arbitrary element in $H_0$ and we denote $h=g^{-1}h_0g$ then
$h\in H$ and $h_0=ghg^{-1}$. By using the cocycle relations for $w$
we then get:
$$
\lambda_{h_0}(w_g)= w_{h_0} (\sigma_{h_0}\otimes
\rho_{h_0})(w_g)w_{h_0}^*
$$
$$
=w_{h_0g} w_{h_0}^* = w_{gh} w_{h_0}^*= w_g (\sigma_g \otimes
\rho_g) (w_h)w_{h_0}^*=w_g \rho_g(w_h)w_{h_0}^*.
$$
Thus, $\lambda_{h_0}(w_g N) = w_gN$, $\forall h_0\in H_0$. Since the
action $\lambda$ of $H_0$ is weak mixing on $P=P \otimes 1$ and
leaves $N=1\otimes N$ globally invariant, by 2.10 it follows that
$w_g\in 1\otimes N$, i.e. $g\in H'$.

4$^\circ$. Assume there exists $v\in \Cal U(P \overline{\otimes} N)$
such that $v^* w_1(g) (\sigma_g\otimes \rho_g)(v)=w_2(g)$ for all
$g\in \Gamma$. Thus, Ad$w_1(g) ((\sigma_g\otimes \rho_g)(v))=v
w_2(g)w_1(g)^*$, where we have identified $w_i(g)$ with $1 \otimes
w_i(g)$. This shows that $\lambda_g = \text{\rm Ad}(1\otimes w_1(g))
\circ (\sigma_g\otimes \rho_g)$ leaves invariant the space
$v(1\otimes N)$. By 2.10 this implies $v\in 1\otimes N$.  \hfill
$\square$

\heading 4. Perturbation of cocycles and the use of property (T)
\endheading

In this section we prove that if a $\Gamma$-action $\sigma$ has good
deformation properties and the group satisfies a weak form of the
property (T), then any cocycle for $\sigma$ checks the untwisting
criterion 3.1. The proof uses the same ``deformation/rigidity''
argument as on  (page 304 of [P1]) and (page 31 of [P2]). \vskip
.05in
\noindent {\bf 4.1. Definition}. An inclusion of discrete
groups $H \subset \Gamma$ has the {\it relative property} (T) of
Kazhdan-Margulis (or is {\it rigid}) if the following holds true (cf
[Ma]; see also [dHV]):

\vskip .05in \noindent $(4.1.1)$. $\exists F_0 \subset \Gamma$
finite and $\varepsilon_0 > 0$ such that if $\pi$ is a unitary
representation of $\Gamma$ on a Hilbert space $\Cal H$ with a unit
vector  $\xi \in \Cal H$ satisfying $\|\pi_g(\xi)-\xi\| \leq
\varepsilon_0$, $\forall g\in F_0$, then $\exists 0\neq \xi_0 \in
\Cal H$ with $\pi_h(\xi_0)=\xi_0, \forall h\in H$. \vskip .05in

Note that in case $H=\Gamma$, $(4.1.1)$ amounts to the usual
property (T) of Kazhdan for $\Gamma$ ([K]). It is easy to see that
if $H=\Gamma$, more generally if $H$ is normal in $\Gamma$, then the
fixed vector $\xi_0$ in $(4.1.1)$ can be taken close to $\xi$. This
fact was shown in ([Jo]) to hold true for arbitrary (not necessarily
normal) inclusions $H \subset \Gamma$ with the relative property
(T). In turn, the characterization of the relative property (T) for
$H \subset \Gamma$ in ([Jo]) is easily seen to be equivalent to the
following property, more suitable for us here:

\vskip .05in \noindent $(4.1.2)$. $\forall \varepsilon>0$ there
exist a finite subset $F(\varepsilon)\subset \Gamma$ and
$\delta(\varepsilon)>0$ such that if $\pi:\Gamma \rightarrow \Cal
U(\Cal H)$ is a unitary representation of $\Gamma$ on the Hilbert
space $\Cal H$ with a unit vector $\xi\in\Cal H$ satisfying
$\|\pi(g)\xi-\xi\|<\delta(\varepsilon)$, $\forall g\in
F(\varepsilon)$, then $\|\pi_h(\xi)-\xi\|<\varepsilon$ , $\forall
h\in H$.

\vskip .05in
The relative property (T) will only be used through the
following:

\proclaim{4.2. Lemma} Let $\sigma$ be an action of a discrete group
$\Gamma$ on a finite von Neumann algebra $(Q,\tau)$. Assume $H
\subset \Gamma$ is a subgroup with the relative property $(\text{\rm
T})$. If $w$ is a cocycle for $\sigma$, then for any $\varepsilon >
0$ there exists a neighborhood $\Omega$ of $w$ in the space
$\text{\rm Z}^1(\sigma)$ of cocycles for $\sigma$ such that $\forall
w'\in \Omega$ $\exists v \in Q$ partial isometry satisfying $w'_h
\sigma_h(v)=vw_h, \forall h\in H$ and $\|v-1\|_2 \leq \varepsilon$.
Moreover, if the action $\text{\rm Ad}(w_h)\circ \sigma_h, h\in H,$
of $H$ on $(Q,\tau)$ is ergodic, then the restriction to $H$ of any
cocycle in $\Omega$ is cohomologous to $w_{|H}$ (as cocycles for
$\sigma_{|H}$).
\endproclaim
\noindent {\it Proof}. With the notations in 4.1, let
$F=F(\varepsilon)$ and $\delta=\delta(\varepsilon^2/4)$. Let $w' \in
\text{\rm Z}^1(\sigma)$ be so that $\|w_g-w'_g\|_2 \leq \delta$,
$\forall g\in F$. Let $\pi$ be the representation of $\Gamma$ on
$L^2Q$ defined by $\pi(g)(\eta) = w'_g \sigma_g(\eta)w^*_g$, $\eta
\in  L^2 Q$, $g\in \Gamma$.

Then we have
$$
\|\pi_g(1)-1\|_2=\|w'_g \sigma_g(1)w^*_g -1\|_2=\|w_g-w'_g\|_2 \leq
\delta, \forall g\in F.$$ Thus, $\|\pi_h(1) - 1\|_2 \leq
\varepsilon^2/4$. Letting $\xi_0 \in L^2Q$ be the element of minimal
norm $\|\cdot \|_2$ in $K=\overline{\text{\rm co}}^w \{\pi_h(1) \mid
h\in H\}$, it follows that $\|\xi_0-1\|_2 \leq \varepsilon^2/4$,
$\|\xi_0\| \leq 1$ (thus $\xi_0 \in Q$) and $w'_h \sigma_h(\xi_0)
w_h^*=\xi_0$, $\forall h\in H$. But then $w'_h \sigma_h(\xi_0)=\xi_0
w_h$, $\forall h\in H$, so if $v \in Q$ denotes the partial isometry
in the polar decomposition of $\xi_0$ then $w'_h \sigma_h(v)=v w_h$,
$\forall h\in H$, and by ([C2]) we have $\|v-1\|_2 \leq
\varepsilon$.

Thus, if we let $\Omega=\{w'\in \text{\rm Z}^1(\sigma)
\mid \|w'_g-w_g\|_2 \leq \delta, \forall g\in F\}$, then $\Omega$
satisfies the desired conditions.

Now, since the relation $w'_h \sigma_h(v)=v w_h$, $\forall g\in H$,
implies that $v^*v$ is in the fixed point algebra of the action
Ad$w_h \circ \sigma_h, h\in H$, of $H$ on $Q$,
if the cocycle $w$ is ergodic then $v^*v=1$
showing that $w'\sim w$ as cocycles for $\sigma_{|H}$.
\hfill $\square$

\vskip .05in \noindent {\bf 4.3. Definition}.  Let $ \sigma $ be an
action of a discrete group $\Gamma$ on a standard probability space
$(X,\mu)$. Denote by $(A,\tau)=(L^\infty(X,\mu), \int \cdot
\text{\rm d} \mu)$ the associated abelian von Neumann algebra and by
$\tilde{\sigma}: \Gamma \rightarrow \text{\rm
Aut}(A\overline{\otimes} A, \tau \otimes \tau)$ the diagonal product
action given by $\tilde{\sigma}_g = \sigma_g \otimes \sigma_g$,
$g\in \Gamma$. The action $\sigma$ is {\it malleable} if the flip
automorphism $\alpha' \in \text{\rm Aut} (A\overline{\otimes} A,
\tau \otimes \tau)$, defined by $\alpha'(a_1 \otimes a_2)=a_2\otimes
a_1$, $a_1, a_2 \in A$, is in the (path) connected component of the
identity in the commutant of $\tilde{\sigma}(\Gamma)$ in $\text{\rm
Aut} (A\overline{\otimes} A, \tau \otimes \tau)$ (cf [P4]; note that
this terminology is used in [P1,2,3] for a slightly stronger
condition).

The action $\sigma$ is {\it s-malleable}  if there exist a
continuous action $\alpha:\Bbb R\rightarrow \text{\rm Aut}
(A\overline{\otimes} A, \tau \otimes \tau)$ and a period 2
automorphism $\beta\in \text{\rm Aut}(A\overline{\otimes} A, \tau
\otimes \tau)$ satisfying
$$
[\alpha, \tilde{\sigma}]=0, \alpha_1(A \otimes 1)=1\otimes A, \tag
4.3.1
$$
$$
[\beta, \tilde{\sigma}]=0, \beta(a \otimes 1)=a\otimes 1, \forall
a\in A, \beta \alpha_t = \alpha_{-t}\beta, \forall t\in \Bbb R, \tag
4.3.2
$$

More generally, an action $\sigma$ of $\Gamma$ on a finite von
Neumann algebra $(P, \tau)$ is {\it s-malleable} if there exist a
continuous action $\alpha: \Bbb R \rightarrow \text{\rm
Aut}(P\overline{\otimes} P, \tau \otimes \tau)$ and a period 2
automorphism $\beta \in \text{\rm Aut}(P\overline{\otimes} P, \tau
\otimes \tau)$ such that if we denote by $\tilde{\sigma}$ the
diagonal product action $\tilde{\sigma}_g = \sigma_g \otimes
\sigma_g$ of $\Gamma$ on $(P\overline{\otimes} P, \tau \otimes
\tau)$ then we have:

\vskip .05in \noindent $(4.3.1')$. $\alpha$ commutes with
$\tilde{\sigma}$ and satisfies:

$$
\tau(x\alpha_1(x))=\tau(x)\tau(y), \forall x,y \in P;
\overline{\text{\rm sp}}^w P \alpha_1(P)=P\overline{\otimes} P.
$$

\vskip .05in \noindent
$(4.3.2')$. $\beta$ commutes with
$\tilde{\sigma}$ and satisfies
$$\beta(x \otimes 1)=x\otimes 1, \forall
x\in P, \beta\alpha_t =\alpha_{-t} \beta, \forall t\in \Bbb R.
$$

\vskip .05in \noindent {\it 4.4. Example}. Let $(X_0, \mu_0)$ be a
standard probability space. Let $\Gamma$ be a countable discrete
group and $K$ a countable set on which $\Gamma$ acts. Let
$(X,\mu)=\Pi_k (X_0, \mu_0)_k$ be the standard probability space
obtained as the product of identical copies $(X_0,\mu_0)_k$ of
$(X_0, \mu_0)$, $k\in K$. Let $\sigma : \Gamma \rightarrow
{\text{\rm Aut}}(X, \mu)$ be defined by
$\sigma(g)((x_k)_k)=(x'_{k})_k$, where $x'_{k}=x_{g^{-1}k}$. We call
$\sigma$ the $(X_0, \mu_0)$-{\it Bernoulli} $(\Gamma
\curvearrowright K)$-{\it action}. We generically refer to such
actions as {\it generalized Bernoulli actions}. In case $K=\Gamma$
and $\Gamma \curvearrowright \Gamma$ is the left multiplication, we
simply call $\sigma$ the $(X_0, \mu_0)$-{\it Bernoulli}
$\Gamma$-{\it action}.

Note that if we denote $(A_0, \tau_0)=(L^\infty X_0, \int \cdot
\text{\rm d}\mu_0)$, then the algebra $(L^\infty X, \int \cdot
\text{\rm d}\mu)$ coincides with $\overline{\otimes}_k
(A_0,\tau_0)$, with the action implemented by $\sigma$ on elements
of the form  $\otimes_k a_k \in \overline{\otimes}_k (A_0, \tau_0)$
being given by $\sigma_g(\otimes_k a_k)=\otimes_k a'_k$,
$a'_k=a_{g^{-1}k}, k\in K$, $g\in \Gamma$.

More generally, if $(N_0,\tau_0)$ is a finite von Neumann algebra
and we denote $(N,\tau)=\overline{\otimes}_k (N_0,\tau_0)_k$, then
for each $g\in \Gamma$, $\otimes_k a_k \in N$ we let
$\sigma_g(\otimes_k a_k)=\otimes_k a'_k$, where $a'_k=a_{g^{-1}k},
k\in K$. Then $\sigma$ is clearly an action of $\Gamma$ on
$(N,\tau)$ which we call the $(N_0, \tau_0)$-{\it Bernoulli}
$(\Gamma \curvearrowright K)$-{\it action}.

It was shown in ([P1], [P2]) that if $N_0$ is either abelian diffuse
or a tensor product of 2 by 2 matrix algebras then any
$(N_0,\tau_0)$-Bernoulli $\Gamma \curvearrowright K$ action is
s-malleable. We reprove these results below, for the sake of
completeness.

\proclaim{4.5. Lemma} Let $\Gamma$ be a countable discrete group
acting on a countable set $K$. Let $(N_0, \tau_0)$ be a finite von Neumann
algebra and $\sigma$ the $(N_0,
\tau_0)$-Bernoulli $\Gamma \curvearrowright K$ action on
$(N,\tau)=\overline{\otimes}_{k\in K} (N_0,\tau_0)_k$.

$1^\circ$. If either $(N_0, \tau_0)$ is diffuse $($i.e.
it has no non-zero minimal projections$)$ and for all $ g\neq
e$ there exists $k\in K$ such that $gk\neq k$, or if $(N_0, \tau_0)$
is arbitrary and for all $g \neq e$ the set $\{ k \in K \mid gk \neq
k \}$ is infinite, then $\sigma$ is a free action.

$2^\circ$. $\sigma$ is weakly mixing iff $\forall K_0 \subset K$
$\exists g \in \Gamma$ such that $gK_0 \cap K_0 = \emptyset$ and iff
any orbit of $\Gamma \curvearrowright K$ is infinite. Moreover, if
any of these conditions is satisfied then $\sigma$ is weakly mixing
relative to any subalgebra of the form $N^0=\overline{\otimes}_k
(N^0_0)_k \subset N$, with $N^0_0 \subset N_0$.

$3^\circ$. $\sigma$ is (strongly) mixing iff $\forall K_0 \subset K$
finite $\exists F \subset \Gamma$ finite such that $gK_0 \cap K_0 =
\emptyset$, $\forall g \in \Gamma \setminus F$, and iff
the stabilizer $\{h \in \Gamma \mid hk=k\}$ of any $k\in K$ is finite.

$4^\circ$. If $N_0$ is either abelian diffuse or a finite factor of
the form $\quad \quad (N_0,\tau_0)=\quad \overline{\otimes}_{l\in L}
(M_{2 \times 2}(\Bbb C), tr)_l$, for some set of indices $L$, then
$\sigma$ is s-malleable.
\endproclaim

\noindent {\it Proof}. 1$^\circ$, $3^\circ$ and the first part of
$2^\circ$ are well known (and easy exercises!). To prove the last
part of $2^\circ$, let $\{\eta_n \mid n \geq 0\} \subset L^2
(N_0,\tau_0)$ be an orthonormal basis over $N^0_0$, in the sense of
1.4, with $\eta_0=1$. Denote by $\{\xi_n\}_n \subset N$ the
(countable) set of elements of the form $\xi_n = \otimes
(\eta_{n_k})_k$, with $n_k\geq 0$ all but finitely many equal to
$0$, at least one being $\geq 1$. It is immediate to see that $\{1\}
\cup \{\xi_n\}_n$ is an orthonormal basis of $L^2N$ over $N^0$ that
checks condition $(2.2.1)$ with respect to the Bernoulli $\Gamma
\curvearrowright K$ action $\sigma$.

To prove $4^\circ$, note first that $N_0$ abelian diffuse implies
$N_0 \simeq L^\infty(\Bbb T, \lambda)$, with $\tau_0$ corresponding to
$\int \cdot \text{\rm d}\mu$. To construct $\alpha$ and $\beta$ it is
then sufficient to construct an action $\alpha_0: \Bbb R \rightarrow
\text{\rm Aut}(A_0 \overline{\otimes} A_0, \tau_0 \otimes \tau_0)$
and  a period 2 automorphism $\beta_0 \in \text{\rm
Aut}(A_0\overline{\otimes} A_0, \tau_0\otimes \tau_0)$ such that
$\alpha_0(A_0 \otimes 1) = 1 \otimes A_0$, $\beta_0(a_0 \otimes 1) =
a_0 \otimes 1$, $\forall a_0 \in A_0$, $\beta_0 \alpha_0(t) =
\alpha_0(-t)\beta_0$, $\forall t \in \Bbb R$. Indeed, because then
the product actions $\alpha_t = \otimes_k (\alpha_0(t))_k$, $\beta =
\otimes_k (\beta_0)_k$ will clearly satisfy conditions $2.1.1$,
$2.1.2$.

Since $A_0$ is diffuse, $(A_0, \tau_0) \simeq (L^\infty(\Bbb T,
\mu), \int \cdot \text{\rm d}\mu)$. Let $u,v$ be Haar unitaries
generating $A_0 \otimes 1$ and $1 \otimes A_0$. Then $u,v$ is a {\it
pair of Haar unitaries} for $(\tilde{A}_0, \tilde{\tau}_0)=(A_0
\overline{\otimes} A_0, \tau_0 \otimes \tau_0)$, i.e. $u,v$ generate
$\tilde{A}_0$ and $\tilde{\tau}_0(u^nv^m)=0$ for all $(n,m)\neq
(0,0)$, where $\tilde{\tau}_0=\tau_0 \otimes \tau_0$. Note that if
$u',v' \tilde{A}_0$ is any other pair of Haar unitaries, then there
exists a unique automorphism $\theta$ of $(\tilde{A}_0,
\tilde{\tau}_0)$ that takes $u$ to $u'$ and $v$ to $v'$, defined by
$\theta(u^nv^m)=(u')^n(v')^m, \forall n,m \in \Bbb Z$. Note also
that if $w\in \Cal U(\tilde{A}_0)$ belongs to an algebra which is
$\tilde{\tau}_0$-independent of the algebra generated by $u$ then
$wu$ is a Haar unitary.

With this in mind, let $h$ be the unique selfadjoint element in
$\tilde{A}_0$ with spectrum in $[0,2\pi]$ such that $e^{ih}=vu^*$.
From the above remark it follows that both $u'=e^{ith}u$ and
$v'=e^{ith}v$ are Haar unitaries in $\tilde{A}_0$. Also,
$v'{u'}^*=vu^*$ so the von Neumann algebra generated by $u',v'$
contains $h$, thus it also contains $u,v$, showing that $u',v'$
generate $\tilde{A}_0$, thus being a pair of Haar unitaries for
$\tilde{A}_0$. Thus, there exists a unique automorphism
$\alpha_0(t)$ of $\tilde{A}_0$ such that $\alpha_0(t)(u)=e^{ith}u$,
$\alpha_0(t)(v)=e^{ith}v$. By definition, we see that
$\alpha_0(t)(vu^*)=vu^*$, thus $\alpha_0(t)(h)=h$ as well. But this
implies that for all $t,t' \in \Bbb R$ we have
$$
\alpha_0(t')(\alpha_0(t)(u))=
\alpha_0(t')(e^{ith}u)=e^{ith}\alpha_0(t')(u)=e^{ith}e^{it'h}u=
\alpha_0(t'+t)(u).
$$
Similarly $\alpha_0(t')(\alpha_0(t)(v))=\alpha_0(t'+t)(v),$ showing
that $\alpha_0(t'+t)=\alpha_0(t')\alpha_0(t)$, $\forall t,t'$.

Further on, since $u^*, uv^*$ is a Haar pair for $\tilde{A}_0$,
there exists a unique automorphism $\beta_0$ of $(\tilde{A}_0,
\tilde{\tau}_0)$ such that $\beta_0(u)=u^*$, $\beta_0(vu^*)=vu^*$.
By definition, $\beta_0$ satisfies $\beta_0^2(u)=u$ and
$\beta_0^2(vu^*)=vu^*$, showing that $\beta_0^2=id$. Finally, since
$\beta_0(h)=h$, we have $\alpha_0(t)(\beta_0(u))= e^{-ith} u^*=
\beta_0(\alpha_0(-t)(u))$. Similarly
$\alpha_0(t)(\beta_0(v))=\beta_0(\alpha_0(-t)(v))$, altogether
showing that $\alpha_0(t) \beta_0=\beta_0 \alpha_0(-t)$. Thus,
$\sigma$ is s-malleable.

Assume now that  $N_0$ is a tensor product of $2$ by $2$ matrix
algebras. Consider first the case $N_0=M_{2\times 2}(\Bbb C)$. We
let $\{e_{ij}\}_{i,j=1,2}$ be a matrix unit for $N_0$. Define the
unitaries $\alpha_0(t) \in N_0\otimes N_0$ by
$$
\pmatrix
1 & 0 & 0 & 0\\
0 &\cos \pi t/2 &\sin \pi t/2 &0\\
0 &-\sin \pi t/2 &\cos \pi t/2 &0\\
0 & 0 & 0 &1\\
\endpmatrix
$$
This clearly defines a continuous unitary representation of $\Bbb
R$. Then we define $\beta_0\in \Bbb C \otimes N_0$ to be the
unitary element:
$$
\pmatrix
1 & 0 & 0 & 0\\
0 & 1 & 0 & 0\\
0 & 0 & -1 & 0\\
0 & 0 & 0 & -1\\
\endpmatrix
$$

An easy calculation shows that $\beta_0 \alpha_0(t) \beta_0^{-1} =
\alpha_0(-t)$ and that Ad$(\alpha_0(1))(N_0 \otimes \Bbb C) \perp
N_0\otimes 1$. By taking tensor products of the above construction,
this shows that if $(N_0, \tau_0) =\overline{\otimes}_l (M_{2\times
2}(\Bbb C), tr)_l$ then there exists a continuous action
$\alpha_0:\Bbb R \rightarrow \text{\rm Aut}(N_0\overline{\otimes}
N_0,\tau_0\otimes \tau_0)$ and a period 2-automorphism $\beta_0$ of
$N_0$ such that $\alpha_0(1)(N_0 \otimes 1)\perp N_0\otimes 1$,
sp$(N_0\otimes 1)\alpha_0(1)(N_0\otimes 1)$ dense in
$N_0\overline{\otimes} N_0$, $N_0\otimes 1$ fixed by $\beta$, $\beta
\alpha_0(t)=\alpha_0(-t)\beta, \forall t\in \Bbb R$.

But then $\alpha=\otimes_k (\alpha_0)_k$, $\beta=\otimes_k
(\beta_0)_k$ clearly satisfy $(4.3.1')$, $(4.3.2')$. \hfill
$\square$

\proclaim{4.6. Lemma} Let $\Gamma$ be a countable discrete group and
$H \subset \Gamma$ a subgroup with the relative property $(\text{\rm
T})$. Let $\sigma$ be a s-malleable action of $\Gamma$ on the finite
von Neumann algebra $(P,\tau)$, with $\sigma_{|H}$ weak mixing. Let
$\tilde{\sigma}$, $\alpha: \Bbb R \rightarrow \text{\rm
Aut}(P\overline{\otimes} P, \tau \otimes \tau)$, $\beta \in
\text{\rm Aut}(P\overline{\otimes} P, \tau \otimes \tau)$ be as in
Definition $4.3$. Let $(N,\tau)$ be a finite von Neumann algebra and
$\rho$ and action of $\Gamma$ on it. If $w$ is a cocycle for the
action $\sigma_g \otimes \rho_g$ of $\Gamma$ on $P\overline{\otimes}
N$ then $w_{|H}$ and $\alpha_1(w)_{|H}$ are equivalent as cocycles
for the action $\tilde{\sigma}_h \otimes \rho_h$ of $H$ on
$P\overline{\otimes} P \overline{\otimes} N$.
\endproclaim
\noindent {\it Proof}. Denote $\tilde{P}=P\overline{\otimes} P$. We
still denote by $\alpha, \beta$ the $N$-amplifications
$\alpha\otimes id_N$ and $\beta \otimes id_N$. Also, denote
$\sigma'_g=\sigma_g \otimes \rho_g$ and
$\tilde{\sigma}_g'=\tilde{\sigma}_g \otimes \rho_g$. Since $\alpha$
is a continuous action by automorphisms commuting with
$\tilde{\sigma}$, $\alpha$ acts continuously on the set of cocycles
$w$ for the action $\tilde{\sigma}'$ of $\Gamma$ on
$\tilde{P}\overline{\otimes} N$.

Note that by part $2^\circ$
of Lemma 2.12 it is sufficient to prove that for any $\varepsilon >
0$ there exists $v_1\in \tilde{P}\overline{\otimes} N$ partial
isometry such that $\|v_1-1\|_2 \leq \varepsilon$ and $w_h
\tilde{\sigma}'_h(v_1) = v_1 \alpha_{1}(w_h)$, $\forall h\in H$.
Indeed, because if $u$ is a unitary in $\tilde{P}\overline{\otimes} N$
such that $uv_1^*v_1=v_1$ then
$$
\|w_h
\tilde{\sigma}'_h(u) - u \alpha_{1}(w_h)\|_2 \leq 2 \|u-v_1\|_2 = 2\|1-v_1^*v_1\|_2
\leq 4\varepsilon,$$
forall $h\in H$, and $2.12.2^\circ$ applies.

To construct such $v_1$, note first that by Lemma 4.2 for the given
$\varepsilon > 0$ there exists $n$ such that if we denote
$t_0=2^{-n}$ then there exists a partial isometry $v_0 \in \tilde{P}
\overline{\otimes} N$ satisfying $w_h \tilde{\sigma}'_h(v_0) = v_0
\alpha_{t_0}(w_h)$, $\forall h\in H$, and $\|v_0-1\|_2 \leq
\varepsilon$.

Let us  now show that if
$$
w_h \tilde{\sigma}'_h(v) = v \alpha_{t}(w_h), \forall h\in H, \tag
4.6.1
$$
for some $0 < t <1$ and a partial isometry $v\in
\tilde{P}\overline{\otimes} N$, then there exists a partial isometry
$v'\in \tilde{P} \overline{\otimes} N$ satisfying $\|v'\|_2=\|v\|_2$
and $w_h \tilde{\sigma}'_h(v') = v' \alpha_{2t}(w_h)$, $\forall h\in
H$. This will of course prove the existence of $v_1$, by starting
with $t=t_0=2^{-n}$ then proceeding by induction until we reach
$t=1$  (after $n$ steps).

Applying $\beta$ to $(4.6.1)$ and using that $\beta$ commutes with
$\tilde{\sigma}'$, $\beta(x)=x, \forall x\in P \overline{\otimes}
N\subset \tilde{P}\overline{\otimes} N$ and $\beta \alpha_t =
\alpha_{-t} \beta$ as automorphisms on $\tilde{P}\overline{\otimes}
N$, we get $\beta(w_h)=w_h$ and
$$
w_h \tilde{\sigma}'_h(\beta(v))=\beta(v) \alpha_{-t}(w_h), \forall
h\in H. \tag 4.6.2
$$
Since $(4.6.1)$ can be read as $v^* w_h =
\alpha_t(w_h)\tilde{\sigma}'_h(v^*)$, from $(4.6.1)$ and $(4.6.2)$
we get the identity

$$
v^* \beta(v) \alpha_{-t}(w_h)=v^* w_h \tilde{\sigma}'_h(\beta(v))
$$
$$
=\alpha_t(w_h)\tilde{\sigma}'_h(v^*)\tilde{\sigma}'_h(\beta(v))
=\alpha_t(w_h)\tilde{\sigma}'_h(v^*\beta(v)),
$$
for all $h\in H$. By applying $\alpha_t$ on both sides of this
equality, if we denote $v'=\alpha_t(\beta(v)^*v)$ then we further
get
$$
{v'}^* w_h = \alpha_{2t} (w_h) \tilde{\sigma}'_h({v'}^*), \forall
h\in H,
$$
showing that $w_h \tilde{\sigma}'_h(v') = v' \alpha_{2t}(w_h)$,
$\forall h\in H$, as desired. On the other hand, the intertwining
relation $(4.6.1)$ implies that $vv^*$ is in the fixed point algebra
$B$ of the action Ad$w_h \circ \tilde{\sigma}'_h$ of $H$ on
$\tilde{P}\overline{\otimes} N$. Since $\tilde{\sigma}'_{|H}$ is
weak mixing on $(1 \otimes P) \otimes 1 \subset
\tilde{P}\overline{\otimes} N$ (because it coincides with $\sigma$
on $1_P \otimes P \otimes 1_N \simeq P$) and because Ad$w_h$ acts as
the identity on $(1 \otimes P) \otimes 1$ and leaves $(P\otimes 1)
\overline{\otimes} N$ globally invariant, it follows that $B$ is
contained in $(P \otimes 1) \overline{\otimes} N$. Thus $\beta$ acts
as the identity on it (because it acts as the identity on both
$P\otimes 1$ and $1\otimes N$). In particular $\beta(vv^*)=vv^*$,
showing that the right support of $\beta(v^*)$ equals the left
support of $v$. Thus, $\beta(v^*)v$ is a partial isometry having the
same right support as $v$, implying that $v'$ is a partial isometry
with $\|v'\|_2 = \|v\|_2$. \hfill $\square$

\heading 5. Cocycle and OE superrigidity results
\endheading

\vskip .1in \noindent {\it 5.1. Definitions}. An infinite subgroup
$H$ of a group $\Gamma$ is {\it w-normal} (resp. {\it wq-normal}) in
$\Gamma$ if there exists an ordinal $\imath$ and a well ordered
family of intermediate subgroups $H = H_0 \subset H_1 \subset ...
\subset H_{\jmath} \subset ... \subset H_{\imath} = \Gamma$ such
that for each $0 < \jmath \leq \imath$, the group $H'_\jmath=\cup_{n
< \jmath} H_n$ is normal in $H_\jmath$ (resp. $H_{\jmath}$ is
generated by the elements $g \in \Gamma$ with $|gH'_{\jmath}g^{-1}
\cap H'_{\jmath}|=\infty$).

For instance, if $H=H_0 \subset H_1 \subset ... \subset H_n =\Gamma$
are all normal inclusions and $H$ is infinite then $H \subset
\Gamma$ is w-normal. An inclusion of the form $H \subset (H *
\Gamma_0) \times \Gamma_1$ with $H, \Gamma_1$ infinite and
$\Gamma_0$ arbitrary is wq-normal but not w-normal.

All statements in this Section will be about inclusions of countable
infinite groups $H \subset \Gamma$ with the relative property (T) of
Kazhdan-Margulis (or $H\subset \Gamma$ is a rigid subgroup,
i.e. it checks condition $(4.1.2)$)
such that $H$ is either wq-normal or w-normal in $\Gamma$.
The class of groups having
infinite wq-normal rigid subgroups is closed to finite index
restriction/extension, to normal extensions and to operatios
such that $\Gamma \mapsto (\Gamma * \Gamma_0) \times \Gamma_1$,
where $\Gamma_0$ is arbitrary and $|\Gamma_1 |=\infty.$ The class of groups
having infinite w-normal rigid subgroups is closed to normal extensions.
Both classes are closed to inductive limits.

Any lattice $\Gamma$ in a Lie group of real rank at least 2, such as
$SL(n,\Bbb Z)$, $n \geq 3$, has property (T) by Kazhdan's celebrated
result ([K]), so $H=\Gamma \subset \Gamma$ is rigid. The normal
inclusions of groups $\Bbb Z^2\subset \Gamma=\Gamma_0 \ltimes \Bbb
Z^2$, with $\Gamma_0 \subset SL(2,\Bbb Z)$ non amenable, are rigid
(cf. [K], [M], [B]), and so are all inclusions $\Bbb Z^n \subset
\Gamma=\Gamma_0 \ltimes \Bbb Z^n$ in ([Va], [Fe]). For each such $H
\subset \Gamma$, the inclusion $H \subset \Gamma'= \Gamma \times H'$
is w-normal rigid, for any group $H'$. Also, if $H \subset \Gamma$
is wq-normal rigid then $H \subset (\Gamma*\Gamma_0) \times
\Gamma_1$ is wq-normal rigid whenever $\Gamma_1$ is an infinite
group, $\Gamma_0$ arbitrary, but it can be shown that if $\Gamma_0$
is non-trivial then $(\Gamma * \Gamma_0) \times \Gamma_1$ has no
w-normal rigid subgroups.

Notice that in the cocycle superrigidity results 5.2-5.5 below the
actions $\sigma$ are not assumed to be free.

\proclaim{5.2. Theorem (Cocycle superrigidity of s-malleable
actions)} Let $\sigma$ be a s-malleable action of a countable
discrete group $\Gamma$ on the standard probability space $(X,\mu)$.
Let $H \subset \Gamma$ be a rigid subgroup such that $\sigma_{|H}$
is weak mixing. Let $\Cal V\in \mycal U_{fin}$.

If $\rho$ is an arbitrary action of $\Gamma$ on a standard
probability space $(Y,\nu)$, then any $\Cal V$-valued cocycle $w$
for the diagonal product action $\sigma_g\times \rho_g, g\in
\Gamma,$ on $(X\times Y, \mu \times \nu)$ is cohomologous to a $\Cal
V$-valued cocycle $w'$ whose restriction to $H$ is independent on
the $X$-variable. If in addition $H$ is w-normal in $\Gamma$, or if
$H$ is wq-normal in $\Gamma$ but $\sigma$ is mixing, then $w'$ is
independent on the $X$-variable on all $\Gamma$, in other words the
inclusion $\text{\rm Z}^1(\rho;\Cal V) \subset \text{\rm Z}^1(\sigma
\times \rho; \Cal V)$ implements an isomorphism between the sets of
equivalence classes $\text{\rm Z}^1(\rho;\Cal V)/\sim$ and
$\text{\rm Z}^1(\sigma \times \rho; \Cal V)/\sim.$

In particular, any cocycle $w$ for $\sigma$ with values in $\Cal V$
is cohomologous to a cocycle $w'$ which restricted to $H$ is a group
morphism of $H$ into $\Cal V$ and if $H$ is w-normal in $\Gamma$, or
if $H$ is wq-normal in $\Gamma$ but $\sigma$ is mixing, then
$\text{\rm Z}^1(\sigma;\Cal V)=\text{\rm Z}^1_0(\sigma;\Cal V)$.
\endproclaim
\noindent {\it Proof}. The last part follows by simply taking $\rho$
to be the $\Gamma$-action on the one point probability space. To
prove the first part, by Proposition 3.5 it is sufficient to
consider the case $\Cal V = \Cal U(N_0)$, with $N_0$ a separable
finite von Neumann algebra. But if $\Cal V=\Cal U(N_0)$ and $\sigma$
s-malleable, then the statement follows from Lemma 4.6 and Theorem
3.1. \hfill $\square$

\proclaim{5.3. Theorem (Cocycle superrigidity of sub s-malleable
actions)} The same statement as in $5.2$ holds true if instead of
$\sigma$ s-malleable with $\sigma_{|H}$ weak mixing, we merely
assume that $\sigma$ is the quotient of a s-malleable action $\Gamma
\curvearrowright^{\sigma'} (X',\mu')$ such that $\sigma'_{|H}$ is
weak mixing on $(X',\mu')$ and weak mixing relative to $L^\infty X$.
\endproclaim
\noindent {\it Proof}. Trivial by 5.2 and 2.11. \hfill $\square$

\proclaim{5.4. Corollary (Cocycle superrigidity of generalized
Bernoulli actions)} Let $\sigma$ be a $(X_0, \mu_0)$-Bernoulli
$\Gamma \curvearrowright K$ action, where $\Gamma$ is a countable
discrete group acting on the countable set $K$ and $(X_0,\mu_0)$ is
a non-trivial standard probability space. Let $H \subset \Gamma$ be
an infinite subgroup with the relative property $(\text{\rm T})$
such that all orbits of $H \curvearrowright K$ are infinite (see
$4.5.2^\circ$). Let $\Cal V$ be a Polish group of finite type.

Then any cocycle $w: X \times \Gamma \rightarrow \Cal V$ for
$\sigma$ with values in $\Cal V$ is cohomologous to a cocycle $w'$
whose restriction to $H$ is a group morphism of $H$ into $\Cal V$.
If in addition $H$ is w-normal in $\Gamma$, or if $H$ is wq-normal
but $\{g\in \Gamma \mid gk=k\}$ is finite, $\forall k\in K$ (see
$4.5.3^\circ$), then $w'$ follows a group morphism on all $\Gamma$.

More generally, if $\rho$ is an arbitrary action of $\Gamma$ on a
standard probability space $(Y,\nu)$ and $w$ is a cocycle for the
diagonal product action $\sigma_g \times \rho_g$, $g\in \Gamma,$ on
$(X\times Y,\mu\times \nu)$, then any $\Cal V$-valued cocycle $w$
for $\sigma \times \rho$ is cohomologous to a $\Cal V$-valued
cocycle $w'$ whose restriction to $H$ is independent on the
$X$-variable, i.e. to a cocycle for $\rho_{|H}$. Moreover, if $H$ is
w-normal in $\Gamma$, or if $H$ is wq-normal but $\{g\in \Gamma \mid
gk=k\}$ finite, $\forall k\in K$, then $w'$ follows independent on
the $X$-variable on all $\Gamma$.
\endproclaim

\proclaim{5.5. Theorem (Cocycle superrigidity: the non-commutative
case)} Let $\sigma$ be a s-malleable action of a countable discrete
group $\Gamma$ on a finite von Neumann algebra $(P,\tau)$. Let
$H\subset \Gamma$ be a rigid subgroup such that $\sigma_{|H}$ is
weak mixing. Let $(N,\tau)$ be an arbitrary finite von Neumann
algebra and $\rho$ an action of $\Gamma$ on $(N,\tau)$. Then any
cocycle $w$ for the diagonal product action $\sigma\otimes \rho$ of
$\Gamma$ on $P \overline{\otimes} N$ is equivalent to a cocycle $w'$
whose restriction to $H$ takes values in $N=1\otimes N$. If in
addition $H$ is w-normal in $\Gamma$, or if $H$ is wq-normal in
$\Gamma$ but $\sigma$ is mixing, then $w'$ takes values in $N$
on all $\Gamma$.
\endproclaim
\noindent {\it Proof}. Let $\sigma'_g=\sigma_g \otimes \rho_g$,
$g\in \Gamma$. Let $\tilde{\sigma}: \Gamma \rightarrow \text{\rm
Aut}(\tilde{P}, \tilde{\tau})$ denote the diagonal product action
$\tilde{\sigma}_g=\sigma_g \otimes \sigma_g$ of $\Gamma$ on
$\tilde{P}=P\overline{\otimes} P$ and $\alpha: \Bbb R \rightarrow
\text{\rm Aut}(\tilde{P}, \tilde{\tau})$, $\beta \in \text{\rm
Aut}(\tilde{P}, \tilde{\tau})$ as in Definition $4.3$. By Lemma 4.6,
$w_{|H} \sim \alpha_1(w)_{|H}$ as cocycles for $\tilde{\sigma}_h
\otimes \rho_h$, $h\in H$. But then Proposition 3.2 implies $w_{|H}$
is equivalent to a cocycle $w'$ with $w'_h \in \Cal U(N)$, $\forall
h\in H$. \hfill $\square$

\vskip .05in

We'll now apply the case $\Cal V=\Lambda$ discrete of the Cocycle
Superrigidity 5.3 to deduce orbit equivalence superrigidity of the
corresponding source actions. The results will in fact hold true for
all source actions that are cocycle superrigid, in the following
sense:

\vskip .05in \noindent {\bf 5.6.0. Terminology}. Let $\mycal U$ be a
family of Polish groups. We say that $\Gamma \curvearrowright^\sigma
(X,\mu)$ is $\mycal U$-{\it cocycle superrigid} if any $\Cal
V$-valued cocycle for $\sigma$ is cohomologous to a group morphism
of $\Gamma$ into $\Cal V$, $\forall \Cal V\in \mycal U$. If $\mycal
U$ is the family of all discrete groups, we simply say that $\sigma$
is {\it cocycle superrigid}. With this terminology, Theorem 5.2
shows that if $\Gamma$ has an infinite subgroup $H\subset \Gamma$
with the relative property (T), $\sigma$ is a weak mixing
s-malleable $\Gamma$-action and either $H$ is w-normal in $\Gamma$,
with $\sigma_{|H}$ weak mixing, or $H$ is wq-normal, with $\sigma$
mixing, then $\sigma$ is cocycle superrigid. By 5.3, for an action
$\Gamma \curvearrowright^\sigma X$ to be cocycle superrigid it is in
fact sufficient that it is the quotient of a cocycle superrigid
action $\Gamma \curvearrowright^{\sigma'} X'$ with the property that
$\Gamma \curvearrowright L^\infty X'$ is weak mixing relative to
$L^\infty X$ (in the sense of Definition 2.9).

\vskip .05in \noindent {\bf 5.6.1. Assumption}. $\Gamma
\curvearrowright^\sigma  X$ is free and, as an action on $L^\infty
X$, it admits an extension to a s-malleable action $\sigma'$ of
$\Gamma$ on a larger abelian von Neumann algebra $(A',\tau')$ with
the property that $\tau'_{|A}=\tau$, $\sigma'$ weak mixing relative
to $L^\infty X\subset A'$ and such that there exists an infinite
rigid subgroup $H\subset \Gamma$ with $\sigma'_{|H}$ weak mixing on
$A'$ and either $H$ is w-normal in $\Gamma$ or $H$ is merely
wq-normal but $\sigma'$ mixing. \vskip .05in

Examples of $(\sigma,\Gamma)$ satisfying 5.6.1 are all
$(X_0,\mu_0)$-Bernoulli $\Gamma \curvearrowright K$ actions of
groups $\Gamma$ having an infinite rigid subgroup $H \subset \Gamma$
with the property that $\Gamma \curvearrowright K$ satisfies
$4.5.1^\circ$ and either $H \curvearrowright K$ checks $4.5.2^\circ$
and $H\subset \Gamma$ is w-normal, or $\Gamma$ checks $4.5.3^\circ$ and
$H$ is wq-normal.

Recall from 1.4 that if $\Cal R$,  $\Cal S$ are countable m.p.
relations on the probability spaces $(X,\mu), (Y,\nu)$ then an {\it
orbit equivalence} ({\it OE}) of $\Cal R$, $\Cal S$ is an
isomorphism of probability spaces $\Delta:(X,\mu) \simeq (Y,\nu)$
which takes $\Cal R$ onto $\Cal S$. A {\it local OE} of $\Cal R$,
$\Cal S$ is a m.p. map $\Delta:(X,\mu) \rightarrow (Y,\nu)$ (a.e.
surjective, but not necessarily a.e. $1$ to $1$) with the property
that for almost all $t\in X$, $\Delta$ gives a bijection between the
$\Cal R$-orbit of $t$ and the $\Cal S$-orbit of $\Delta(t)$ (see
1.4.2). For instance, if $\Gamma \curvearrowright^\sigma (X,\mu)$ is
a free m.p action, $B \subset L^\infty X$ is a
$\sigma(\Gamma)$-invariant von Neumann subalgebra of $L^\infty X$
such that $\sigma_{|B}$ is still free and we write $B=L^\infty Y$,
for some probability space $(Y,\nu)$ with the property that
$\tau_\mu$ equals $\tau_\nu$ on $B$, then the associated quotient
map $\Delta:(X,\mu) \rightarrow (Y,\nu)$ (see 1.1.1) implements a
local OE of $\Cal R_\sigma$, $\Cal R_\theta$, where $\theta_g$ is
the action of $\Gamma$ on $(Y,\mu)$ implemented by the restriction
of $\sigma_g$ to $L^\infty Y=B $. An {\it embedding} of $\Cal R$
into $\Cal S$ is an isomorphism $\Delta:(X,\mu) \simeq (Y,\nu)$ that
takes $\Cal R$ onto a subequivalence relation of $\Cal S$.

\proclaim{5.6. Theorem (OE Superrigidity)} Let $\Gamma
\curvearrowright^\sigma X$ be a free, weakly mixing, cocycle
superrigid action, e.g. an action satisfying $5.6.1$. Assume
$\Gamma$ has no finite normal subgroups. Let $\theta$ be an
arbitrary free ergodic m.p. action of a countable discrete group
$\Lambda$ on a standard probability space $(Y,\nu)$ and $\Delta:\Cal
R_\sigma \simeq \Cal R_\theta^t$ an orbit equivalence, for some $t
> 0$.

Then $n=1/t$ is an integer and there exist a subgroup $\Lambda_0
\subset \Lambda$ of index $[\Lambda: \Lambda_0]=n$, a subset $Y_0
\subset Y$ of measure $\nu(Y_0)=1/n$ fixed by $\theta_{|\Lambda_0}$,
an automorphism $\alpha \in [\Cal R_{\theta}]$ and a group
isomorphism $\delta: \Gamma \simeq \Lambda_0$ such that $\alpha
\circ \Delta$ takes $X$ onto $Y_0$ and conjugates the actions
$\sigma, \theta_0 \circ \delta$, where $\theta_0$ denotes the action
of $\Lambda_0$ on $Y_0$ implemented by $\theta$.

Moreover, if $\Gamma$ is assumed ICC but the free action $\Gamma
\curvearrowright^\sigma X$ is merely a quotient of a weakly mixing
cocycle superrigid action (e.g. of an action satisfying $5.6.1$),
then $\Gamma \curvearrowright X$ is OE superrigid, i.e. if $\Lambda
\curvearrowright Y$ is a free ergodic m.p. action of an arbitrary
countable group $\Lambda$ and $\Delta:X \simeq Y$ is an orbit
equivalence of $\Gamma \curvearrowright X, \Lambda \curvearrowright
Y$ then there exist $\alpha \in [\Lambda]$ and $\delta:\Gamma \simeq
\Lambda$ such that $\Delta_0=\alpha \circ \Delta$ conjugates
$\sigma, \theta \circ \delta$.
\endproclaim

\proclaim{5.7. Theorem (Superrigidity of local OE)} Let $\Gamma
\curvearrowright^\sigma (X,\mu)$ be a free weakly mixing cocycle
superrigid action, e.g. an action satisfying $5.6.1$. Assume
$\Gamma$ has no finite normal subgroups. Let $\theta$ be an
arbitrary free ergodic m.p. action of a countable discrete group
$\Lambda$ on a standard probability space $(Y,\nu)$. Let $\Delta:
(X,\mu) \rightarrow (Y,\nu)^t$ be a local OE of $\Cal R_\sigma$,
$\Cal R_\theta^t$, for some $t > 0$.

Then $n=t^{-1}$ is an integer and if $\theta'$ denotes the pull back
of $\theta$ to a $\Lambda$-action on $(X^n, \mu_{X^n})$ (see $1.4$)
then there exist a subgroup $\Lambda_0 \subset \Lambda$ of index
$[\Lambda: \Lambda_0]=n$, a subset $X_0 \subset X^n$ of measure
$\mu_{X^n}(X_0)=1/n$ fixed by $\theta'_{|\Lambda_0}$, an
automorphism $\alpha \in [\Cal R_{\sigma^n}]$ and a group
isomorphism $\delta: \Gamma \simeq \Lambda_0$ such that $\alpha$
takes $X$ onto $X_0$ and conjugates the actions $\sigma, \theta_0
\circ \delta$, where $\theta_0$ denotes the action of $\Lambda_0$ on
$X_0$ implemented by $\theta'$.

In particular, if $\Delta$ is a local OE of $\Cal R_\sigma, \Cal
R_\theta$ then there exist $\alpha \in [\Cal R_\sigma]$ and
$\delta:\Gamma \simeq \Lambda$ such that $\theta \circ \delta$ is a
quotient of $\alpha \sigma \alpha^{-1}$ (via $\Delta$).
\endproclaim

\proclaim{5.8. Theorem (Superrigidity of embeddings)} Let $\Gamma
\curvearrowright^\sigma (X,\mu)$ be a free weakly mixing cocycle
superrigid action, e.g. an action satisfying $5.6.1$. Assume
$\Gamma$ has no finite normal subgroups. Let $\Lambda
\curvearrowright^\theta (Y,\nu)$ be an arbitrary action and $\Delta:
(X,\mu) \simeq (Y,\nu)^t$ an embedding of $\Cal R_\sigma$ into $\Cal
R^t_\theta$, for some $t>0$, such that any $\Gamma$-invariant finite
subequivalence relation of $\Cal R_\theta^t$ is contained in $\Cal
R_\sigma$ (when identifying $\Cal R_\sigma$ with a subequivalence of
$\Cal R_\theta^t$, via $\Delta$).

Then $t \leq 1$ and there exist an isomorphism $\delta$ of $\Gamma$
onto a subgroup $\Lambda_0$ of $\Lambda$ and $\alpha \in [\Cal
R_{\theta}]$ such that $\Delta_0=\alpha\circ \Delta$ takes $X$ onto
a $\Lambda_0$-invariant subset $Y_0\subset Y$, with $\mu(Y_0)=t$,
and conjugates the action $\sigma$ with the action
$\theta_{|\Lambda_0}$ of $\Lambda_0$ on $Y_0$, with respect to the
isomorphism $\delta:\Gamma \simeq \Lambda_0$.
\endproclaim

\proclaim{5.9. Corollary} Let $\Lambda \curvearrowright (X,\mu)$ be
a free ergodic m.p. action and assume there exists $t>0$ such that
$\Cal R^t_\Lambda \supset \Cal R_\Gamma$ for some  $\Gamma
\curvearrowright (X,\mu)$ satisfying $5.6.1$, $\Gamma$ without
finite normal subgroups and such that any $\Gamma$-invariant finite
subequivalence relation of $\Cal R_\Lambda^t$ is contained in $\Cal
R_\Gamma$. Then $\mycal F(\Cal R_\Lambda)=\{1\}$. Thus, if $\Cal R$
is a countable m.p. equivalence relation with $\mycal F(\Cal R) \neq
\{1\}$ and such that $\Cal R\supset \Cal R_\Gamma$ for some $\Gamma
\curvearrowright X$ satisfying the above properties then $\Cal R^t$
cannot be implemented by a free group action, $\forall t
> 0$.
\endproclaim

\proclaim{5.10. Corollary} Let $\Gamma$ be a group having an
infinite wq-normal rigid subgroup and no finite normal subgroup. Let
$(X,\mu)=(\{0,1\}, \mu_0)^\Gamma$, with
$s=\mu_0(\{0\})/\mu_0(\{1\})\neq 1$. Let $\Cal R_0$ be the
hyperfinite equivalence relation on $\{0,1\}^\Gamma$ given by:
$(t_g)_g \sim (t'_g)_g$ iff there exists a finite subset $F\subset
\Gamma$ such that $t_g=t'_g$, $\forall g\in \Gamma \setminus F$ and
$\Pi_{g \in F} \mu_0(t_g)= \Pi_{g \in F} \mu_0(t'_g)$. Let $\Cal R$
be the countable m.p. equivalence relation generated by $\Cal R_0$
and the Bernoulli $\Gamma$-action $\Gamma \curvearrowright X$ (which
leaves $\Cal R_0$ invariant). Then $\Cal R^t$ cannot be implemented
by a free group action, $\forall t>0$. Moreover, $\mycal F(\Cal
R)\supset s^{\Bbb Z}$ and if $\Gamma$ is of the form $\Gamma_0
\ltimes \Bbb Z^2$, with $\Gamma_0$ a finitely generated non
virtually cyclic subgroup of $SL(2,\Bbb Z)$, then $\mycal F(\Cal
R)=s^{\Bbb Z}$.
\endproclaim

We deduce Theorems 5.6-5.8 and their Corollaries 5.9, 5.10 from the
Cocycle Superrigidity 5.2-5.3 and a {\it General Principle}
(Proposition 5.11 below) showing that any ``untwister'' of a Zimmer
cocycle associated with an orbit equivalence of actions (defined
below) gives rise to a natural ``conjugator'' of the two actions. We
in fact prove a more general such principle, dealing also with
embeddings and local OE of equivalence relations, and their
amplifications. Prior results along this line have been obtained in
(4.2.9, 4.2.11 of [Z2]) and (3.3 of [Fu2], 2.4 of [Fu3]). The proof
uses von Neumann algebra framework, but is otherwise quite
straightforward.

To formulate this result, let us recall the definition of the {\it
Zimmer cocycle} $w_\Delta$ associated with a morphism $\Delta$ of
equivalence relations implemented by group actions, with the target
action free. Thus, let $\sigma$ be a free $\Gamma$-action on
$(X,\mu)$ and $\theta$ a free $\Lambda$-action on $(Y,\nu)$. Let
$Y_0\subset Y$ be a subset of positive measure and assume $\Delta:
(X,\mu) \rightarrow (Y_0,\nu_{Y_0})$ is a m.p. morphism of $\Cal
R_\sigma$ into $\Cal R_\theta^{Y_0}$, with $N_0\subset X$ a subset
of measure 0 such that $\Delta$ takes the $\Cal R_\sigma$-orbit of
any $t\in X\setminus N_0$ into the $\Cal R_\theta^{Y_0}$-orbit of
$\Delta(t)$. Let $w=w_\Delta: X \times \Gamma \rightarrow \Lambda$
be defined as follows: For given $t\in X\setminus N_0$, $g\in
\Gamma$, $w(t,g)$ is the unique $h\in \Lambda$ such that
$\Delta(\sigma_g^{-1}(t))=\theta_h^{-1}(\Delta (t))$. It is
immediate to check that $w: X \times \Gamma \rightarrow \Lambda$
this way defined is measurable and satisfies the (right) cocycle
relation $(2.1.1)$ (see [Z2]).

The ``principle'' below shows that if $\Gamma
\curvearrowright^\sigma (X,\mu)$ is free and weak mixing, $ \Lambda
\curvearrowright^\theta (Y,\nu)$ is free and $\Delta$ is a (local)
OE (resp. an embedding) of the corresponding equivalence relations
$\Cal R_\sigma, \Cal R_\theta$, then a measurable map $v: X
\rightarrow \Lambda$ implementing an equivalence between $w_\Delta$
and a group morphism $\delta: \Gamma \rightarrow \Lambda$ gives rise
to a natural automorphism $\alpha=\alpha_v$ in the full group
$[\Cal R_\sigma]$ (resp. $[\Cal R_\theta]$) that (virtually)
conjugates the $\Gamma$-actions $\Delta \sigma \Delta^{-1}, \theta
\circ \delta$.

\proclaim{5.11. Proposition} Let $\Gamma \curvearrowright^\sigma
(X,\mu)$, $\Lambda \curvearrowright^\theta (Y,\nu)$ be free m.p.
actions, with $\sigma$ weak mixing, and $\Delta:(X,\mu) \rightarrow
(Y_0,\nu_{Y_0})$ either an embedding or a local OE of $\Cal
R_\sigma, \Cal R^{Y_0}_\theta$, for some $Y_0\subset Y$. Denote
$A=L^\infty X$, $B=L^\infty Y$, $M=A \rtimes \Gamma$, $P=B \rtimes
\Lambda$, $p=\chi_{Y_0}$, and $u_g\in M$, $v_h\in P$ the canonical
unitaries implementing $\sigma,\theta$.

Let $w=w_\Delta$ denote the Zimmer cocycle associated to $\Delta$.
Assume $\delta: \Gamma \rightarrow \Lambda$ is a group morphism and
$v: X \rightarrow \Lambda$ a measurable map such that for each $g\in
\Gamma$ we have $v(t)^{-1} w(t, g) v(g^{-1}t)=\delta(g)$, $\forall
t\in X$ (a.e.). Let $\{X_h\}_h$ be the partition of $X$ into
measurable subsets such that $v(t)=h$ for $t\in X_h$, $h\in
\Lambda$. Denote $q_h= \chi_{X_h}\in L^\infty X$, $h\in \Lambda$.
Let $b\overset \text{\rm def} \to =\Sigma_h q_h v_h$, which is
viewed as an element in $L^2P$ when $\Delta$ is an embedding, via
the inclusion $M\subset P$ (see 1.4), and is viewed as an element in
$L^2 M^n$ when $\Delta$ is a local OE, via the inclusion $P\subset
M^n$ (see 1.4.3), where $n=\nu(Y_0)^{-1}$.

Then $u_g b = b v_{\delta(g)}, \forall g\in \Gamma$, $b$ is a scalar
multiple of a partial isometry, $K=ker(\delta)$ is a finite normal
subgroup of $\Gamma$ and we have:

$(i)$. If $\Delta$ is a local OE then $bb^*=\Sigma_{k\in K} u_k$,
$b^*b=|K|q$ with $q\in A^n$ a projection invariant to the pull back
$\theta'_{h}, h\in \Lambda_0=\delta(\Gamma)$, and
$n=[\Lambda:\Lambda_0]$ is an integer. Moreover, if we denote
$A^K=L^\infty (X/K)$ the fixed point subalgebra of $A$ under the
action $\sigma_{|K}$, then $\alpha=|K|^{-1}\text{\rm Ad}(b^*)$
implements an isomorphism of $A^K$ onto $A^nq$ which conjugates the
$\Gamma/K$-action $\sigma_0=\sigma_{|A^K}$ and the
$\Lambda_0$-action $\theta'_{|\Lambda_0}$ with respect to the
identification $\Gamma/K \simeq \Lambda_0$ implemented by $\delta$.
In particular, if $\Gamma$ has no finite normal subgroups then
$b^*b=q, bb^*=p$, $\alpha=\text{\rm Ad}b^*$ has graph in $\Cal
R_\sigma^n=\Cal R_{\theta'}$  and  conjugates $\sigma,
\theta'_{|\Lambda_0}$, with respect to the isomorphism $\delta:
\Gamma \simeq \Lambda_0\subset \Lambda$.

$(ii)$. If $\Delta$ is an embedding and we identify $A^n=B$, $M^n
\subset P$ and $\Cal R_\sigma \subset \Cal R^X_\theta$ via $\Delta$,
then there exists a finite $\sigma(\Gamma)$-invariant subequivalence
relation $\Cal K\subset \Cal R_{\theta}^{X}$ containing $\Cal
R_{\sigma(K)}$, such that if we denote by $m$ the cardinality of the
orbits of $\Cal K$, $A^{\Cal K}\simeq L^\infty(X/\Cal K)$ the
quotient (or fixed point algebra) under $\Cal K$ then
$\alpha=m^{-1}\text{\rm Ad}(b^*)$ implements an isomorphism of
$A^{\Cal K}$ onto $A^nq$ which conjugates the $\Gamma/K$-action
$\sigma_0=\sigma_{|A^{\Cal K}}$ onto the $\Lambda_0$-action
$\theta_0$ implemented by $\theta_{|\Lambda_0}$ on $Bq=A^nq$, with
respect to the identification $\Gamma/K \simeq \Lambda_0$
implemented by $\delta$. If in addition any finite
$\sigma(\Gamma)-$invariant subequivalence relation of $\Cal
R_\theta^X$ is contained in $\Cal R_\sigma$ then $\Cal K=\Cal
R_{\sigma(K)}$, $b^*b=|K| q$, $bb^*=\Sigma_{k\in K} u_k$,
$\alpha=|K|^{-1}\text{\rm Ad}(b^*)$ and $\alpha$ implements an
isomorphism of $A^K$ onto $A^nq$. Also, if $\Gamma$ has no
non-trivial finite normal subgroups, then $b$ is a partial isometry
in $P$ normalizing $A^n=B$ with left support $p=\chi_X$ that
conjugates the $\Gamma$-action $\sigma$ and the $\Lambda_0$-action
$\theta_0$.
\endproclaim

Let us first prove 5.6-5.10 assuming Proposition 5.11, then we prove
the latter.

\vskip .05in \noindent {\it Proof of Theorem 5.6}. The first part is
trivial by Theorem 5.7. To prove the last part, assume $\Gamma
\curvearrowright^\sigma X$ is a quotient of a free weakly mixing
cocycle superrigid action $\Gamma \curvearrowright^{\sigma'}
(X',\mu')$. Let $A'=L^\infty (X',\mu')$, $A=L^\infty X \subset A'$,
$M=A' \rtimes \Gamma$ and $N=A\rtimes \Gamma\subset M$, with
$\{u_g\}_g \in M$ the canonical unitaries implementing $\sigma'$ on
$A'$ and $\sigma$ on $A\subset A'$. Let $\{v_h\}_{h\in
\Lambda}\subset N$ be the canonical unitaries implementing $\theta$
on $L^\infty Y=L^\infty X$ (identification via $\Delta$). By
Proposition 5.11, there exists a unitary element $u'=\Sigma_h q_h
v_h\in M$ normalizing $A'$ such that $u' u_g = v_{\delta(g)} u',
\forall g\in \Gamma$. Letting $u$ be the expectation
$E_N(u')=\Sigma_h E_A(q_h)v_h\neq 0$ of $u'$ on $N$, since $u_g, v_h
\in N$, it follows that we still have $u u_g = v_{\delta(g)} u,
\forall g\in \Gamma$. Thus ${u'}^*u\in \{u_g\}_g' \cap M$. But
$\sigma'$ weak mixing implies $\{u_g\}_g' \cap M \subset L\Gamma$
and if in addition $\Gamma$ is ICC then $\{u_g\}_g' \cap M \subset
L\Gamma=L\Gamma ' \cap L\Gamma = \Bbb C$. Thus $\Sigma_h E_A(q_h)v_h
\in \Bbb C u'$, implying that $E_A(q_h)=q_h, \forall h$, i.e. $u'
\in N$. But then $\alpha=\text{\rm Ad}u'$ will do. \hfill $\square$

\vskip .05in \noindent {\it Proof of Theorem 5.7}. Let $m$ be the
smallest integer $\geq 1$ such that $m \geq t$. Let
$\tilde{\Lambda}=\Lambda \times (\Bbb Z/m\Bbb Z)$ and denote by
$\tilde{\theta}$ the action of $\tilde{\Lambda}$ on $\tilde{Y}=Y
\times (\Bbb Z/m\Bbb Z)$ given by the product of the actions
$\theta$ and the translations by elements in $\Bbb Z/m\Bbb Z$. Also,
denote $s=m/t$ and $\tilde{\theta}'$ the $\Delta^s$-pull back of
$\tilde{\theta}$, where $\Delta^s:(X,\mu)^s \simeq (\tilde{Y},
\tilde{\nu})$ is the $s$-amplification of $\Delta$. Thus, we may
regard $\Cal R_\theta^t$ as $\Cal R_{\tilde{\theta}}^{Y_0}$, for
some subset $Y_0\subset \tilde{Y}$ of measure $t/m$, and $\Delta^s$
as a local OE of $\Cal R_\sigma^s=\Cal R_{\tilde{\theta}'}$, $\Cal
R_{\tilde{\theta}}$. By Proposition $5.11(i)$, it follows that $m/t$
is an integer. By the choice of $m$ this implies either $t\leq 1$
(and so $m=1$) or $t=m$ is an integer $\geq 2$. The second case
implies $\Cal R_\sigma=\Cal R_{\tilde{\theta}'}$, and since $\sigma$
is cocycle superrigid by 5.11$(i)$ the actions $\sigma,
\tilde{\theta}'$ are conjugate. But $\sigma$ is weak mixing and
$\tilde{\theta}'$ is not, contradiction. Thus $t \leq 1$ and the
statement follows now trivially from the cocycle superrigidity of
$\sigma$ and 5.11$(i)$. \hfill $\square$

\vskip .05in \noindent {\it Proof of Theorem 5.8}. Let $m$ be the
smallest integer $\geq 1$ such that $m \geq t$. Let
$\tilde{\Lambda}=\Lambda \times (\Bbb Z/m\Bbb Z)$ and denote by
$\tilde{\theta}$ the action of $\tilde{\Lambda}$ on $\tilde{Y}=Y
\times (\Bbb Z/m\Bbb Z)$ given by the product of the actions
$\theta$ and the translations by elements in $\Bbb Z/m\Bbb Z$. Thus,
we may regard $\Cal R_\theta^t$ as $\Cal R_{\tilde{\theta}}^{Y_0}$,
for some subset $Y_0\subset \tilde{Y}$ of measure $t/m$. Let
$s=m/t$. Thus, we may apply $5.11(ii)$ to get a subset $Y_0'\subset
\tilde{Y}$, with $\mu(Y_0')=\mu(Y_0)=t/m$, a subgroup $\Lambda_0
\subset \tilde{\Lambda}$, such that $\tilde{\theta}_{|\Lambda_0}$
leaves $Y_0'$ invariant, and $\alpha \in [\Cal R_{\tilde{\theta}}]$
that takes $Y=\Delta(X)$ onto $Y_0'$, such that $\alpha \circ
\Delta$ conjugates $\sigma, \theta_0' \circ \delta$, where
$\theta_0$ is the action implemented by
$\tilde{\theta}_{|\Lambda_0}$ on $Y_0'$. In particular $\theta_0$ is
weak mixing (because $\sigma$ is). This implies that if we denote
$q=\chi_{Y_0'}$ then the finite dimensional subspace $L^\infty (\Bbb
Z/m\Bbb Z)q$ of $L^\infty Y q$, which is clearly invariant to
$\theta_0$, must be reduced to the scalars. Thus $\tau(q) \leq 1/m$,
in other words $t \leq 1$.

\hfill $\square$

\noindent {\it Proof of Corollary 5.9.} If $\mycal F(\Cal R_\Lambda)
\neq \{1\}$ then after amplifying $\Cal R_\Lambda$ by a sufficiently
large number we may assume $\Cal R_\Gamma \subset \Cal R^t_\Lambda$,
for some $t> 1$. But this contradicts $5.8$. \hfill $\square$

\vskip .05in

\noindent {\it Proof of Corollary 5.10}. By ([P2]) we have $\mycal
F(\Cal R)\supset s^\Bbb Z$, with equality whenever $\Gamma$ is of
the required form. If $\Cal R^t=\Cal R_\Lambda$ for some free
ergodic $\Lambda$-action then $(\Cal R_\Lambda)^{1/t} = \Cal R
\supset \Cal R_\Gamma$, where $\Cal R_\Gamma$ is the equivalence
relation given by the Bernoulli $\Gamma$-action $\Gamma
\curvearrowright X=X_0^\Gamma$. Thus, if we can show that any finite
$\Gamma$-invariant subequivalence relation $\Cal R'\subset \Cal R$
is included into $\Cal R_\Gamma$ then the rest of the statement
follows from Corollary 5.9. With the notations in Sec. 1.4, let
$M=L(\Cal R)$ be the von Neumann algebra of $\Cal R$ with $u_g \in
M, g\in \Gamma$, the canonical unitaries implementing $\Gamma
\curvearrowright X$. Then $\Cal R'$ $\Gamma$-invariant is equivalent
to $B=L(\Cal R')$ being Ad$u_g$-invariant, $g\in \Gamma$. Note that
$A=L^\infty X$ is contained in $L(\Cal R')$ and $\Cal R'$ is finite
iff there exist finitely many partial isometries $\{v_i\}_i\in B$
normalizing $A$ and such that $E_A(v_i^*v_j)=0$ for $i\neq j$,
$B=\Sigma_i v_i A$. If we denote $N=L(\Cal R_\Gamma) = A \vee
\{u_g\}_g$, this also implies that $\Cal H_0=\Sigma_i v_i N$ is a
$N$-bimodule. If $\Cal R' \not\subset \Cal R_\Gamma$, then there
exists $i$ such that $v_i \not\in N$. It follows that there exists a
projection $0\neq q\in Av_i^*v_i$ such that $v_iq=wu_h\in B$, for
some $h\in \Gamma$ and $w\in L(\Cal R_0)$ normalizing $A$ with
$E_A(w)=0$. But the action $\sigma_g(x)=u_gxu_g^*$, $x\in L(\Cal
R_0)$, $g\in \Gamma$, is mixing (see e.g. 2.4.3 in [P1], or 1.6.2 in
[P2]). Thus, $\forall \delta>0$, $\exists g\in \Gamma$ such that
$\|E_N(\sigma_g(w^*)v_j)\|_2 \leq \delta$, $\forall j$. Since
$w=v_iqu_g^*\in \Cal H_0$ and $\delta$ is arbitrary, this shows that
$w=0$, contradiction. \hfill $\square$

\vskip .05in

\noindent {\it Proof of Proposition 5.11}. For each $g\in \Gamma$,
$h\in \Lambda$, we denote $X^g_h=\{t\in X \mid w(t,g)=h\}$. By the
definition of $w$ it follows that for each $g\in \Gamma$,
$\{X^g_h\}_h$ gives a (a.e.) partition of $X$ into measurable
subsets. Also, if $\sigma$ is free, then both in the case $\Delta$
is an embedding or local OE, the sets $\{X_h^g\}_g$ are disjoint.
Indeed, this is because in both cases the map $\Delta$ is 1 to 1 on
($\mu$-almost) each orbit of $\Cal R_\sigma$, so that if $t\in X^g_h
\cap X^{g'}_h$ then
$\Delta(g^{-1}t)=h^{-1}\Delta(t)=\Delta({g'}^{-1}t)$, implying that
$g^{-1}t={g'}^{-1}t$, thus $g=g'$ by the freeness of $\sigma$.
Denote $p^g_h=\chi_{X^g_h}\in L^\infty X$.

We begin by re-writing the equivalence of the $\Lambda$-valued
cocycles $w, \delta$ for the $\Gamma$-action $\sigma$ in von Neumann
algebra framework, by using the ``dictionary'' in Section 2. Thus,
we view $\Lambda$ as the (discrete) group of canonical unitaries
$\{v^0_h \mid h \in \Lambda \}$ of the group von Neumann algebra
$L\Lambda$ associated with $\Lambda$ and $w_g=w(\cdot, g) \in
\Lambda^X \subset L^\infty X \overline{\otimes} L\Lambda$.
Consequently, $w_g = \Sigma_h p_h^g \otimes v^0_h$ and $v=\Sigma_h
q_h \otimes v^0_h$. Note that $\Sigma_h q_h =p$ and $\Sigma_h
p^g_h=p, \forall g\in \Gamma$.

If we still denote by $\sigma$ the action of $\Gamma$ on $L^\infty X
\overline{\otimes} L\Gamma$ given by  $\sigma_g \otimes id, g\in
\Gamma$, then the cocycle relation for $w$ becomes
$w_{g_1}\sigma_{g_1}(w_{g_2})=w_{g_1g_2}$, $\forall g_1,g_2\in
\Gamma$. Also, the equivalence $w \sim \delta$ given by $v$ becomes

$$
w_g \sigma_g(v)=v (1 \otimes v^0_{\delta(g)}), g\in \Gamma,
$$
which in turn translates into the identities:
$$
\Sigma_h p_h^g \sigma_g(q_{h^{-1}k})= q_{k\delta(g)^{-1}}, \forall k
\in \Lambda, g\in \Gamma. \tag 5.11.1
$$

We first prove that the set of conditions $(5.11.1)$ is equivalent
to the following intertwining relation, viewed in $P$ when $\Delta$
is an embedding, respectively in $M^n$ when $\Delta$ is a local OE:

$$
u_g b = b v_{\delta(g)}, g \in \Gamma. \tag 5.11.2
$$
By using the definitions and the fact that $p^g_hv_h = p^g_h u_g$,
(cf. 1.4), we get in $P$ (resp. in $M^n$) the equalities
$$
u_g b = (\Sigma_{h_1} p^g_{h_1} v_{h_1})(\Sigma_{h_2} q_{h_2}
v_{h_2}) \tag 5.11.3
$$
$$
= \Sigma_{h_1,h_2} p_{h_1}^g \sigma_{g}(q_{h_2}) v_{h_1h_2}
=\Sigma_k (\Sigma_h p_h^g \sigma_g (q_{h^{-1}k})) v_k
$$
and also
$$
b v_{\delta(g)} = \Sigma_k q_{k\delta(g)^{-1}} v_k. \tag 5.11.4
$$
Thus, $(5.11.2)$ holds true if and only if we have the identities
$$
\Sigma_h p_h^g \sigma_g (q_{h^{-1}k})=q_{k\delta(g)^{-1}}, \forall
k\in \Lambda, g\in \Gamma, \tag 5.11.5
$$
which are exactly the same as conditions $(5.11.1)$. We have thus
shown that $b$ intertwines the representations $\{u_g\}_g,
\{v_{\delta(g)}\}_g$ of $\Gamma$, i.e. $u_g b = b v_{\delta(g)}$,
$\forall g\in \Gamma$. In particular, $bb^*$ commutes with
$\{u_g\}_g$ and $b^*b$ commutes with $\{v_{\delta(g)}\}_g$. To see
from this that $ker(\delta)$ is finite note that for each $g \in
ker(\delta) \vartriangleleft \Gamma$ we have $u_gb=bv_e=b$, hence
$u_gs=s$, where $s$ is the support projection of $bb^*$. This shows
that the left regular representation of $\ker(\delta)$ contains the
trivial representation of $ker(\delta)$, implying that $ker(\delta)$
is finite (see e.g. [D]).

Note now that if $q'\in B$ (resp. $q'\in A^n$) is a projection fixed
by $\{v_{\delta(g)} \mid g\in \Gamma \}$ then $bq'=\Sigma_h
q_h\theta_h(q') v_h$ is still an intertwiner between $\{u_g\}_g$ and
$\{v_{\delta(g)}\}_g$. Thus $[bq'b^*, u_g]=0$, $[q'b^*bq',
v_{\delta(g)}]=0$, $\forall g\in \Gamma$. Since $u_g, v_h$ normalize
$B$ (resp. $A^n$), this implies that for all $g\in \Gamma$ we have
$$
[\Sigma_h q_h\theta_h(q'), u_g]=[E_B(bb^*), u_g]=0, $$ $$ [\Sigma_h
\theta_{h^{-1}}(q_h)q', v_{\delta(g)}]=[E_B(b^*b),v_{\delta(g)}]=0,
$$
in the embedding case, and

$$
[\Sigma_h q_h\theta_h(q'), u_g]=[E_{A^n}(bb^*), u_g]=0, $$ $$
[\Sigma_h \theta_{h^{-1}}(q_h)q',
v_{\delta(g)}]=[E_{A^n}(b^*b),v_{\delta(g)}]=0,
$$
in the local OE case, where we have still denoted (and will do so
hereafter) by $\theta_h$ the automorphism Ad$(v_h)=\theta'_h$ on
$A^n$, for simplicity.

Applying this to non-zero spectral projections $q'$ of $\Sigma_h
\theta_{h^{-1}}(q_h)$, by the ergodicity of $\sigma_g=\text{\rm
Ad}(u_g), g\in \Gamma$, on $A$ it follows that $\Sigma_h
q_h\theta_h(q')$ is equal to $p=\chi_X$. Thus $q_h \leq
{\theta}_h(q')$, or equivalently $\theta_{h^{-1}}(q_h) \leq q'$, for
all $h \in \Lambda$. Summing up over $h\in \Lambda$ it follows that
the support projection $q$ of $\Sigma_h \theta_{h^{-1}}(q_h)$ is
majorized by $q'$. This shows $b^*b=\Sigma_h
\theta_{h^{-1}}(q_h)=cq$ for some scalar $c$, necessarily an integer
with $c=\tau(p)/\tau(q)$. Thus, $c^{-1/2}b$ is a partial isometry
with right support $q\in A^n$. On the other hand, taking in the
above $q'$ to be a projection in $Bq$ (resp. in $A^nq$) fixed by
$\{v_{\delta(g)} \mid g\in \Gamma\}$, this also shows that $q'=q$,
i.e. $\{v_{\delta(g)}\}_g$ must act ergodically on $Bq$ (resp. $A^n
q$).

Let us first finalize the proof of case $\Delta$ is a local OE, i.e.
when $M^n \supset P$, $A^n \supset B$. In this case, $\sigma$ weak
mixing on $L^\infty X=A=A^np$ implies that $\{u_g\}_g' \cap M$ is
contained in the center of the group von Neumann algebra
$L\Gamma=\{u_g\}_g'' \subset M$, so in particular the projection
$e=c^{-1}bb^*$ lies in $\Cal Z(L\Gamma)$. Thus, $u_ge \mapsto
v_{\delta(g)}q$ is an equivalence between the left regular
representations of $\Gamma/ker(\delta)$ and
$\delta(\Gamma)=\Lambda_0$ (with respect to the identification of
these two groups implemented by $\delta$), spatially implemented by
$c^{-1/2}b$.

Noticing that all the coefficients of $bb^*=\Sigma_h (\Sigma_l
q_{hl}\theta_h(q_l))v_h$ are projections, it follows that if we
denote $K=\ker(\delta)$ then $c=|K|$, $b=\Sigma_{k\in K} u_k$ and
$e=|K|^{-1}\Sigma_{k\in K} u_k$. Moreover, since $qM^nq=A^n q \vee
\{v_h q \mid h \in \Lambda_0\}$, we have $qv_hq=0, \forall h\in
\Lambda \setminus \Lambda_0$ implying that $\Lambda_0 =
\delta(\Gamma)$ has index $[\Lambda:\Lambda_0]=n = \tau(q)^{-1}$.

Assume now that $\Delta$ is an embedding, i.e. $M^n \subset P$,
$A^n=B$. Let $z=E_M(bb^*)\in M$ and $s$ the support projection of
$z$. Since $bb^*$ commutes with $\{u_g\}_g$, it follows that $z\in
L\Gamma '\cap M=\Cal Z(L\Gamma)$ (because $\sigma$ is weak mixing on
$A$). Since $\{v_h\mid h\in \Lambda_0\}$ is in its standard
representation on $q(L^2(L\Lambda_0))\simeq \ell^2\Lambda_0$ and the
isomorphism $u_g s \mapsto v_{\delta(g)}q$ is spatially implemented,
it follows that $z$ must be a multiple of $s$ which in turn must
equal a minimal projection in the algebra generated by $u_g, g\in K
=ker(\delta)$, lying in the center of this algebra. On the other
hand, since $bb^*=\Sigma_h (\Sigma_l q_{hl}\theta_h(q_l))v_h$ and
$\{\theta_h(q_{l})\}_l$ are mutually orthogonal, it follows that the
Fourier coefficients $a_h$ of $bb^*=\Sigma_h a_h v_h$ are
projections. Using that $\{p^g_h\}_h$ are mutually orthogonal, this
implies the Fourier coefficients of $z=E_M(bb^*)$ in $\{u_g\}_g$ are
projections as well. This shows that we must have $z=\Sigma_{k\in K}
u_k$.

On the other hand, by the form of $b$ we have $Ab \supset bA$ and
$b^*Ab=Aq$, thus $bb^*Abb^* \subset Abb^*$ implying that the support
projection $e$ of $bb^*$ implements a $\tau_\mu$-preserving
conditional expectation $E_0$ of $A$ onto a subalgebra $A_0\subset
A$ by $E_0(a)=\tau_M(e)^{-1}E_A(eae), a\in A$. Thus, if we denote by
$\Cal K\subset \Cal R^{Y_0}_\theta$ the sub-equivalence relation
implemented by the partial isometries $a_hv_h$ then $\Cal K$
contains $\Cal R_{\sigma(K)}$, $A_0=L^\infty (X/\Cal K)$, all orbits
of $\Cal K$ have $m$ points ($\mu$-a.e.), $m=\tau_M(e)^{-1}$ and by
the definitions $\alpha$ intertwines the restriction of the
$\Gamma/K$-action implemented by $\sigma$ on $A_0=A^{\Cal K}$ with
the $\Lambda_0=\delta(\Gamma)$ action implemented by
$\theta_{|\Lambda_0}$ on $A^nq=Bq$. The last part of the statement
is now trivial. \hfill $\square$

\heading 6. Final remarks
\endheading

\noindent {\bf 6.1}. Proposition 5.11 shows that all OE rigidity
results 5.6-5.8 hold in fact true without the assumption that the
source group $\Gamma$ has no finite normal subgroups, provided we
replace ``conjugate'' by ``virtual conjugate'' (in the sense of
[Fu1,2]) in the conclusions. Moreover, if we are only interested in
getting a ``virtual conjugacy'' conclusion, then by 5.11 we do not
need the assumption {\it any $\Gamma$-invariant finite
subequivalence relation of $\Cal R_\theta^t$ is contained in $\Cal
R_\sigma$} in the statement of Superrigidity of Embeddings 5.8. On
the other hand, note that in order to get a ``conjugacy'' conclusion
in 5.8, the above condition is unavoidable, as shown by the
following example: Let $\Lambda=\Gamma \times K_0$, with $K_0$ a
finite group, and $\Lambda \curvearrowright^{\theta'} X$ a free
action. Let $\Lambda \curvearrowright^\theta (X/K_0 \times K_0)$ be
the product of the action $\Gamma \curvearrowright X/K_0$
implemented by $\theta'$ and the action of $K_0$ on itself by left
translation. It is trivial to see that $\theta,\theta'$ are orbit
equivalent. If $\sigma=\theta'_{|\Gamma}$ then the inclusion $\Gamma
\subset \Gamma \times K_0=\Lambda$ implements an emebedding of
equivalence relation $\Cal R_\sigma\subset \Cal R_{\theta'}$. Thus
$\Cal R_\sigma \subset \Cal R_\theta$. However, if $\Gamma$ is
w-rigid, $K_0\neq 1$ and $\sigma$ is say a Bernoulli
$\Gamma$-action, then it is easy to see that $\sigma$ cannot be
conjugate to the restriction of $\theta$ to a subgroup $\Lambda_0$
of $\Gamma \times K_0$.

\vskip .05in \noindent {\bf 6.2}. One can obtain a cocycle
superrigidity result similar to 5.2 in which all assumptions are the
same except that ``s-malleable'' is replaced by the weaker
assumption ``malleable'' (as defined in the first part of 4.3), but
where the target group is in the smaller class of Polish groups of
{\it finite type I}, $\mycal U_{I,fin}$, i.e. Polish groups that can
be embedded as closed subgroups of the unitary groups of finite type
I von Neumann algebras (e.g. compact groups, or residually finite
discrete groups). Moreover, the same is true if $\sigma$ is merely
{\it sub-malleable}, i.e. if it can be extended to a malleable
action with the relative weak mixing condition satisfied. The proof
of this fact is simpler than that of Theorem 5.2, closely mimicking
proofs in ([P1], [PSa]).

\vskip .05in \noindent {\bf 6.3}. Alex Furman noticed in ([Fu4])
that Gaussian actions (see e.g. [CCJJV] for the definition) are also
s-malleable. The proof is the same as the proof of s-malleability of
their non-commutative version, the Bogoliubov actions in (1.6.3 of
[P2]). Gaussian actions are weak mixing whenever they come from
orthogonal representations with no finite dimensional invariant
subspaces. It would be interesting to produce more examples of
malleable actions, or more generally of quotients of malleable
actions. It seems to us that the action obtained by taking a group
embedding $\Gamma \subset \Cal G$ into a compact (Lie) group $\Cal
G$ with $\overline{\Gamma}=\Cal G$ and letting $\Gamma$ act on $\Cal
G$ with its Haar measure by (left) translation (page 1085 in [Fu2])
should be a (quotient of) malleable action. But being compact, these
actions are not weak mixing.

\vskip .05in \noindent {\bf 6.4}. As already pointed out in ([P4]),
if $\sigma$ is a Bernoulli action of an infinite
Kazhdan group $\Gamma$ on
the probability space $(X,\mu)$ and $K\subset \text{\rm Aut}(X,\mu)$
is a finite abelian group commuting with $\sigma(\Gamma)$ then the
``quotient action'' $\sigma^K$ of $\Gamma$ on $X/K$ implemented by
$\sigma$ is free mixing but not sub malleable, in the sense of $6.2$
above. Indeed, by ([P4]) we have H$^1(\sigma^K)=\text{\rm
Hom}(\Gamma \times K, \Bbb T)$ while if $\sigma^K$ would be sub
malleable then by ([P4], [PSa]) we would have
H$^1(\sigma^K)=\text{\rm Hom}(\Gamma, \Bbb T)$. This is a
contradiction, since for $K$ non-trivial $\text{\rm Hom}(\Gamma
\times K, \Bbb T)=\text{\rm Hom}(\Gamma, \Bbb T) \times \hat{K}$ has
more elements than the (finite) group $\text{\rm Hom}(\Gamma, \Bbb
T)$.

Another class of actions that are not sub malleable are the {\it
rigid} actions, as defined in (5.10.1 of [P5]). The prototype such
action is obtained by taking a rigid inclusion of groups of the form
$\Bbb Z^n \subset \Bbb Z^n \rtimes \Gamma$ (e.g. $n=2$ and $\Gamma
\subset SL(2,\Bbb Z)$ non-amenable), and letting $\sigma$ be the
action of $\Gamma$ on $\Bbb T^n=\hat{\Bbb Z^n}$ induced by the
action of $\Gamma$ on $\Bbb Z^n$. In fact, it it can be easily shown
using same arguments as in (Sec. 6 of [P3]) that $\Gamma
\curvearrowright \Bbb T^n$ cannot be realized as a quotient of a
malleable action.

\vskip .05in \noindent {\bf 6.5}. It would be interesting to have an
abstract characterization of the Polish groups of finite type. Note
in this respect that any such group $\Cal V$ is isomorphic to a
close subgroup of the unitary group of a (separable) Hilbert space
and admits a complete metric which is both left and right invariant.
Are these conditions sufficient to insure that $\Cal V\in \mycal
U_{fin}$? As far as ``classic'' groups are concerned, it is
interesting to recall an old result of Kadison and Singer ([KaSi];
see also [D1]), improving on an earlier result of von Neumann and
Segal ([vNS]), showing that if a connected locally compact group $G$
can be faithfully represented into the unitary group of a finite von
Neumann algebra then  $G= K \times H$ where $K$ is connected compact
and $H$ is a vector group (I am grateful to Dick Kadison and Raja
Varadarajan for pointing out to me this result). Thus, by
Proposition 3.5, any Polish group that contains the homeomorphic
image of a connected locally compact group $G$ which is not of this
form is not in the class $\mycal U_{fin}$ either. Thus, the class
$\Cal U_{fin}$ of ``target groups'' in our cocycle superrigidity
results 5.2/5.3 is essentially disjoint from linear algebraic
groups, which are the target groups in Zimmer's cocycle
superrigidity ([Z1,2]). As pointed out by Alex Furman, the fact that
$G={\text{\rm GL}}(n)$ and other linear algebraic groups are not in
$\mycal U_{fin}$ follows also as a combination of Zimmer's result
and the ``hereditary principle'' 3.5.

\vskip .05in \noindent {\bf 6.6}. In its general form, the cocycle
superrigidity result 5.2 shows that any cocycle for a diagonal
product $\Gamma$-action $\sigma \times \rho$, with $\Gamma$ baring
some mild {\it rigidity} (of Kazhdan-type) and $\sigma$ some {\it
malleability} property, is ``absorbed'' by the $\rho$ action. The
validity of such ``principle'' to situations more general than the
ones covered by 5.2 may lead to further applications, notably to the
calculation of the Feldman-Moore higher cohomology groups H$^n(\Cal
R_\sigma), n\geq 2,$ of such $\Gamma$-actions $\sigma$ (see [FM] for
the definition of H$^n$). In this respect, we expect that group
actions satisfying condition 5.6.1 should satisfy H$^2(\Cal
R_\sigma)=\text{\rm H}^2(\Gamma)$ ($\Bbb T$-valued 2-cocycles for
$\Gamma$). In this same vein, it would be of great interest to
extend the class of target groups $\Cal V$ covered by the Cocycle
Superrigidity 5.2, 5.3 from $\mycal U_{fin}$ to a larger class. It
seems reasonable to believe that such larger class may include all
Polish subgroups of the unitary group $\Cal U(\Cal H)$ on Hilbert
space, in particular all separable locally compact groups. In other
words, the following may hold true: {\it If $\Gamma$ has an infinite
w-normal rigid subgroup $H$ and $\Gamma \curvearrowright X$ is
s-malleable and weak mixing on $H$ then any $\Cal V$-valued cocycle
for $\sigma$, with $\Cal V \subset \Cal U(\Cal H)$ a closed
subgroup, is cohomologous to a group morphism of $\Gamma$ into $\Cal
V$}.

\vskip .05in \noindent {\bf 6.7}. It would be extremely interesting
to characterize the class $\Cal C\Cal S$ of groups $\Gamma$ for
which any s-malleable weak mixing $\Gamma$-action $\Gamma
\curvearrowright^\sigma X$ (e.g. $\Gamma \curvearrowright
[0,1]^\Gamma$) is $\Cal U_{fin}$-cocycle superrigid, i.e. for which
any measurable $\Cal V$-valued cocycle over $\sigma$ is cohomologous
to a group morphism of $\Gamma$ into $\Cal V$, for any $\Cal V \in
\Cal U_{fin}$. The class $\Cal C\Cal S$ should be much larger than
the class of groups considered in this paper. It cannot, of course,
contain the free groups (see e.g. Sec. 3 in [P4]), but may contain
all non-amenable groups of the form $\Gamma = \Gamma_0 \times
\Gamma_1$ with both $\Gamma_i$ infinite, thus relating with the
rigidity results of Monod-Shalom ([MoSh1,2]) and Hjorth-Kechris
([HjKe]).

\vskip .05in \noindent {\bf 6.8}. Let $\{\Gamma_i\}_    {i\in I}$ be a
countable (at most) family of groups such that each $\Gamma_i$ has a
wq-normal rigid subgroup $H_i\subset \Gamma_i$. Denote $\Gamma=*_i
\Gamma_i$ the free product of the groups $\Gamma_i$. As pointed out
in ([P4]), if $|I|\geq 2$ then $\Gamma$ does not have any wq-normal
rigid subgroup. Nevertheless, the following version of cocycle
superrigidity does hold true for these free product groups: Let
$\sigma$ be an action of $\Gamma$ on the standard probability space
$(X,\mu)$ and $\Cal V \in \mycal U_{fin}$. Assume that for each
$i\in I$, $\sigma_{|\Gamma_i}$ is s-malleable and either mixing or
weak mixing with $H_i$ w-normal in $\Gamma_i$. If $w: X \times
\Gamma \rightarrow \Cal V$ is a measurable cocycle for $\sigma$ then
there exist group morphisms $\delta_i : \Gamma_i \rightarrow \Cal V$
and measurable maps $v_i:X \rightarrow \Cal V$ such that $w_g =
v_i^*\gamma_i(g)\sigma_g(v_i)$, $\forall g\in \Gamma_i$, $\forall
i\in I$. When combined with Proposition 5.11, this can be used to
obtain an OE rigidity result of Bass-Serre type, in the spirit of
results in ([IPeP]).

\vskip .05in \noindent {\bf 6.9}. If one restricts the proofs in
Sec. 2-4, leading to the proof of the Cocycle Superrigidity 5.2/5.3,
to the case the groups act on commutative algebras, then the
arguments can be easily translated into measure theoretical terms.
We opted for a von Neumann algebra presentation because the ideas
behind the proofs are so much ``von Neumann algebra'' in spirit and
because of the non-commutative generalizations it allows (e.g. 5.5).
Another reason was that we needed this setting for Proposition 5.11,
which we proved using von Neumann algebra analysis. In this respect,
we mention that shortly after the initial version of this paper was
circulated, Stefaan Vaes in ([V]) and then Alex Furman ([Fu4]) were
able to give alternative, genuine measure theory proof to the ``OE
case'' of 5.11, using 3.3 in [Fu2]. Their expository notes also
contain measure theoretical presentations of our proof of 5.2.

On the other hand, note that we only used the Cocycle Superrigidity
5.2 towards establishing rigidity results for embeddings and local
OE of actions. It should be possible to derive from 5.2 rigidity
results for arbitrary morphisms (as defined in 1.4.2) between
equivalence relations implemented by free m.p. actions. Such
morphisms still give rise to cocycles, which by 5.2 can be untwisted
whenever the source action satisfies 5.6.1. But it is not clear how
to interpret the ``untwister'' as some kind of ``generalized
conjugator'' in this generality. In fact, one needs to first
understand what rigidity conclusion one seeks to derive for
arbitrary morphisms between these equivalence relations.
The following type of morphisms, generalizing both local OE and embeddings,
could be easier to
analize: If $\Gamma\curvearrowright X, \Lambda \curvearrowright Y$
are free ergodic m.p. actions, then a morphism
$\Delta:(X,\mu)\rightarrow (Y,\nu)$ between $\Cal R_\sigma, \Cal R_\theta$
is a {\it local embedding} if it is 1 to 1 (but
not necessarily onto) on $\mu$-almost every $\Gamma$-orbit. The
``ideal'' superrigidity statement for a local
embedding $\Delta$ should generalize both 5.7, 5.8. It should
(roughly) show that there exist $\Lambda_0\subset \Lambda$ and
$\alpha\in [\Gamma]$, $\beta \in [\Lambda]$ such that the
$\Delta$-pull back of $\beta\Lambda_0 \beta^{-1}$ is conjugate to a
quotient of $\alpha \Gamma \alpha^{-1}$.

\vskip .05in \noindent {\bf 6.10}. Several applications of the
Cocycle Superrigidity 5.2/5.3 have been obtained since the initial
circulation of this paper in Dec. 2005. This includes Theorem 5.7
(solving 5.7.2$^\circ$ in [P2]) and Corollary 5.10 in the present
version of the paper, and a result by Stefaan Vaes and the author
showing that if $\Gamma \curvearrowright X$ satisfies 5.6.1,
$K\curvearrowright X$ is a compact action commuting with it and so
that $\Gamma \curvearrowright X/K$ is still free, then $\Gamma
\curvearrowright X/K$ is OE superrigid ([PV]). Also, Alex Furman
used 5.2 to show that if $\Gamma,\Lambda$ are lattices in a higher
rank semisimple Lie group $\Cal G$ then the action $\Gamma
\curvearrowright \Cal G/\Lambda$ cannot be realized as a quotient of
a Bernoulli $\Gamma$-action, more generally of a s-malleable weak
mixing action ([Fu4]). This solves a problem posed in ([P3], see
comments after 7.7) and shows that the case $\Gamma$ higher rank
lattice of the OE Superrigidity 5.6 is already covered by the OE
Superrigidity in ([Fu2]). Finally, Simon Thomas in ([T])
used the Cocycle
Superrigidity 5.2 to answer some open problems in descriptive set
theory (Borel equivalence relations), showing for instance that the
universal countable Borel equivalence relation $E_\infty$
cannot be implemented by a free action of a countable group
(see [JaKeL] for definitions).

\head  References\endhead

\item{[A1]} S. Adams: {\it Indecomposability of treed equivalence relations}
Israel J. Math., {\bf 64} (1988), 362-380.

\item{[A2]} S. Adams: {\it Indecomposability of equivalence
relations generated by word hyperbolic groups}, Topology {\bf 33}
(1994), 785-798.

\item{[B]} M. Burger: {\it Kazhdan constants for} $SL(3,\Bbb Z)$,
J. reine angew. Math., {\bf 413} (1991), 36-67.

\item{[CCJJV]} Cherix, Cowling, Jolissaint, Julg, Valette: ``Groups
with Haagerup property'', \newline
Birkh$\ddot{\text{\rm a}}$user Verlag, Basel Berlin
Boston, 2000.

\item{[C1]} A. Connes: {\it Une classification des facteurs de type} III,
Ann. \'Ec. Norm. Sup {\bf 1973}, 133-252.

\item{[C2]} A. Connes: {\it Classification of injective factors},
Ann. of Math., {\bf104} (1976), 73-115.

\item{[C3]} A. Connes: {\it A type II$_1$ factor with countable
fundamental group}, J. Operator Theory {\bf 4} (1980), 151-153.

\item{[CFW]} A. Connes, J. Feldman, B. Weiss: {\it An amenable equivalence
relation is generated by a single transformation}, Ergodic Theory
Dynamical Systems {\bf 1} (1981), 431-450.

\item{[CJ]} A. Connes, V.F.R. Jones: {\it A} II$_1$ {\it factor
with two non-conjugate Cartan subalgebras}, Bull. Amer. Math. Soc.
{\bf 6} (1982), 211-212.

\item{[Co]} Y. de Cornulier: {\it Relative Kazhdan property},
Ann. Sci. Ecole Norm. Sup. {\bf 39} (2006), 301-333.

\item{[D1]} J. Dixmier: ``Les C$^*$-Alg\'ebres et leurs repr\'esentations'',
Gauthier-Villars, Paris 1969.

\item{[D2]} J. Dixmier: {\it Sous anneaux ab\'eliens maximaux dans
les facteurs de type fini}, Ann. of Math. {\bf 59} (1954), 279-286.

\item{[Dy1]} H. Dye: {\it On groups of measure preserving
transformations} I, Amer. J. Math, {\bf 81} (1959), 119-159.

\item{[Dy2]} H. Dye: {\it On groups of measure preserving
transformations}, II, Amer. J. Math, {\bf 85} (1963), 551-576.

\item{[FMo]} J. Feldman, C.C. Moore: {\it Ergodic equivalence
relations, cohomology, and von Neumann algebras I, II}, Trans. Amer.
Math. Soc. {\bf 234} (1977), 289-324, 325-359.

\item{[Fe]} T. Fernos: {\it Relative Property} (T) {\it and linear groups},
Ann. Inst. Fourier, {\bf 56} (2006), 1767-1804.

\item{[FiHi]} D. Fisher, T. Hitchman: {\it Cocycle superrigidity and
harmonic maps with infinite dimensional targets}, math.DG/0511666,
preprint 2005.

\item{[Fu1]} A. Furman: {\it Gromov's measure
equivalence and rigidity of higher rank lattices}, Ann. of Math.
{\bf 150} (1999), 1059-1081.

\item{[Fu2]} A. Furman: {\it Orbit equivalence
rigidity}, Ann. of Math. {\bf 150} (1999), 1083-1108.

\item{[Fu3]} A. Furman: {\it Outer automorphism groups of some
ergodic equivalence relations}, Comment. Math. Helv. {\bf 80}
(2005), 157 - 196.

\item{[Fu4]} A. Furman: {\it On Popa's Cocycle Superrigidity
Theorem}, to appear.

\item{[F]} H. Furstenberg: {\it
Ergodic behavior of diagonal measures and a theorem of Szemeredi on
arithmetic progressions}, J. d'Analyse Math. {\bf 31} (1977)
204-256.

\item{[G1]} D. Gaboriau: {\it Cout des r\'elations d'\'equivalence
et des groupes}, Invent. Math. {\bf 139} (2000), 41-98.

\item{[G2]} D. Gaboriau: {\it Invariants $\ell^2$ de r\'elations
d'\'equivalence et de groupes},  Publ. Math. I.H.\'E.S. {\bf 95}
(2002), 93-150.

\item{[Ge]} S.L. Gefter: {\it On cohomologies of ergodic actions
of a T-group on homogeneous spaces of a compact Lie group}
(Russian), in ``Operators in functional spaces and questions of
function theory'', Collect. Sci. Works, Kiev, 1987, pp 77-83.

\item{[GeGo]} S.L. Gefter, V.Y. Golodets: {\it Fundamental
groups for ergodic actions and actions with unit fundamental
groups}, Publ RIMS {\bf 6} (1988), 821-847.

\item{[HVa]} P. de la Harpe, A. Valette: ``La propri\'et\'e T
de Kazhdan pour les groupes localement compacts'', Ast\'erisque {\bf
175}, Soc. Math. de France (1989).

\item{[HjKe]} G. Hjorth, A. Kechris: ``Rigidity theorems
for actions of product groups and countable Borel equivalence
relations'', Memoirs of AMS {\bf 177}, No. 833, 2005.

\item{[IPeP]} A. Ioana, J. Peterson, S. Popa: {\it Amalgamated free
products of w-rigid factors and calculation of their symmetry
groups}, math.OA/0505589, to appear in Acta Math.

\item{[JaKeL]} S. Jackson, A. Kechris, G. Hjorth: {\it Countable
Borel equivalence relations}, J. Math. Logic {\bf 1} (2002), 1-80.

\item{[Jo]} P. Jolissaint: {\it On Property} (T)
{\it for pairs of topological groups}, l'Ens. Math. {\bf 51} (2005),
31-45.

\item{[J]} V.F.R. Jones : {\it Index for subfactors}, Invent.
Math. {\bf 72} (1983), 1-25.

\item{[KaSi]} R.V. Kadison, I.M. Singer: {\it Some remarks
on representations of connected groups}, Proc. Nat. Acad. Sci. {\bf
38} (1952), 419-423.

\item{[K]} D. Kazhdan: {\it Connection of the dual space of a
group with the structure of its closed subgroups}, Funct. Anal. and
its Appl. {\bf1} (1967), 63-65.

\item{[Ma]} G. Margulis: {\it Finitely-additive invariant measures
on Euclidian spaces}, Ergodic. Th. and Dynam. Sys. {\bf 2} (1982),
383-396.

\item{[MoS1]} N. Monod, Y. Shalom:
{\it Cocycle superrigidity and bounded cohomology for negatively
curved spaces}, J. Diff. Geom. {\bf 67} (2004), 395-455.

\item{[MoS2]} N. Monod, Y. Shalom:
{\it Orbit equivalence rigidity and bounded cohomology}, Ann. of
Math. {\bf 164} (2006), 825-878.

\item{[MvN1]} F. Murray, J. von Neumann:
{\it On rings of operators}, Ann. Math. {\bf 37} (1936), 116-229.

\item{[MvN2]} F. Murray, J. von Neumann:
{\it Rings of operators IV}, Ann. Math. {\bf 44} (1943), 716-808.

\item{[vNS]} J. von Neumann, I.E. Segal: {\it A theorem on unitary
representations of semisimple Lie groups}, Ann. of Math.
{\bf 52} (1950), 509-516.

\item{[OW]} D. Ornstein, B. Weiss: {\it Ergodic theory of
amenable group actions I. The Rohlin Lemma} Bull. A.M.S. (1) {\bf 2}
(1980), 161-164.

\item{[P1]} S. Popa: {\it Some rigidity results for
non-commutative Bernoulli shifts}, J. Fnal. Analysis {\bf 230}
(2006), 273-328.

\item{[P2]} S. Popa: {\it Strong Rigidity of} II$_1$ {\it Factors
Arising from Malleable Actions of $w$-Rigid Groups} I, Invent. Math.
{\bf 165} (2006), 369-408 (math.OA/0305306).

\item{[P3]} S. Popa: {\it Strong Rigidity of} II$_1$ {\it Factors
Arising from Malleable Actions of $w$-Rigid Groups} II, Invent.
Math. {\bf 165} (2006), 409-452 (math.OA/0407137).

\item{[P4]} S. Popa: {\it Some computations of $1$-cohomology groups
and construction of non orbit equivalent actions}, J. Inst. Math.
Jussieu {\bf 5} (2006), 309-332 (math.OA/0407199).

\item{[P5]} S. Popa: {\it On a class of type II$_1$ factors with
Betti numbers invariants}, Ann. of Math. {\bf 163} (2006), 809-889
(math.OA/0209310).

\item{[P6]} S. Popa: ``Classification of subfactors and of their
endomorphisms'', CBMS Lecture Notes, {\bf 86}, Amer. Math. Soc.
1995.

\item{[P7]} S. Popa: {\it Markov traces on  universal Jones
algebras and subfactors of finite index}, Invent. Math. {\bf 111}
(1993), 375-405.

\item{[P8]} S. Popa: {\it Correspondences}, INCREST preprint
1986 (unpublished),\newline www.math.ucla.edu/~popa/preprints.html

\item{[PSa]} S. Popa, R. Sasyk:
{\it On the cohomology of Bernoulli actions}, Erg. Theory Dyn. Sys.
{\bf 27} (2007), 241-251 (math.OA/0310211).

\item{[PV]} S. Popa, S. Vaes: {\it Strong rigidity of generalized Bernoulli
actions and computations of their symmetry groups}, math.OA/0605456.

\item{[Sh]} Y. Shalom: {\it Measurable group theory}, in
``Proceedings of the 2004 European Congress of Mathematics''.

\item{[Si]} I.M. Singer: {\it Automorphisms of finite factors},
Amer. J. Math. {\bf 77} (1955), 117-133.

\item{[T]} S. Thomas: {\it Popa Superrigidity and Countable Borel
Equivalence Relations}, to appear.

\item{[V]} S. Vaes: {\it Rigidity results for Bernoulli
actions and their von Neumann algebras} (after Sorin Popa)
S\'eminaire Bourbaki, expos\'e {\bf 961}, to appear in Ast\'erisque
(math.OA/ \newline 0605456).

\item{[Va]} A. Valette: {\it Group pairs with relative property} (T)
{\it from arithmetic lattices}, Geom. Dedicata {\bf 112} (2005),
183-196.

\item{[Z1]} R. Zimmer: {\it Strong rigidity for ergodic actions of
seimisimple Lie groups}, Ann. of Math. {\bf 112} (1980), 511-529.

\item{[Z2]} R. Zimmer: ``Ergodic Theory and Semisimple Groups'',
Birkhauser, Boston, 1984.

\item{[Z3]} R. Zimmer: {\it Extensions of ergodic group actions},
Ill. J. Math. {\bf 20} (1976), 373-409.

\item{[Z4]} R. Zimmer: {\it Superrigidity, Ratner's theorem and the
fundamental group}, Israel J. Math. {\bf 74} (1991), 199-207.

\enddocument